\documentclass[11pt]{amsart}
\title[Integral points on symmetric varieties] {Integral points on
  symmetric varieties and Satake compatifications} 
\author{Alexander Gorodnik, Hee Oh and Nimish
  Shah} \address{Mathematics 253-37\\ Caltech\\Pasadena, CA 91106}
\email{gorodnik@caltech.edu} \address{Mathematics 253-37\\
  Caltech\\Pasadena, CA 91106} 
\curraddr{Math Department\\151 Thayer St.\\Brown
  University\\Providence, RI 02912} \email{heeoh@math.brown.edu}

\address{School of Mathematics, TIFR\\ 1 Homi Bhabha Road\\ Mumbai,
             400005, India} \email{nimish@math.tifr.res.in}

\thanks{The first and the second authors partially supported by NSF
  0400631 and NSF 0333397, 0629322 respectively}

\usepackage{epsfig} \usepackage{graphicx} \usepackage{amsmath,amscd,amssymb}
\theoremstyle{plain} 

\newtheorem{Thm}[equation]{Theorem}
\newtheorem{thm}[equation]{Theorem}
\newtheorem{que}[equation]{Question}

\newtheorem{Cor}[equation]{Corollary}
\newtheorem{cor}[equation]{Corollary}
\newtheorem{Prop}[equation]{Proposition}
\newtheorem{Lem}[equation]{Lemma} 
\newtheorem{lem}[equation]{Lemma}
\newtheorem{rem}[equation]{Remark}
\newtheorem{Def}[equation]{Definition}
\newtheorem{Ex}[equation]{Example} \numberwithin{equation}{section}

\newcommand{\q}{\mathbb{Q}}
\newcommand{\Q}{\mathbb{Q}}
\newcommand{\e}{\varepsilon}
\newcommand{\z}{\mathbb{Z}}
\newcommand{\Z}{\z}

\newcommand{\N}{\mathbb{N}}
\newcommand{\K}{\mathbb{K}}
\renewcommand{\c}{\mathbb{C}}
\newcommand{\C}{\mathbb{C}}
\newcommand{\br}{\mathbb{R}}
\newcommand{\R}{\br}

\newcommand{\cC}{\mathcal C}

\newcommand{\cO}{\mathcal O}

\newcommand{\cW}{\mathcal W}
\newcommand{\G}{\Gamma}

\newcommand{\Cal}{\mathcal}

\newcommand{\cl}[1]{\overline{#1}}
\newcommand{\inv}{^{-1}}

\providecommand{\abs}[1]{\left\lvert#1\right\rvert}
\providecommand{\norm}[1]{\left\lVert#1\right\rVert}
\providecommand{\inpr}[2]{\left\langle #1,\, #2\right\rangle}

\newcommand{\GL}{\operatorname{GL}}
\newcommand{\vol}{\operatorname{Vol}}

\newcommand{\SL}{\operatorname{SL}}
\newcommand{\SO}{\operatorname{SO}}
\newcommand{\la}[1]{\mathfrak{\lowercase{#1}}}

\newcommand{\Stab}{\operatorname{Stab}}
\newcommand{\spn}{\operatorname{span}}

\newcommand{\Ad}{\operatorname{Ad}}

\newcommand{\ad}{\operatorname{ad}}
\newcommand{\tr}{\operatorname{tr}}
\newcommand{\diag}{\operatorname{diag}}
\newcommand{\supp}{\operatorname{supp}}
\newcommand{\Pic}{\operatorname{Pic}}

\newcommand{\ignore} [1] {}
\newcommand{\del}{\partial}
\newcommand{\trn}{\,{}^t\!}

\newif\ifdraft\drafttrue
\begin{document}

\begin{abstract}
  Let $V$ be an affine symmetric variety defined over $\Q$.  We
  compute the asymptotic distribution of the angular components of the
  integral points in $V$.  This distribution is described by a family
  of invariant measures concentrated on the Satake boundary of $V$.
  In the course of the proof, we describe the structure of the Satake
  compactifications for general affine symmetric varieties and compute
  the asymptotic of the volumes of norm balls.
\end{abstract}

\maketitle

\tableofcontents

\section{Introduction}

\label{sec:intro}

Let $V=\{x\in\R^n:\, f_1(x)=\cdots=f_s(x)=0\}$ with
$f_i\in\Z[x_1,\ldots,x_n]$ be an affine variety. It is a fundamental
problem of Diophantine geometry to understand the set of integral
points $V(\Z)$ in $V$. In particular, when the number of integral
points is infinite, one may ask

\begin{que}\label{que1}
  Given a norm $\norm{\cdot}$ on $\R^n$, determine the asymptotic of
  \begin{align*}
    N_T(V):=\#\{x\in V(\Z) : \norm{x}< T\}
  \end{align*}
  as $T\to\infty$.
\end{que}

We are interested in a more refined question:

\begin{que}\label{que2}
  For a radial cone $\mathcal{C}$ in $\R^n$ centered at the origin,
  determine the asymptotic of
  \begin{align*}
    N_T(V,\mathcal{C}) :=\#\{x\in V(\Z)\cap\mathcal{C}:\norm{x}< T\}
  \end{align*}
  as $T\to\infty$.
\end{que}

One can also state an analogous question in terms of convergence of
measures.  We define a probability measure $\mu_T$ on the unit sphere
$S^{n-1}$ in $\R^n$:
\begin{align*}
  \mu_T:=\frac{1}{N_T(V)}\sum_{x\in V(\Z):0<\norm{x}<T}
  \delta_{\pi(x)},
\end{align*}
where $\pi:\R^n\setminus\{0\}\to S^{n-1}$ denotes the radial
projection, and $\delta_z$ denotes the Dirac measure at $z$.  As
$T\to\infty$, these measures characterize the asymptotic distribution
of the angular components of points of $V(\Z)$.

\begin{que}\label{que3}
  Determine the weak$^*$ limits of the measures $\mu_T$ as
  $T\to\infty$.
\end{que}

Recall that a sequence of measures $\{\mu_i\}$ on $S^{n-1}$ converges
to $\mu$ in weak$^*$ topology if $\int_{S^{n-1}}\phi\, d\mu_i
\to\int_{S^{n-1}}\phi\,d\mu$ for every $\phi\in C(S^{n-1})$.

In this paper we give a complete solution to Questions~\ref{que2}
and~\ref{que3} when $V$ is an affine symmetric variety.  In this case,
Question~\ref{que1} was answered by Duke, Rudnick, Sarnak \cite{DRS}
and Eskin, McMullen \cite{EM}, though explicit asymptotics in terms of
$T$ were not computed in general.  Later Eskin, Mozes and Shah developed an approach
using the ergodic theory on homogeneous spaces, based on the work
of Dani, Margulis (\cite{DM1}, \cite{DM2}) and Ratner \cite{Ra} on unipotent flows.

We note that the method of \cite{EM} shows that
$N_T(V, \mathcal C)$ is asymptotic to the volume of
the set $\mathcal C_T:=\{x\in V\cap \mathcal C: \|x\|<T\}$,
provided the family $\{\mathcal C_T\}$ satisfies the so-called
well roundedness property. However it is a highly nontrivial task
to check whether the sets $\mathcal C_T$ are well rounded, and this
is precisely where the main technical difficulties of this paper lies
(see section \ref{ip} for more discussion on this point).

 The above questions \ref{que2} and \ref{que3} are motivated
by the conjectures of Manin \cite{bm,fmt}, Peyre \cite{p}, and
Batyrev, Tschinkel \cite{bt}, which describe the distribution of
rational points on projective Fano varieties (see
Remark~\ref{r:tshinkel} and \cite{Ts}). Recently, Chambert-Loir and Tschinkel
proposed an analogous conjecture for integral points.  We expect that
our results will support this conjecture (see Section~\ref{s:tsch}).

We illustrate our main results by the following example of a quadratic
surface.  We refer to Section~\ref{sec:ex} for further examples.

\begin{Ex}
  \label{quad}
  \rm Fix $n\ge 4$ and $k\in \z\setminus\{0\}$. Let $Q$ be an integral
  non-degenerate indefinite quadratic form in $n$ variables such that
  $Q(x)=k$ has at least one integral solution. Let $\Omega$ be a Borel
  subset of $S^{n-1}$ such that the interior of $\Omega$ intersects
  $\{Q=0\}$, and the boundary of $\Omega$ has measure zero with
  respect to a smooth measure on $\{Q=0\}\cap S^{n-1}$.  Then setting
  $\mathcal{C}=\R^+\cdot\Omega$, we have
  \begin{align*}
    N_T(\{Q=k\},\mathcal{C}) \sim_{T\to\infty} \vol(\{x\in \mathcal
    C:\, Q(x)=k ,\, \norm{x} < T\}) \sim_{T \to \infty} d_{\mathcal C}
    \cdot T^{n-2},
  \end{align*}
  where the volume is computed with respect to a suitably normalized
  $\SO(Q)$-invariant measure on $\{Q=k\}$ and $d_{\mathcal C}>0$ is a
  computable constant.
  
  \begin{center}
    \includegraphics[clip=true,width=6cm,height=6cm]{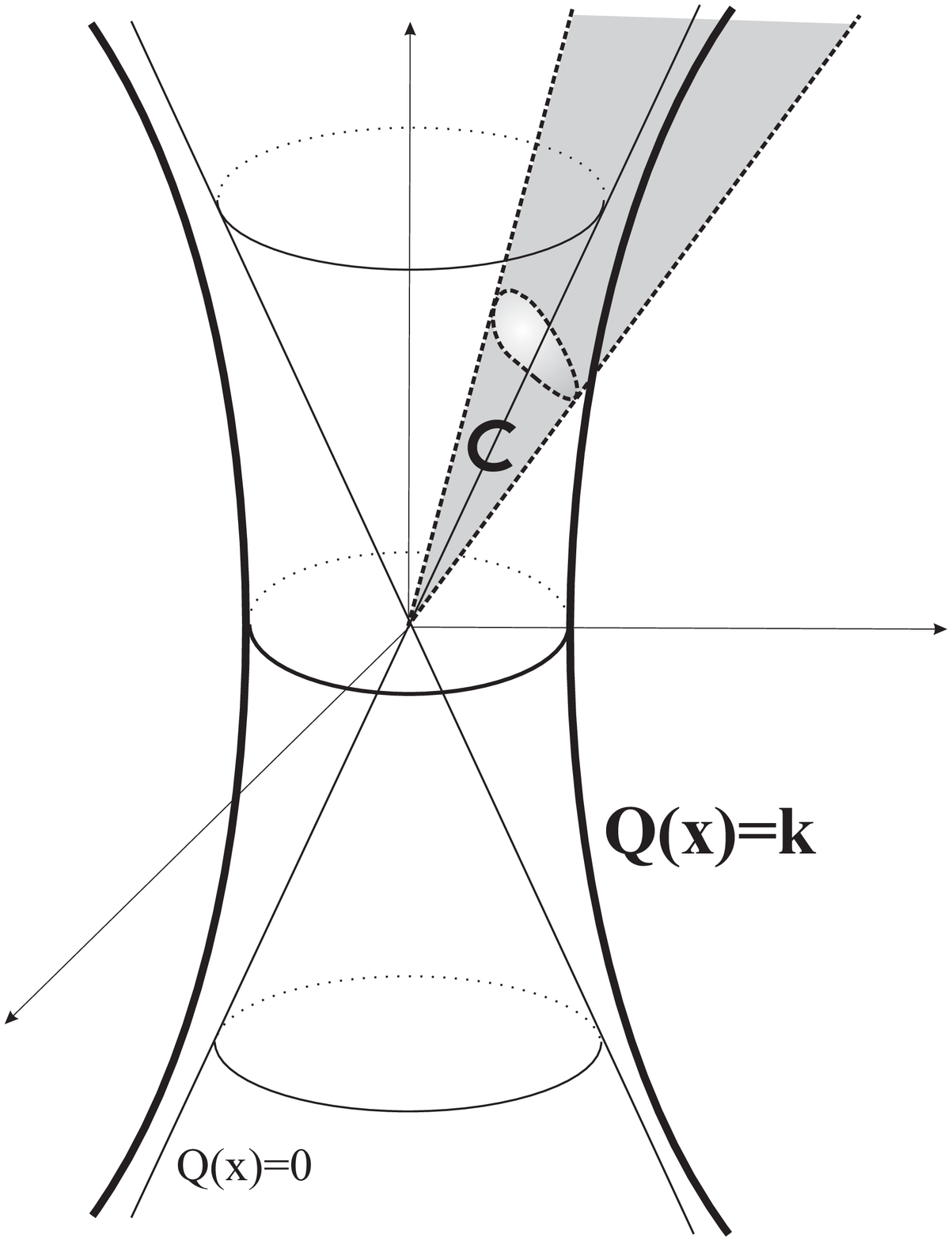}
  \end{center}
    
  We denote by $\nu$ the unique $\SO(Q)$-invariant measure on
  $\{x\in\R^n\setminus\{0\}:\,Q(x)=0\}$ normalized so that
  \begin{align*}
    \nu(\{x\in\R^n\setminus\{0\}: Q(x)=0,\norm{x}<1\})=1.
  \end{align*}
  Define the measure $\mu$ on $S^{n-1}$ by
  \begin{align*}
    \mu:=\pi_*(\nu|_{B_1})
  \end{align*}
  where $B_1=\{x\in\R^n: \norm{x}<1\}$.  Then $\mu_T\to\mu$ as
  $T\to\infty$.
\end{Ex}

\subsection{Main results}

Let $\mathbf G$ be a connected $\q$-simple algebraic $\q$-group
isotropic over $\R$ with a given $\R$-irreducible $\q$-rational
representation $\iota:{\mathbf G}\to \GL({\mathbf W})$.  Suppose that
there exists $v_0\in {\mathbf W}(\q)$ such that ${\mathbf G}v_0$ is
Zariski closed and that $V:={\mathbf G}(\br)^\circ v_0$ is an affine
symmetric (real) variety, that is, the stabilizer $\mathbf H$ of $v_0$
in $\mathbf G$ is the set of fixed point of an involution $\sigma$ of
$\mathbf G$. Examples of affine symmetric varieties are provided by
Proposition~\ref{p:exist}.

\subsubsection*{Notation} $G={\mathbf G}(\br)^\circ$, $W={\mathbf
  W}(\R)$, and $H={\mathbf H}(\br)\cap G$.  Then $V\cong G/H$.
\emph{We assume that $V(\Z)\ne \emptyset$ and that $\mathbf H$ has no
  nontrivial $\q$-characters.}

We fix a basis, say ${\mathcal B}$, of $W$, and define
\begin{equation} \label{eq:S(W)} S(W)=\big\{\sum_{{\mathbf
      e}\in{\mathcal B}} x_{\mathbf e}{\mathbf e}\in W: \sum_{{\mathbf
      e}\in {\mathcal B}} x_{\mathbf e}^2= 1\big\}.
\end{equation}
Let $\pi:W\setminus\{0\}\to S(W)$ denote the radial projection.

We define the {\it Satake boundary} $V^\infty$ of $V$:
\begin{align*}
  V^\infty:=\{\lim \pi(v):\, \text{$v\in V$, $v\to \infty$}\} =
  \overline{\pi(V)}-\pi(V).
\end{align*}
For example, when $V=\{x\in\R^n:\,Q(x)=k\}$ is a quadratic surface,
\begin{align*}
  V^\infty=\{x\in S^{n-1} :\, Q(x)=0\}.
\end{align*}
The map $\pi$ embeds $V$ homeomorphically into $\overline{\pi(V)}$ as
an open dense subset (see Proposition~\ref{p:homeo}), and we call
$\overline{\pi(V)}$ the {\it Satake compactification\/} of $V$. For a
Riemannian symmetric space, this compactification was introduced by
Satake in \cite{sat}.  We note that Satake \cite{sat} considered only
the special case when $\iota$ is a representation on the space of
bi-linear forms.

The action of $G$ on $W$ induces, via $\pi$, a $G$-action on $S(W)$.
In Section~\ref{sec:satake} we will prove that $V^\infty$ is a union
of finitely many $G$-orbits, which are locally closed.

Given a measure $\mu$ on $V^\infty$ which is a linear combination of
smooth positive measures on some $G$-orbits, we say that $\mu$ is {\em
  concentrated on\/} the union of these $G$-orbits.

Let $\mathcal{C}$ be a Borel cone in $W$ centered at the origin.  In
order to have meaningful results about the sets
$V(\Z)\cap\mathcal{C}$, it is necessary to assume that the
intersection $V\cap\mathcal{C}$ is large in a suitable sense.  As we
will see below, the ``size'' of $V(\Z)\cap\mathcal{C}$ depends quite
sensitively on the set of $G$-orbits in $V^\infty$ that $\mathcal{C}$
intersects.
 
A Borel cone $\mathcal{C}\subset W$ centered at $0$ is called {\em
  admissible\/} if the closure of $\cC$ has nonempty intersection with
$V^\infty$, and the boundary of $\cC$ is of zero measure with respect
to the smooth measure class on each of the finitely many $G$-orbits on
$V^\infty$.  A Borel cone $\cC\subset W$ centered at $0$ is called
{\em generic\/} if the closure of $\cC$ and the interior of $\cC$
intersect the same collection of $G$-orbits in $V^\infty$.
  
The following theorem gives a natural generalization of
Example~\ref{quad}.  We fix any norm $\norm{\cdot}$ on $W$ and set
\begin{align*}B_T=\{v\in W: \norm{v}<T\}.\end{align*}

\begin{thm}\label{mt}
  For every admissible generic cone $\mathcal{C}\subset W$,
  \begin{align*}
    \#(V(\Z)\cap B_T\cap \mathcal{C})\sim_{T\to\infty} \vol(V\cap
    B_T\cap \mathcal{C})\sim_{T\to\infty} d_{\mathcal{C}}\cdot
    T^{a_\mathcal{C}}(\log T)^{b_\mathcal{C}-1},
  \end{align*}
  where the volume is computed with respect to a suitably normalized
  $G$-invariant measure on $V$, and $d_{\mathcal{C}}>0$,
  $a_\mathcal{C}\in\Q^+$, $b_\mathcal{C}\in\N$ are computable
  constants.
\end{thm}

Given a Euclidean norm $\norm{\cdot}$ on $W$ and $v_0\in V^\infty$,
the cone $\mathcal{C}=\{v:\, \norm{\pi(v)-v_0}<\e\}$ is admissible and
generic for all sufficiently small $\e>0$ (see
subsection~\ref{subsec:proof:cor:dioph}).  Hence, we get the following
application of Theorem~\ref{mt} to Diophantine approximation.

\begin{Cor}\label{cor:dioph}
  Let $\norm{\cdot}$ be a Euclidean norm on $W$ and $v_0\in V^\infty$.
  Then for any sufficiently small $\e>0$, there exist $c=c(v_0,\e)>0$,
  $a=a(v_0)\in \Q^+$, $b=b(v_0)\in\N$ such that
  \begin{align*}
    \#\{ v\in V(\z)\cap B_T: \norm{\pi(v)- v_0}<\e\} &\sim_{T\to
      \infty} \vol(\{v\in V \cap
    B_T: \norm{\pi(v)-v_0}<\e\})\\
    &\sim_{T\to \infty} c\cdot T^{a}(\log T)^{b-1}.
  \end{align*}
\end{Cor}

Now we describe the structure of the Satake boundary $V^\infty$ of
$V$.  Let $K$ be a maximal compact subgroup of $G$ compatible with $H$
and $\mathfrak{a}$ a Cartan subalgebra corresponding to the pair $K$
and $H$, so that the Cartan decomposition $G=K\exp(\mathfrak{a})H$
holds (\cite[Ch.~7]{sch}).  We fix a system of simple roots
$\Delta_\sigma$ of $\la{a}$ and denote by
$\mathfrak{a}^{+}\subset\la{a}$ the closed positive Weyl chamber. One
can choose a subset $\mathcal{W}$ of the normalizer of $\mathfrak{a}$
in $K$ such that we have a decomposition
\begin{align*}
  G=K\exp(\mathfrak{a}^{+})\mathcal{W}H
\end{align*}
where the $\la{a}^+$-component and the $\mathcal{W}$-component of each
element of $G$ are uniquely defined (see Section~\ref{sec:basic}).

For any subset $J\subset \Delta_\sigma$, we set
\begin{align*}
  V^\infty_{J} = \left\{\lim \pi(k\exp(a)wv_0):
    \begin{tabular}{l}
      $k\in K$, $a\in \mathfrak{a}^{+}$, $w\in\mathcal{W}$,\\
      $\alpha(a)\to\infty$ for $\alpha\in\Delta_\sigma \setminus J$,\\
      $\alpha(a)$ is bounded for $\alpha\in J$.
    \end{tabular}
  \right\}.
\end{align*}
Note that $V^\infty_{\Delta_\sigma}=\pi(V)$ and
\begin{equation}\label{eq:vd_decomp}
  V^\infty=\bigcup_{J\subsetneq \Delta_\sigma} V^\infty_{J}.
\end{equation}
Every set $V^\infty_J$ is a union of finitely many $G$-orbits (see
Theorem~\ref{th:sat}).

Denoting by $2\rho$ the sum (with multiplicities) of all positive
roots of $\la{a}$ and by $\lambda_\iota$ the highest weight of the
representation $\iota$ with respect to $\la{a}^+$ (see
Section~\ref{sec:rep}), we have decompositions
\begin{align} \label{eq:defi0} 2\rho:=\sum_{\alpha\in\Delta_\sigma
  }u_\alpha \alpha\quad\text{and}\quad \lambda_\iota :=\sum_{\alpha\in
    \Delta_\sigma} m_\alpha \alpha.
\end{align}
Note that $u_\alpha>0$, $m_\alpha>0$, and $u_\alpha,m_\alpha\in\Q$
(\cite[p.~85]{ov}). Define
\begin{align}
  \label{eq:defi}
  a_\iota & = \max \left\{\frac{u_\alpha}{m_\alpha}: \alpha\in
    \Delta_\sigma
  \right\}, \nonumber \\
  I_\iota &= \big\{\alpha\in \Delta_\sigma:\frac{u_\alpha}{m_\alpha}
  <a_\iota\big\},  \\
  b_\iota &=\#(\Delta_\sigma\setminus I_\iota) \; \geq 1. \nonumber
\end{align}

\begin{thm}\label{th:con}
  For an admissible cone $\mathcal{C}\subset W$, the limits
  \begin{align*}
    \lim_{T\to\infty} \frac{\#(V(\Z)\cap B_T
      \cap\mathcal{C})}{T^{a_\iota}(\log T)^{b_\iota-1}}
    \quad\text{and}\quad\lim_{T\to\infty} \frac{\vol(V\cap B_T\cap
      \mathcal{C})}{T^{a_\iota}(\log T)^{b_\iota-1}}
  \end{align*}
  exist and are equal, where the volume is computed with respect to a
  suitably normalized $G$-invariant measure on $V$.  Moreover, if
  $\mathcal{C}^\circ \cap V^\infty_{I_\iota}\ne\emptyset$, then the
  limits are strictly positive.
\end{thm}

We also extend the result about convergence of measures in
Example~\ref{quad} to general affine symmetric spaces.  There is an
explicitly given $G$-invariant measure $\nu_\iota$ on $\R^+\cdot
V^\infty_{I_\iota}$, normalized such that $\nu_\iota(B_1)=1$. The
measure $\nu_\iota$ is homogeneous of degree $a_\iota$, and we have a
decomposition
\begin{align*}
  d\nu_\iota(t\cdot\theta) =t^{a_\iota-1}\, dt d\theta, \quad
  t\in\R^+,\; \theta \in V^\infty_{I_\iota},
\end{align*}
where $dt$ is a Lebesgue measure on $\R^+$ and $d\theta$ is a smooth
measure on $V^\infty_{I_\iota}$.  We define the probability measure
$\mu_\iota$ on $V_{I_\iota}^\infty$ by
$\mu_\iota=\pi_*(\nu_\iota|_{B_1})$ or equivalently,
\begin{align*}
  d\mu_\iota(\theta)=\norm{\theta}^{-a_\iota}d\theta,\quad \theta \in
  V^\infty_{I_\iota}.
\end{align*}
Note that the norm $\norm{\cdot}$ need not be constant on $S(W)$ for
our fixed choice of the sphere $S(W)$.  Later on, we give an explicit
formula for $\mu_\iota$ (see (\ref{eq:mu_iota})).

\begin{thm}\label{th:m}
  As $T\to\infty$, we have
  \begin{align*} \mu_T\to \mu_\iota.\end{align*}
\end{thm}

Theorem~\ref{th:m} also holds for representations $\iota$ which are
not irreducible (see Remark~\ref{rem:rep}).

In the case of the group variety (i.e., when $G=L\times L$ and
$H=\{(l,l):\, l\in L\}$), these equidistribution results
(Theorem~\ref{th:con} and Theorem~\ref{th:m}) were first proved by
Maucourant \cite{mau}, although not in terms of the Satake boundary.

\newcommand{\Ht}{\text{H}}

\begin{rem}\label{r:tshinkel} 
  \rm It is interesting to compare Theorem~\ref{th:m} with a result in
  \cite{gmo} (see also \cite{stt}), which describes distribution of
  rational points of bounded height. Let $\mathbf G$ be a connected
  adjoint absolutely simple algebraic $\Q$-group and $\iota:{\mathbf
    G}\to \GL_N$ an absolutely irreducible representation defined over
  $\Q$.  We denote by $\bar\iota$ the corresponding map from ${\mathbf
    G}$ to the projective space ${\mathbb P}^{N^2-1}$ and by
  $\Ht=\Ht_\infty\cdot\prod_{p\mbox{\tiny :prime}} \Ht_p$ a height
  function on ${\mathbb P}^{N^2-1}(\Q)$ where $\Ht_\infty$ is a norm
  on $\R^{N^2}$.  Let $G={\mathbf G}(\R)^\circ$, $G(\Q) =G\cap
  {\mathbf G}(\Q)$ and $G(\Z)=G\cap {\mathbf G}(\Z)$. As explained in
  Section~\ref{sec:group}, $V=\iota(G)$ is an affine symmetric
  variety, and we have a decomposition
  \begin{align*}
    \overline{\bar\iota(G)} = \bigcup_{I\subset \Delta_\sigma}
    V_I^\infty.
  \end{align*}
  It follows from \cite{gmo} that for every $\phi\in
  C(\overline{\bar\iota(G)})$,
  \begin{align*}
    \lim_{T\to\infty}\frac{1}{\#\{\gamma\in G(\Q):
      \Ht(\bar\iota(\gamma))<T\}}\sum_{\gamma\in G(\Q):
      \Ht(\bar\iota(\gamma))<T} \phi(\bar\iota(\gamma))=
    \int_{V^\infty_{\Delta_\sigma}} \phi(\omega)
    \frac{d\omega}{\Ht_\infty(\omega)^{a_\iota'}}
  \end{align*}
  where $d\omega$ is a $G$-invariant measure on
  $V^\infty_{\Delta_\sigma}=\bar\iota(G)$ and $a_\iota'>a_\iota$.  On
  the other hand, it follows from Theorem~\ref{th:m} that for every
  $\phi\in C(\overline{\bar\iota(G)})$,
  \begin{align*}
    \lim_{T\to\infty} \frac{1}{\#\{\gamma\in G(\Z):
      \Ht_\infty(\bar\iota(\gamma))<T\}} \sum_{\gamma\in G(\Z):
      \Ht_\infty(\bar\iota(\gamma))<T}
    \phi(\bar\iota(\gamma))=\int_{V^\infty_{I_\iota}} \phi(\theta)
    \frac{d\theta}{\Ht_\infty(\theta)^{a_\iota}}.
  \end{align*}
\end{rem}

\vspace{0.5cm}

Note that for affine symmetric varieties of higher rank, $\mu_\iota=
\lim_{T\to\infty}\mu_T$ is concentrated on $V_{I_\iota}^\infty$, which
might have empty interior in $V^\infty$ (see Section~\ref{sec:ex} for
examples). In particular, Theorem~\ref{th:m} does not imply
Theorem~\ref{mt}. To prove Theorem~\ref{mt}, we need to investigate
the accumulation of integral points on all components of $V^\infty$.

\begin{Def} \rm
  \label{d:lambda_conn}
  A subset $I\subset \Delta_\sigma$ is called {\it
    $\lambda_\iota$-connected\/} if the Dynkin diagram of
  $\{\lambda_\iota\}\cup I$ is connected. In other words, if
  $I\cup\{\lambda_\iota\}=S_1\cup S_2$, $S_i\neq\emptyset$, then
  $S_1\not\perp S_2$ with respect to the identification of
  $\la{a}^\ast$ with $\la{a}$ via the Killing form.
\end{Def}

We show (see Theorem~\ref{th:sat}) that
\begin{align*}
  V^\infty =\bigsqcup_{\text{$\lambda_\iota$-connected $I\subsetneq
      \Delta_\sigma$}} V^\infty_I,
\end{align*}
and for $\lambda_\iota$-connected $I,J\subset\Delta_\sigma$,
\begin{align}
  \label{eq:bdry}
  I\subsetneq J \Longleftrightarrow V_I^\infty\subset \del
  (V_J^\infty).
\end{align}

For $\lambda_\iota$-connected $I\subsetneq\Delta_\sigma$, we set
\begin{align}
  \nonumber a_\iota(I) & :=\max \left\{\frac{u_\alpha}{m_\alpha}:
    \alpha\in
    \Delta_\sigma \setminus I \right\}>0,\\
  \label{eq:Iiota}
  I_\iota(I) & :=I\cup \big\{\alpha\in \Delta_\sigma \setminus I:\,
  \frac{u_\alpha}{m_\alpha}
  <a_\iota(I)\big\}\subsetneq \Delta_\sigma,\\
  \nonumber b_\iota(I) &:=\#(\Delta_\sigma \setminus I_\iota(I)) \;
  \geq 1.
\end{align}

Note that $I_\iota(\emptyset)=I_\iota$. We will show in
Proposition~\ref{p:I_iota} that $I_\iota(I)$ is
$\lambda_\iota$-connected. We consider the lexicographical order on
the set of pairs $(a,b)\in\R\times \R$.  Note that for
$\lambda_\iota$-connected subsets $I$ and $J$ of $\Delta_\sigma$,
\begin{align*}
  I\subset J \implies I_\iota(I)\subset I_\iota(J) \implies
  (a_\iota(I),b_\iota(I))\geq (a_\iota(J),b_\iota(J)).
\end{align*}

For $\Omega\subset W$ with $\cl{\pi(\Omega\setminus\{0\})}\cap
V^\infty\neq\emptyset$, we define
\begin{align}
  \nonumber \Theta_\Omega &:=\{I\subsetneq \Delta_\sigma: \text{$I$ is
    $\lambda_\iota$-connected and
    $\cl{\pi(\Omega\setminus\{0\})}\cap V_I^\infty\ne \emptyset$}\},\\
  \label{eq:iota0}
  (a_\iota(\Omega),b_\iota(\Omega))&:=\max\{(a_\iota(I),b_\iota(I)):\,
  I\in \Theta_\Omega\},\\
  \nonumber \Theta_\iota(\Omega) &:=\{I_\iota(I):\, I\in\Theta_\Omega,
  (a_\iota(I),b_\iota(I)) = (a_\iota(\Omega),b_\iota(\Omega))\}\subset
  \Theta_\Omega.
\end{align}
Roughly speaking, we show that the asymptotic number of points in
$V(\Z)$, with norm less than $T$, whose images under $\pi$ accumulate
on $\Omega$, is of the order $T^{a_\iota(\Omega)}(\log
T)^{b_\iota(\Omega)-1}$.

It might happen that in Theorem~\ref{th:con}, both of the limits are
zero. This simply means that the normalization term $T^{a_\iota}(\log
T)^{b_\iota-1}$ is not suitable.  We prove a refined version of
Theorem~\ref{th:con}.

\begin{thm}\label{th:cone_o}
  For every admissible cone $\mathcal{C}\subset W$, the limits
  \begin{align*}
    \lim_{T\to\infty} \frac{\#(V(\Z)\cap B_T
      \cap\mathcal{C})}{T^{a_\iota(\mathcal{C})}(\log
      T)^{b_\iota(\mathcal{C})-1}}
    \quad\text{and}\quad\lim_{T\to\infty} \frac{\vol(V\cap B_T\cap
      \mathcal{C})}{T^{a_\iota(\mathcal{C})}(\log
      T)^{b_\iota(\mathcal{C})-1}}
  \end{align*}
  exist and are equal, where the volume is computed with respect to a
  suitably normalized $G$-invariant measure on $V$.  Moreover, if
  $\mathcal{C}^\circ \cap V^\infty_{I}\ne\emptyset$ for some
  $I\in\Theta_\iota(\mathcal{C})$, then the limits are strictly
  positive.
\end{thm}

Moreover if $\mathcal{C}$ is generic, then for any $I\in\Theta_{\cC}$,
$\mathcal{C}^\circ\cap V^\infty_I\ne \emptyset$; and by
\eqref{eq:bdry}, we get $\cC^\circ\cap
V^\infty_{I_\iota(I)}\neq\emptyset$.  Therefore
Theorem~\ref{th:cone_o} implies Theorem~\ref{mt}.

Let $\mathcal{D}_\iota=\{I_\iota(I): \lambda_\iota\text{-connected }
I\subsetneq \Delta_\sigma\}$.  To state the next result, we will
introduce a family of smooth measures $\mu_{I}$ on $V_I^\infty$ for
$I\in\mathcal{D}_\iota$. For each $I\in{\mathcal D}_\iota$, here is an
explicitly given $G$-invariant measure $\nu_{I}$ on $\R^+\cdot
V^\infty_{I}$, which is a finite union of $G$-orbits. Note that
$a_I=\frac{u_\alpha}{m_\alpha}$ is constant for $\alpha\in
\Delta_\sigma\setminus I$, and the measure $\nu_I$ is homogeneous of
degree $a_I$. We have a decomposition
\begin{align*}
  d\nu_I(t\cdot\theta) =t^{a_I-1}\, dt d\theta, \quad t\in\R^+,\;
  \theta \in V^\infty_{I},
\end{align*}
where $dt$ is a Lebesgue measure on $\R$ and $d\theta$ is a smooth
measure on $V^\infty_{I}$.  We define $\mu_I=\pi_*(\nu_I|_{B_1})$ or
equivalently,
\begin{align*}
  d\mu_I(\theta)=\norm{\theta}^{-a_I}d\theta,\quad \theta \in
  V^\infty_{I}.
\end{align*}
An explicit formula for $\mu_I$ is given in (\ref{eq:mu_alpha}).  For
$\mathcal{I}\subset \mathcal{D}_\iota$, we define
\begin{align*}
  \mu_\mathcal{I}=\sum_{I\in\mathcal{I}} \mu_I \quad \text{and} \quad
  \nu_\mathcal{I}=\sum_{I\in\mathcal{I}} \nu_I.
\end{align*}

\begin{thm}\label{th:m_o}
  For every $\phi\in C(S(W))$ with $\supp \phi\cap
  V^\infty\neq\emptyset$, we have
  \begin{align*}
    \lim_{T\to\infty} \frac{1}{T^{a_\iota(\phi)}(\log
      T)^{b_\iota(\phi)-1}} \sum_{x\in V(\Z):0<\norm{x}<T}
    \phi(\pi(x)) = c \int_{S(W)} \phi\; d\mu_{\Theta_\iota(\phi)},
  \end{align*}
  where $c>0$ depends only on $V(\Z)$, and
  $a_\iota(\phi)=a_{\iota}(\supp \phi)$, $b_\iota(\phi)=b_\iota(\supp
  \phi)$, $\Theta_\iota(\phi)= \Theta_\iota(\supp \phi)$.
\end{thm}

Note that if $\supp\phi\cap V^\infty\neq\emptyset$, then
$\phi(\pi(x))=0$ for all but finitely many $x\in V(\Z)$.

The measures $\mu_\mathcal{I}$ are analogues of the Patterson-Sullivan
measures, which are concentrated on the visual boundary of a
Riemannian symmetric space.

In order to prove Theorems~\ref{th:cone_o} and~\ref{th:m_o}, we
compare the asymptotic distribution of integral points to the
corresponding {\em continuous\/} asymptotic distribution, which is
given in the following theorem.

\begin{thm}\label{th:asympt0}
  For every $f\in C_c(W\setminus\{0\})$ with $\pi(\supp f)\cap
  V^\infty\neq\emptyset$,
  \begin{align*}
    \lim_{T\to\infty} \frac{1}{T^{a_{\iota}(f)}(\log
      T)^{b_{\iota}(f)-1}}\int_{G/H} f(gv_0/T)d\mu(g)=\int_W f\,
    d\nu_{\Theta_\iota(f)},
  \end{align*}
  where $a_\iota(f)=a_{\iota}(\pi(\supp f))$,
  $b_\iota(f)=b_\iota(\pi(\supp f))$, and $\Theta_\iota(f)=
  \Theta_\iota(\pi(\supp f))$.
\end{thm}

Note that the limit measure $\nu_{\Theta_\iota(f)}$ is homogeneous of
degree $a_\iota(f)$.

Also note that if $f\in C_c(W\setminus\{0\})$ and $\pi(\supp f)\cap
V^\infty=\emptyset$, then for all large $T$, $f(gv_0/T)=0$ for all
$g\in G$.

\begin{rem}\label{r:not_sim} 
  \rm The condition that the group ${\mathbf G}$ is $\Q$-simple can be
  relaxed.  In fact, it suffices to assume that every $\Q$-simple
  factor ${\mathbf G}_0$ of $\mathbf G$ is isotropic over $\R$ and
  ${\mathbf G}={\mathbf G}_0{\mathbf H}$.  A small modification in the
  proof is required only in Section~\ref{sec:equ} (see
  Remark~\ref{rem:mod}).
\end{rem}

\subsection{Ingredients of the proof}\label{ip}
A common strategy for estimating
the number of integral points in various domains
involves two steps:
\begin{enumerate}
\item establishing suitable regularity of domains and their volumes;
\item comparing the number of integral points with the volumes of the domains.
\end{enumerate}

The second step, discussed in Section \ref{sec:equ},
is essentially done using standard techniques developed in \cite{DRS, EM}
in view of the equidistribution theorem (Theorem \ref{mix}) available
in the symmetric setting.
Checking the first step for the domains defined by
the intersection of a cone with the norm balls contains the main technical difficulties of the paper
and requires, in particular, the analysis of the structure of
the Satake boundary (Section  \ref{sec:satake})
and asymptotic estimates for renormalized volumes with respect to the invariant measure (Sections
\ref{sec:invariant-measures}  and \ref{sec:vol}).
For instance, one of the reasons we are working in the setting of a symmetric space, rather
than of a general homogeneous space, is the lack of the structure theory
for Satake compactification needed to establish (1), since Theorem \ref{mix} is
available in a more general setting of homogeneous spaces as obtained in \cite{EMS}.

We mention that
checking the {\it well-roundedness} property of \cite{EM} for the domains
amounts to carrying out the first step.
Setting, for a radial cone $\mathcal C$ and $T>0$,
$$\mathcal C_T:=\{x\in V\cap \mathcal C: \norm{x}<T\},$$
 $\{\mathcal C_T: T\gg 1\}$
is in general {\it not} well-rounded.
We introduce the notion of an admissible generic cone
in terms of its intersection with the Satake boundary of $V$. Showing that
the family $\{\mathcal C_T: T\gg 1\}$ is well-rounded
for an admissible generic cone $\mathcal C$
is a consequence of 
 two main ingredients of the paper: 
\begin{itemize} 
 \item[(i)]  Tube Lemma (Coro~\ref{c:nbhd0}); in showing this lemma,
we needed to generalize Satake's theory in an affine symmetric setting.
\item[(ii)] The computation of the asymptotic limit
of the invariant measures on $V$ (Theorem \ref{th:asympt0}).
\end{itemize}

We emphasize that in order to compute the explicit volume asymptotic of $\mathcal C_T$,
we need only the part (ii). However in order to obtain
that $N_T(V, \mathcal C)\sim \text{Vol}(\mathcal C_T)$ (not to mention the explicit
asymptotic), we need to use both (i) and (ii).

%The solutions of Questions~\ref{que2} and~\ref{que3} use the structure
%of the Satake boundary of $V$.  
In the rest of this section, we explain the
generalization of Satake's result \cite{sat} on Riemannian symmetric
spaces, which summarizes properties of the decomposition
(\ref{eq:vd_decomp}) and Tube lemma.

\begin{thm}\label{th:sat}
  \begin{enumerate}
  \item[(a)] For every $J\subset \Delta_\sigma$, $V^\infty_J$ is a
    union of finitely many $G$-orbits.
  \item[(b)] For every $J\subset \Delta_\sigma$,
    \begin{align*}
      \overline{V^\infty_J}=\bigcup_{I\subset J} V^\infty_I.
    \end{align*}
  \item[(c)] For every $J\subset \Delta_\sigma$,
    $V^\infty_J=V^\infty_I$ where $I$ is the largest
    $\lambda_\iota$-connected subset of $J$.
  \item[(d)] For distinct $\lambda_\iota$-connected subsets
    $I,J\subset \Delta_\sigma$, $V^\infty_I\cap V^\infty_J=\emptyset$.

  \end{enumerate}
\end{thm}

Similar constructions can be also found in \cite{sat}, \cite{osh}, and
\cite{cp}, but they are not suitable for our purpose: Satake
\cite{sat} considered only representations on the space of bilinear
forms, the Oshima compactification \cite{osh} is defined abstractly,
and the de~Concini--Procesi compactification \cite{cp} is defined with
respect to the Zariski topology and applies only to a specific type of
representations which have regular highest weights.

While Theorem~\ref{th:sat} is not used in the proofs of
Theorem~\ref{th:con} and Theorem~\ref{th:m}, it is essential for the
proofs of Theorem~\ref{mt}, Theorem~\ref{th:cone_o}, and
Theorem~\ref{th:m_o}. The crucial observation is the following
corollary that describes neighborhoods of subsets in $V^\infty$.  For
$I\subset \Delta_\sigma$, we set
$$\la{a}_I=\ker I=
\{a\in\la{a}:\alpha(a)=0,\,\forall \alpha\in I\},$$ and
$\la{a}_I^+=\la{a}_I\cap\la{a}^+$, which is a face of the closed Weyl
chamber $\la{a}^+$.
 
\begin{cor}[Tube lemma]\label{c:nbhd0}
  Let $I$ be a $\lambda_\iota$-connected subset of $\Delta_\sigma$ and
  $\Omega$ a compact subset of $S(W)$ contained in $\bigcup_{J\supset
    I} V^\infty_J$. Then there exists a compact set $U\subset
  \mathfrak{a}^{+}$ such that
  \begin{align*}
    \Omega\cap \pi(V) &\subset \pi
    (K\exp(U+\mathfrak{a}_I^{+})\mathcal{W}v_0).
  \end{align*}
\end{cor}

Note that by Theorem~\ref{th:sat}(b), $\bigcup_{J\supset I}
V^\infty_J$ is the smallest open subset of $\overline{\pi(V)}=
V^\infty \cup \pi(V)$ which contains $V^\infty_I$ and is a union of
cells.  The following picture illustrates the structure of the Satake
boundary for an affine symmetric variety associated to an
$\R$-irreducible representation of $\SL_3$ (see
Section~\ref{sec:group}) with the highest weight $\lambda_\iota$ which
is not orthogonal to the simple roots $\alpha$ and $\beta$.  The
shaded regions correspond to neighborhoods of points in the components
of $V^\infty$.
\begin{center}
  \includegraphics[clip=true,height=5cm]{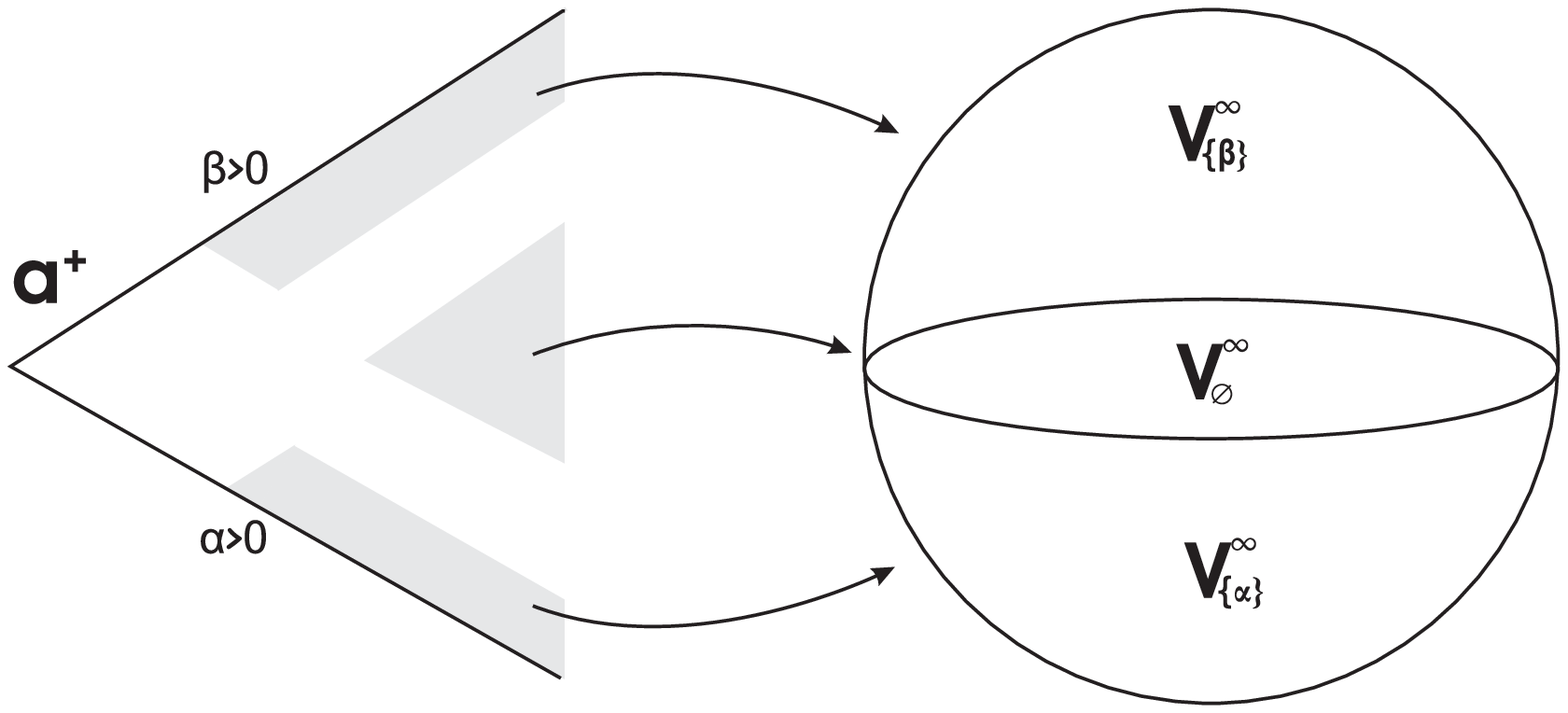}
\end{center}
On the other hand, if $\lambda_\iota\perp \alpha$, we get
\begin{center}
  \includegraphics[clip=true,height=5cm]{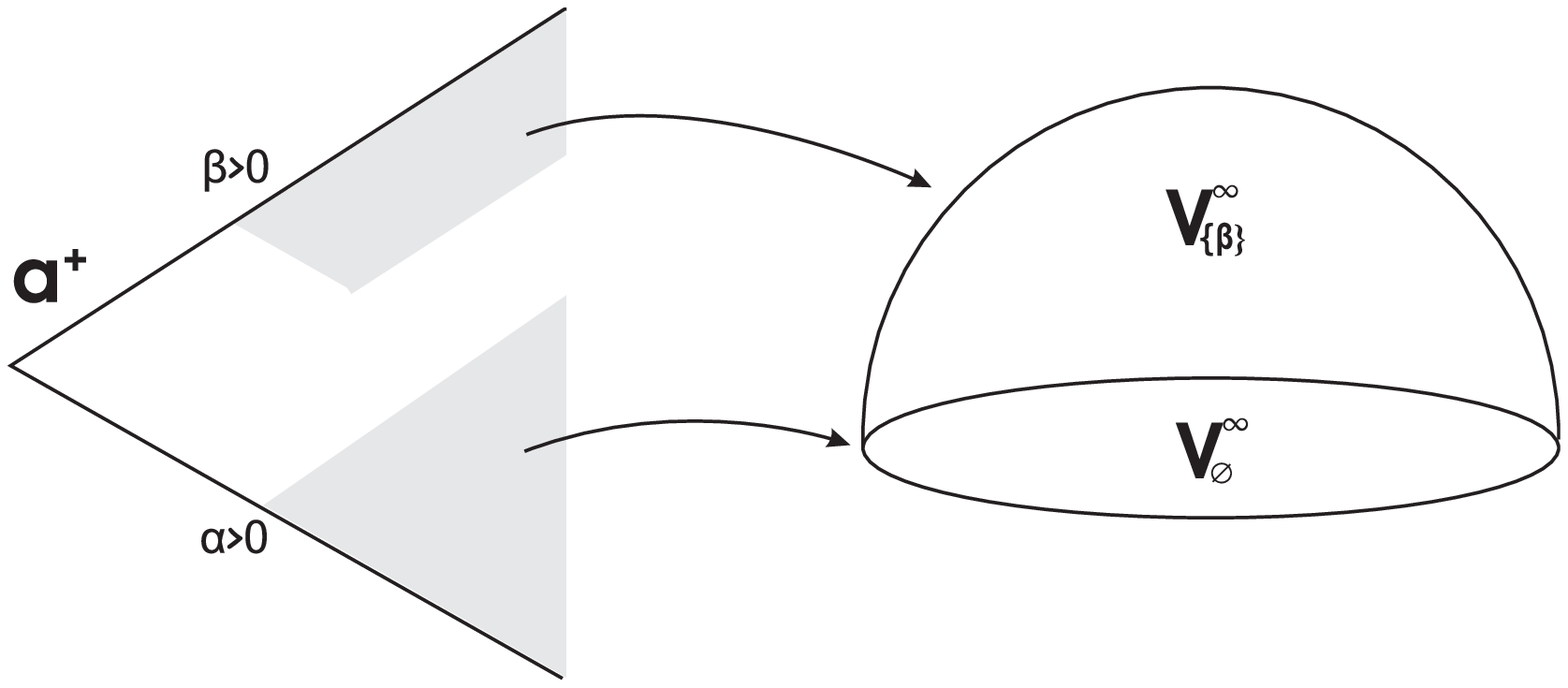}
\end{center}

As in the earlier works on counting integral points \cite{DRS,EM}, the
basic dynamical or ergodic theoretic ingredient in our proof is the
result on limiting distributions of translates of closed $H$-orbits in
$G/G(\Z)$ as established in \cite{EM} (see Theorem~\ref{mix}).

\subsection{Organization of the paper}
Section~\ref{sec:ex} contains examples.  In Section~\ref{sec:basic},
we review some basic properties of affine symmetric spaces and
representation.  The structure of the Satake boundary, including the
proofs of Theorem~\ref{th:sat} and Corollary~\ref{c:nbhd0}, is
discussed in Section~\ref{sec:satake}.  Explicit formula for the
measures $\nu_I$ are obtained in Section~\ref{sec:invariant-measures}.
Section~\ref{sec:vol} contains results on volume asymptotics
(Theorem~\ref{th:asympt0}), which are described via $G$-invariant
measures on the boundary. Finally, the main theorems are proved in
Section~\ref{sec:equ}.  In Section~\ref{s:tsch}, we state a version of
our main result using the language of arithmetic geometry.

\subsection{Acknowledgment}
The authors would like to thank Gopal Prasad for providing us some
important arguments used in the proof of Proposition~\ref{p:max_stab}.

\section{Examples} \label{sec:ex}

\subsection{Quadric} 
We start with an example of a rank-one symmetric space where the
structure of the Satake boundary is quite simple.  Let $Q$ be an
integral non-degenerate quadratic form on $\R^{p+q}$ of signature
$(p,q)$, $p,q\ge 1$, $p+q\ge 4$, and $k\in \N$. We are interested in
the distribution of integral points lying on the the quadratic surface
$V:=\{Q=k^2\}$.  To simplify notation, we assume that
\begin{equation}\label{eq:Q}
  Q(x_1,\ldots,x_{p+q})=x_1^2+\cdots+x_p^2-x_{p+1}^2-\cdots-x_{p+q}^2
\end{equation}
Let $v_0=(k,0,\ldots,0)\in V$ and $H=\Stab_G(v_0)$. Note that $H$ is
the set of fixed points of the involution of $G$:
\begin{align*}
  \sigma:g\mapsto \diag(-1,1,\ldots,1)\cdot g\cdot
  \diag(-1,1,\ldots,1),
\end{align*}
which commutes with the Cartan involution $\theta(g)= \trn g^{-1}$.
We have the Cartan decomposition $G=K\exp(\la{a}^+)\mathcal{W}H$ where
\begin{align*}
  K=\SO(p)\times \SO(q),\quad \la{a}^+=\R_{\ge 0}\cdot
  (E_{p+q,1}+E_{1,p+q}),
\end{align*}
and $\mathcal{W}=\{e\}$ if $q>1$ and $\mathcal{W}=\diag(\pm
1,1,\ldots,1)$ if $q=1$.  (Here $E_{ij}$ is the matrix with $1$ at
$(i,j)$-entry and $0$ at the other entries.)  For $p>1$, $V=Gv_0$ and
for $p=1$, $V=Gv_0\cup G(-v_0)$.  Setting
$v_0^\infty=(1/\sqrt{2},0,\ldots,0,1/\sqrt{2})$, we have
\begin{align*}
  V^\infty=\left\{
    \begin{tabular}{ll}
      $Kv_0^\infty$ & for $p,q>1$,\\
      $Kv_0^\infty\sqcup -Kv_0^\infty$ & otherwise.
    \end{tabular}
  \right.
\end{align*}
Note that in all cases, $V^\infty\simeq S^{p-1}\times S^{q-1}$, where
we set $S^0=\{\pm 1\}$.  One can check that $a_\iota=p+q-2$ and
$b_\iota=1$.  If $V(\Z)\ne \emptyset$, the limiting distribution of
the points $\{\pi(v):\,v\in V(\Z),\norm{v}<T\}$ as $T\to\infty$ is
given by the probability measure $\frac{dv}{\norm{v}^{p+q-2}}$ where
$dv$ is the suitably normalized invariant measure on $S^{p-1}\times
S^{q-1}$ (see (\ref{eq:mu_iota})); note that $\norm{\cdot}$ can be any
given norm on $\R^{p+q}$.

\subsection{Determinant surface} \label{sec:det} For $k\in
\Z\setminus\{0\}$, let
\begin{align*}
  V=\{v\in\operatorname{M}(n,\R):\,\det (v)=k\}.
\end{align*}
Fix $v_0\in V$.  Note that $V$ is a homogeneous space of
$G=\SL(n,\R)\times \SL(n,\R)$ for the action
\begin{align*}
  (g_1,g_2)\cdot v\mapsto g_1 v (v_0^{-1}g_2^{-1}v_0),\quad
  (g_1,g_2)\in G,\, v\in V,
\end{align*}
and $H=\Stab_G(v_0)$ is the diagonal embedding of $\SL_n(\R)$ in $G$,
which is the fixed point set of the involution
$\sigma(g_1,g_2)=(g_2,g_1)$ commuting with the standard Cartan
involution.  We have Cartan decomposition $G=K\exp(\mathfrak{a}^+)H$
(note that $\mathcal{W}=\{e\}$), where
\begin{align*}
  K&=\SO(n)\times\SO(n),\\
  \la{a}^+&=\{(a,-a):a=\diag(s_1,\ldots,s_n),\sum_i s_i=0,s_i-s_j\ge 0
  \text{ if } i<j\}.
\end{align*}
The simple roots are $\alpha_i=s_i-s_{i+1}$, $ i=1,\ldots,n-1$, the
fundamental weights are $\omega_i=\sum_{j=1}^i s_i$, and
$i=1,\ldots,n-1$, and the highest weight is $\lambda_\iota=2\omega_1$.
Hence, the $\lambda_\iota$-connected subsets of the set of simple
roots are $I_0=\emptyset$ and $I_j=\{\alpha_1,\ldots,\alpha_j\}$,
$j=1,\ldots, n-1$. We have
\begin{align*}
  V^\infty&=\{v\in S(\operatorname{M}(n,\R)):\det (v)=0\},\\
  V^\infty_{I_j}&=\{v\in S(\operatorname{M}(n,\R)) :
  \operatorname{rank}(v)=j+1\}.
\end{align*}
Since
\begin{align*}
  2\rho=2\sum_j j(n-j)\alpha_j\quad\text{and}\quad
  \lambda_\iota=2\omega_1=2\sum_j \frac{n-j}{n}\alpha_j,
\end{align*}
we have $a_\iota=n^2-n$, $b_\iota=1$, $I_\iota=I_{n-2}$. Hence, the
results from Section~\ref{sec:intro} (see Remark~\ref{r:not_sim})
imply that for any admissible cone $\mathcal{C}\subset
\operatorname{M}(n,\R)$ that contains a degenerate matrix in its
interior,
\begin{align*}
  \#\{v\in \operatorname{M}(n,\Z)\cap \mathcal{C}:\,
  \det(v)=k,\norm{v}<T\}\sim_{T\to\infty} c(\mathcal{C},k)\cdot
  T^{n^2-n},
\end{align*}
where $c(\mathcal{C},k)>0$, and the measures
\begin{align*}
  T^{-(n^2-n)}\sum_{v\in V(\Z):\norm{v}<T} \delta_{\pi(v)}
\end{align*}
converge as $T\to\infty$ to a finite smooth measure concentrated on
the set of matrices of rank $n-1$ in $S(\operatorname{M}(n,\R))$.

\subsection{Space of symmetric matrices}\label{sec:quad}
Let $V$ be the space of real symmetric matrices of signature $(p,q)$
of determinant $(-1)^q$. Put $J=\diag(\underbrace{1,\dots,1}_p,
\underbrace{-1,\dots,-1}_q)\in V$.  Then
\begin{align*}
  V=\{g J \trn g:\, g\in\SL(p+q,\R)\}\simeq \SL(p+q,\R)/\SO(p,q).
\end{align*}
Let $n=p+q$, $G=\SL(n,\R)$ and $H=\SO(p,q)$. Note that $V$ is the
orbit of $J$ for the representation $\iota$ of $G$ on the space $W$ of
symmetric $n\times n$ matrices given by
\begin{align*}
  g\cdot w\mapsto g w \trn g, \quad g\in G,\; w\in W.
\end{align*}
Also, $H$ is the the set of fixed points of the involution
$\sigma:g\mapsto J \trn g^{-1}J$, which commutes with the Cartan
involution $\theta:g\mapsto \trn g^{-1}$. We have Cartan decomposition
$G=K\exp(\la{a}^+)\mathcal{W}H$ where $K=\SO(n)$, $\la{a}^+$ is the
standard Weyl chamber in $G$, and $\mathcal{W}$ is the subset of the
monomial matrices which gives coset representatives for
\begin{align*}
  N_K(\la{a})/N_{K\cap H}(\la{a})Z_K(\la{a})\simeq S_{n}/(S_{p}\times
  S_{q})
\end{align*}
where $S_n$ denotes the group of symmetries on $n$ elements.  The
simple roots $\alpha_i$ and the fundamental weights $\omega_i$ are
defined as in Section~\ref{sec:det}, and the highest weight is given
by $\lambda_\iota=2\omega_1$. In particular, it follows that the
$\lambda_\iota$-connected sets are $I_0=\emptyset$,
$I_j=\{\alpha_1,\ldots,\alpha_j\}$, $j=1,\ldots, n-1$, and we have
\begin{align*}
  V^\infty_{I_j}&=\{v\in S(W):\, \operatorname{sign}(v)=(r,s), r+s=j, r\le p, s\le q\},\\
  V^\infty&=\{v\in S(W):\, \operatorname{sign}(v)=(r,s), r+s<n, r\le
  p, s\le q\}.
\end{align*}
Note that in this case, the sets $V^\infty_{I_j}$ are unions of
several orbits of $G$ if $p,q>0$. For example, $V^\infty_{I_{n-1}}$ is
a union of two open orbits which consist of matrices of signature
$(p-1,q)$ and $(p,q-1)$ respectively. One can check (as in
Section~\ref{sec:det}) that $a_\iota=(n^2-n)/2$, $b_\iota=1$,
$I_\iota=I_{n-2}$.  Hence, the results of Section~\ref{sec:intro}
imply that for every admissible cone $\mathcal{C}\subset W$ which
contains a degenerate symmetric matrix in its interior,
\begin{align*}
  \#\{r\in V(\Z):\, \norm{r}<T\}\sim_{T\to\infty} c(\mathcal{C})\cdot
  T^{\frac{n^2-n}{2}},
\end{align*}
where $c(\mathcal{C})>0$, and the measures
\begin{align*}
  T^{-\frac{n^2-n}{2}}\sum_{v\in V(\Z):0<\norm{v}<T} \delta_{\pi(v)}
\end{align*}
converge as $T\to\infty$ to a measure concentrated on the set of
matrices of signature $(p-1,q)$ and $(p,q-1)$ in $S(W)$.

\subsection{Group variety}\label{sec:group}
Let $\bf G$ be a connected $\Q$-simple algebraic group isotropic over
$\R$ and $\iota:{\mathbf G}\to\text{GL}(W)$ an $\R$-irreducible
$\Q$-representation of $\mathbf G$.  We consider the distribution of
integral points in the variety ${\mathbf V}:=\iota({\mathbf G})$. Note
that ${\mathbf V}(\R)$ consists of finitely many orbits of $G={\mathbf
  G}(\R)^\circ$. For simplicity, we make the computation for the orbit
$V=\iota(G)$.

Let $K$ be a maximal compact subgroup of $G$, $\la{a}$ a Cartan
subalgebra associated to $K$, and $\la{a}^+$ a positive Weyl chamber.
We denote by $\Delta$ the set of simple roots of $G$ with respect to
$\mathfrak{a}^+$, and let $\omega_\alpha$, $\alpha\in\Delta$, be the
set of fundamental weights.  We consider the action of
$\tilde{G}=G\times G$ on $V$:
\begin{align*}
  (g_1,g_2)\cdot v\mapsto g_1 v g_2^{-1},\quad (g_1,g_2)\in
  \tilde{G},\, v\in V.
\end{align*}
Then $V\simeq \tilde{G}/H$, where $H=\{(g,g):g\in G\}$.  We have
Cartan decomposition $\tilde{G}=\tilde{K}\exp(\tilde{\la{a}}^+)H$
(note that $\mathcal{W}=\{e\}$), where
\begin{align*}
  \tilde{K}=K\times K\quad\text{and}\quad
  \tilde{\la{a}}^+=\{(a,-a):a\in\la{a}^+\}.
\end{align*}
Note that in this case every $V_I^\infty$ is a single $G$-orbit.  Let
$\rho$ and $\tilde\rho$ be half of the sums of positive roots for $G$
and $\tilde{G}$, and $\lambda_\iota$ and $\tilde{\lambda}_\iota$ be
the highest weights for $\la{a}$ and $\tilde{\la{a}}$ respectively.
Since $\tilde\rho=2\rho$ and $\tilde\lambda_\iota=2\lambda_\iota$, the
parameters $a_\iota$, $b_\iota$, $I_\iota$ are computed as in
(\ref{eq:defi}), and the distribution of integral points is described
by the results from Section~\ref{sec:intro} (see
Remark~\ref{r:not_sim}).  Let us consider ``generic'' case, that is,
\begin{equation}\label{eq:rl}
  2\rho=\sum_{\alpha\in\Delta} u_\alpha \alpha\quad\text{and}\quad
  \lambda_\iota=\sum_{\alpha\in\Delta} m_\alpha\alpha=\sum_{\alpha\in\Delta} n_\alpha\omega_\alpha
\end{equation}
with all $n_\alpha>0$ and $\frac{u_\alpha}{m_\alpha}\ne
\frac{u_\beta}{m_\beta}$ for all $\alpha\ne \beta$. Then the Satake
boundary $V^\infty$ is a union of $2^{\dim\la{a}}$ orbits of $G$, and
there are exactly $\dim\la{a}$ open orbits
$V^\infty_{\Delta\setminus\{\alpha\}}$, $\alpha\in\Delta$, but the
measures
\begin{align*}
  T^{-\max_\alpha (u_\alpha/m_\alpha)}\sum_{v\in V(\Z):\norm{v}<T}
  \delta_{\pi(v)}
\end{align*}
converge as $T\to\infty$ to a measure concentrated on the single open
$G$-orbit $V^\infty_{\Delta\setminus\{\alpha_0\}}$ such that
$\frac{u_{\alpha_0}}{m_{\alpha_0}}=\max_\alpha
\frac{u_\alpha}{m_\alpha}$ (compare with non-generic case in
Section~\ref{sec:det}). The number of integral points in $V$ with norm
less than $T$ is of order $T^{\max_\alpha (u_\alpha/m_\alpha)}$ as
$T\to\infty$, and the number of points whose projections accumulate on
the open $G$-orbit $V^\infty_{\Delta\setminus\{\alpha\}}$ is of order
$T^{u_\alpha/m_\alpha}$ as $T\to\infty$.

\subsection{General affine symmetric variety}
Let $G$ be a connected noncompact semisimple Lie group and $H$ a
symmetric subgroup.  We fix a Cartan decomposition
$$G=K\exp(\la{a}^+)\mathcal{W}H .$$ By Proposition~\ref{p:exist}, given
an integral dominant weight $\lambda$ of $\la{a}^+$, there exists an
$\R$-irreducible $H$-spherical representation $\iota:G\to\GL(W)$ with
the highest weight $\lambda_\iota$ being a multiple of $\lambda$. Then
$W$ contains a symmetric variety $V\simeq G/H$.  The structure of the
Satake boundary $V^\infty$ of $V$ is determined by the combinatorial
data (\ref{eq:rl}) of $2\rho$ and $\lambda_\iota$.  Assume that $G$
and $H$ are the groups of real points of algebraic semisimple
$\Q$-groups $\mathbf G$ and $\mathbf H$ such that $\mathbf G$ is
$\Q$-simple and $\mathbf H$ has no nontrivial $\Q$-characters, and
that $\iota$ is defined over $\Q$. Then if $V(\Z)\ne \emptyset$, the
distribution of integral points $V(\Z)$ is determined by (\ref{eq:rl})
as well. We mention two examples.

The ``generic'' case (i.e., all $n_\alpha>0$ and
$\frac{u_\alpha}{m_\alpha}\ne \frac{u_\beta}{m_\beta}$ for $\alpha\ne
\beta$) is quite similar to the discussion in Section~\ref{sec:group}
except that when $V$ is not a group variety, the sets $V_I^\infty$ may
be unions of several $G$-orbits.

It is well known that $2\rho$ is an integral dominant weight and all
$n_\alpha>0$.  Hence, by Proposition~\ref{p:exist}, there exists an
$H$-spherical representation with the highest weight
$\lambda_\iota=2\ell\rho$ for some $\ell\in\N$. Moreover, if $\mathbf
G$ is an inner form, then the corresponding representation is defined
over $\Q$ (see Remark~\ref{r:q}).  We compute: $a_\iota=1/\ell$,
$b_\iota=\dim\la{a}$, $I_\iota=\emptyset$. Hence, the number of
integral points in $V$ with norm less than $T$ is of order
$T^{1/\ell}(\log T)^{\dim \la{a}-1}$, and the measures
\begin{align*}
  \frac{1}{T^{1/\ell}(\log T)^{\dim \la{a}-1}}\sum_{v\in
    V(\Z):\norm{v}<T} \delta_{\pi(v)}
\end{align*}
converge as $T\to\infty$ to a measure $\mu_\iota$ supported on
$V^\infty_\emptyset$. Note that $K$ acts transitively on
$V^\infty_\emptyset$ (see Proposition~\ref{p:orbit}), and
$\mu_\iota=\norm{v}^{-a_\iota}dv$ where $dv$ is a suitably normalized
$K$-invariant measure on $V^\infty_\emptyset$
(cf.~\eqref{eq:mu_iota}).  On the other hand, for $f\in C(S(W))$ such
that $\supp f\cap V^\infty\subset V^\infty_{\Delta_\sigma-\{\alpha\}}$
for some $\alpha\in\Delta_\sigma$, we have
\begin{align*}
  T^{-1/\ell}\sum_{v\in V(\Z):\norm{v}<T} f(\pi(v))\to
  \mu_{\Delta_\sigma-\{\alpha\}}(f)
\end{align*}
where $\mu_{\Delta_\sigma-\{\alpha\}}$ is a measure concentrated on
$V^\infty_{\Delta_\sigma-\{\alpha\}}$.

\section{Affine symmetric spaces and representations}\label{sec:basic}

\subsection{Affine symmetric spaces} (see
\cite[Ch.~7]{sch},\cite[Part~II]{hs},\cite{os},\cite{ros})

Let $G$ be a connected noncompact semisimple Lie group with finite
center and $\mathfrak{g}$ the Lie algebra of $G$.  A closed subgroup
$H$ of $G$, with the Lie algebra $\mathfrak{h}\subset \mathfrak{g}$,
is called {\it symmetric\/} if $\mathfrak{h}$ is the set of fixed
points of an involution $\sigma$ of $\mathfrak{g}$.  Then the factor
space $G/H$ is called an {\it affine symmetric space}.

There exists a Cartan involution $\theta$ of $\la{g}$ which commutes
with $\sigma$. We denote by $K$ the maximal compact subgroup of $G$
that corresponds to $\theta$ and by $\mathfrak{k}$ its Lie algebra.
We have decompositions
\begin{align*}
  \mathfrak{g}=\mathfrak{h}\oplus\mathfrak{q}\quad\text{and}\quad
  \mathfrak{g}=\mathfrak{k}\oplus\mathfrak{p}
\end{align*}
into $+1$ and $-1$ eigenspaces of $\sigma$ and $\theta$ respectively.

There exists a Cartan subalgebra $\la{c}$ of $\la{g}$ stable under
$\theta$ and $\sigma$ such that $\la{b}:=\la{c}\cap\la{p}$ is a
maximal abelian subalgebra of $\la{p}$, $\la{c}\cap\la{q}$ is a
maximal abelian subalgebra of $\la{q}$, and
$\la{a}:=\la{c}\cap\la{p}\cap \la{q}$ is a maximal abelian subalgebra
of $\la{p}\cap \la{q}$.  We call $\la{b}$ a {\it Cartan subalgebra\/}
associated to $\theta$ and $\la{a}$ a {\it Cartan subalgebra\/}
associated to $(\theta,\sigma)$.  We denote by $\Sigma_\C\subset
\la{c_\C}^*$, $\Sigma\subset \la{b}^*$, and $\Sigma_\sigma\subset
\la{a}^*$ the root systems. One can choose a set of positive roots
$\Sigma_\C^+\subset \Sigma_\C$ so that
$\Sigma^+=\Sigma_\C^+|_\la{b}\setminus \{0\}$ and
$\Sigma_\sigma^+=\Sigma^+|_\la{a}\setminus \{0\}$ are systems of
positive roots in $\Sigma$ and $\Sigma_\sigma$.

Let $\Delta_\C\subset\Sigma_\C^+$ denote the system of simple roots.
Then
\begin{equation}\label{eq:rest}
  \Delta=\Delta_\C|_{\mathfrak{b}}\setminus \{0\}
  \quad\text{and}\quad 
  \Delta_\sigma=\Delta|_{\mathfrak{a}}\setminus \{0\}.
\end{equation}
are systems of simple roots for $\Sigma$ and $\Sigma_\sigma$
respectively.  We also set $\Delta_0=\{\alpha\in\Delta :
\alpha|_\mathfrak{a}=0\}$.

The space $\mathfrak{a}$ is a Cartan subalgebra associated to $\theta$
of the reductive Lie algebra $(\mathfrak{k}\cap\mathfrak{h})\oplus
(\mathfrak{p}\cap\mathfrak{q})$, which is the set of fixed points of
the involution $\sigma\theta$. We denote by
$\Sigma_{\sigma,\theta}\subset \Sigma_\sigma$ the corresponding root
system and choose a set positive roots
$\Sigma_{\sigma,\theta}^+\subset\Sigma_{\sigma,\theta}$ such that
$\Sigma_{\sigma,\theta}^+\subset\Sigma^+_{\sigma}$.

The Weyl groups of $\Sigma_\sigma$ and $\Sigma_{\sigma,\theta}$ are
given by
\begin{align*}
  \mathcal{W}_\sigma=N_K(\mathfrak{a})/Z_K(\mathfrak{a})\quad\text{and}\quad
  \mathcal{W}_{\sigma,\theta}=N_{K\cap H}(\mathfrak{a})/Z_{K\cap
    H}(\mathfrak{a}),
\end{align*}
and one can choose a set $\mathcal{W}\subset N_K(\mathfrak{a})\cap
N_K(\mathfrak{b})$ of coset representatives of
$\mathcal{W}_\sigma/\mathcal{W}_{\sigma,\theta}$

Denoting by $\la{a}^+$ be the closed positive Weyl chamber for
$\Sigma^+_{\sigma}$, we have Cartan decomposition:
\begin{equation}\label{eq:cartan}
  G=K\exp(\mathfrak{a})H= K\exp(\mathfrak{a}^{+})\mathcal{W}H.
\end{equation}
Note that for any $g\in G$, $\mathfrak{a}^{+}$-component of $g$ and
the $\mathcal{W}$-component of $g$ are uniquely defined.

For a root $\alpha\in\Sigma_\sigma\cup \{0\}$, we denote by
$\mathfrak{g}_\alpha$ the corresponding root space associated to
$\la{a}$. Also for a root $\tilde\alpha\in \Sigma$, we denote by
$\la{g}_{\tilde\alpha}(\la{b})$ the corresponding root space
associated to $\la{b}$.

Let $\inpr{\cdot}{\cdot}$ denote the Killing form on $\la{g}$. We
consider a positive definite symmetric bilinear form $B$ on
$\mathfrak{g}$:
\begin{equation}\label{eq:BB}
  B(X, Y)=-\inpr{X}{\theta(Y)} =\operatorname{Tr}(\operatorname{ad} X \circ \operatorname{ad}
  (\theta(Y)).
\end{equation}
Note that
\begin{align}
  &B(\mathfrak{g}_\alpha,\mathfrak{g}_\beta)=0\quad\text{for all $\alpha\ne\beta\in \Sigma_\sigma\cup\{0\}$},\nonumber\\
  &B^\theta=B^\sigma=B.\label{eq:inv}
\end{align}

\begin{rem}
  \label{r:b_gamma} \rm For $\tilde\beta\in\Sigma$, take any $X\in
  \la{g}_{\tilde\beta}(\la{b})$ such that $B(X,X)=1$, and put
  $b_{\tilde\beta}=[X,-\theta(X)]$. Then
  $\theta(b_{\tilde\beta})=-b_{\tilde\beta}$. Since
  $\theta(\la{g}_{\tilde\beta})=\la{g}_{-\tilde\beta}$, we have
  $b_{\tilde\beta}\in \la{g}_0(\la{b})$. Hence
  $b_{\tilde\beta}\in\la{b}$. Moreover for all $b\in\la{b}$,
  \begin{align}
    \label{eq:b_gamma}
    \big\langle b, b_{\tilde\beta}\big\rangle
    =\inpr{b}{[X,-\theta(X)]} =\inpr{[b,X]}{-\theta(X)}
    =\tilde\beta(b)\inpr{X}{-\theta(X)} =\tilde\beta(b).
  \end{align}
  Since the Killing form restricted to $\la{b}$ is nondegenerate,
  $\la{b}_{\tilde\beta}$ is the unique element, say $b_\ast$, of
  $\la{b}$ such that $\inpr{b}{b_\ast} = \tilde\beta(b)$ for all
  $b\in\la{b}$.
\end{rem}

For each $\alpha\in\Sigma_\sigma$, the root space
$\mathfrak{g}_\alpha$ is invariant under the involution
$\sigma\theta$, and it decomposes into $(\pm 1)$-eigenspaces of
$\sigma\theta$:
\begin{align*}
  \mathfrak{g}_\alpha=\mathfrak{g}_\alpha^+\oplus
  \mathfrak{g}_\alpha^-;
\end{align*}
we define
\begin{align*}
  l^\pm_\alpha=\dim \mathfrak{g}_\alpha^\pm, \quad \text{and}\quad
  l_\alpha=l_\alpha^++l_\alpha^-.
\end{align*}
We have
\begin{align*}
  2\rho=\sum_{\alpha\in\Sigma^+_\sigma} l_\alpha \alpha.
\end{align*}

A Haar measure on $G/H$ is given by the formula
\begin{equation}\label{eq:volume}
  \int_{G/H} f\, d\mu = \int_K\sum_{w\in\mathcal{W}} \int_{\mathfrak{a}^{+}}
  f(k\exp(a)wH)\xi(a)\,dadk,\quad f\in C_c(G/H),
\end{equation}
where $da$ and $dk$ denote Haar measures on $\mathfrak{a}$ and $K$,
and
\begin{equation}\label{eq:xi}
  \xi(a)=\prod_{\alpha\in \Sigma^+_\sigma}
  (\sinh\alpha(a))^{l_\alpha^+}(\cosh\alpha(a))^{l_\alpha^-}.
\end{equation}
To match (\ref{eq:volume}) with the integral formula given in
\cite[Ch.~7]{sch}, we note that the function $|\xi|$ is invariant
under the Weyl group.

\begin{rem} \label{r:a_alpha} \rm For $\alpha\in\Sigma_\sigma$, let
  $X\in \la{g}_\alpha^+\cup \la{g}_\alpha^-$ such that $B(X,X)=1$ and
  put $a_\alpha:=[X,-\theta(X)]\in \la{g}_0$. Then
  $\sigma(a_\alpha)=-a_\alpha$ and $\theta(a_\alpha)=-a_\alpha$.
  Therefore $a_\alpha\in\la{a}$ and by \eqref{eq:b_gamma}, we have
  $\inpr{ a}{a_\alpha}=\alpha(a)$ for all $a\in\la{a}$.  Hence
  $a_\alpha$ is the unique element of $\la{a}$ satisfying the last
  equation.
\end{rem}

\subsection{Representations} (see
\cite[Ch.~IV]{gjt},\cite{sat},\cite{cp})\label{sec:rep}
\newcommand{\Wc}{W_0} Let $\iota: G\to \GL(W)$ be an irreducible over
$\R$ representation of $G$ on a real vector space $W$.  We denote by
$\mathfrak{g}_\C$, $\mathfrak{c}_\C$, $\mathfrak{h}_\C$ the
complexifications of $\mathfrak{g}$, $\mathfrak{c}$, $\mathfrak{h}$.
Note that $\sigma$ extends to an involution of $\mathfrak{g}_\C$, and
$\mathfrak{h}_\C$ is the subalgebra of the fixed points of $\sigma$ in
$\mathfrak{g}_\C$.

Let $\Wc$ be a complex $\la{g}$-irreducible subspace of $\c\otimes W$.
Then either $\c\otimes W=\Wc$ or $\c\otimes W=\Wc\oplus \bar \Wc$,
where bar denotes the standard complex conjugation on $\c\otimes W$.
Note that if $\c \otimes W$ is not complex irreducible, then the map
\begin{equation} \label{eq:isom} \Wc\ni v\mapsto (v+\bar v)\in W
\end{equation}
is a $\la{g}$-equivariant isomorphism over $\R$; and hence in this
case $W$ can be treated as vector space over $\c$ with $\c$-linear
action of $\la{g}$. By abuse of notation, the representation of
$\la{g}$ on $\Wc$ over $\c$ will also be denoted by $\iota$.

We denote by $\Lambda_\iota\in \mathfrak{c}_\C^*$ the highest weight
of $\iota$ with respect to the ordering defined by $\Delta_\C$. Then
all other weights of $\la{c}_\C$ with respect to $\iota$ are of the
form
\begin{equation}\label{eq:lambda0}
  \lambda=\Lambda_\iota-\sum_{\alpha\in\Delta_\C}
  n_\alpha(\lambda)\alpha
\end{equation}
for some non-negative integers $n_\alpha(\lambda)$.

The action of $\mathfrak{a}$ on $W$ is diagonalizable (over $\R$) and
\begin{align*}
  W=\oplus_{\lambda\in\Phi_\iota} W^\lambda,
\end{align*}
where $\Phi_\iota\subset \mathfrak{a}^*$ is the set of weights and
\begin{align*}
  W^\lambda=\{w\in W: a\cdot w=\lambda(a)w,\, \forall
  a\in\mathfrak{a}\}
\end{align*}
denotes the weight space with weight $\lambda$. Given $w\in W$, we
have a decomposition
\begin{align*}
  w=\sum_{\lambda\in\Phi_\iota} w^\lambda, \quad w^\lambda\in
  W^\lambda.
\end{align*}

The weight $\lambda_\iota:=\Lambda_\iota|_\la{a}$ is the maximal
element of $\Phi_\iota$ with respect to the ordering defined by
$\Delta_\sigma$.  All the other weights $\lambda\in\Phi_\iota$ are of
the form
\begin{equation}\label{eq:lambda}
  \lambda=\lambda_\iota-\sum_{\alpha\in\Delta_\sigma} n_\alpha(\lambda)\alpha
\end{equation}
for some non-negative integers $n_\alpha(\lambda)$. Let
\begin{align*}
  \supp \lambda=\{\alpha\in\Delta_\sigma: n_\alpha(\lambda)>0\}.
\end{align*}
For a subset $I$ of $\Delta_\sigma$ and a vector $w\in W$, we set
\begin{align*}
  w^I=\sum_{\lambda: \, \supp \lambda \subset I} w^\lambda
  \quad\text{and}\quad W^I=\sum_{\lambda: \supp \lambda\subset I}
  W^\lambda.
\end{align*}

Recalling Definition~\ref{d:lambda_conn}, a subset of $\mathfrak{a}^*$
is called {\it connected\/} if it is not a union of nonempty subsets
orthogonal with respect to the form $B$; that is, if its {\em Dynkin
  diagram}\/ is connected. We say that $I\subset \Delta_\sigma$ is
$\lambda_\iota$-connected, if $I\cup\{\lambda_\iota\}$ is a connected
subset of $\Delta_\sigma$.

\begin{Prop}\label{p:satake}
  For any $\lambda\in \Phi_\iota$, $\supp(\lambda)\cup
  \{\lambda_\iota\}$ is connected, and for every
  $\lambda_\iota$-connected $I\subset \Delta_\sigma$ there exists
  $\lambda\in \Phi_\iota$ such that\/ $\supp(\lambda)=I$.
\end{Prop}

\begin{proof}
  The similar statement for the set of weight of $\mathfrak{b}$ was
  shown in \cite[Sec. 2]{sat}, and the proof applies to our situation
  with minor changes. The key fact is that there exists an involution
  $\alpha\mapsto \alpha'$ of the set $\Delta-\Delta_0$ such that
  \begin{align*}
    -\alpha^\sigma=\alpha'+\sum_{\beta\in\Delta_0}
    n_{\alpha,\beta}\beta,\quad \alpha\in \Delta\setminus\Delta_0,
  \end{align*}
  for some $n_{\alpha,\beta}\in\Z_{\ge 0}$ (see \cite[Lemma
  7.2.3]{sch}).  Using that the proposition holds for the weights of
  $\mathfrak{b}$, one can complete the proof as in \cite{sat}.
\end{proof}

\vspace{0.1cm}

\begin{rem} \label{r:K} We set $\K=\R$ when $W\otimes\c=W_0$, and
  $\K=\c$ when $W\otimes\c=W_0\oplus \bar{W}_0$. Then $W$ can be treated
  as a $\K$-vector space with $\K$-linear action of $\la{g}$. % In both
%   cases, $W^{\lambda_\iota}$ is a one-dimensional $\K$-subspace.
\end{rem}

\vspace{0.1cm}

\subsection*{$H$-spherical representations}

Let $W^H$ denote the space of fixed points of $H$ on $W$. If $W^H\neq
0$, then the representation $\iota$ is called {\em $H$-spherical}.

\begin{Lem}[{\cite[Lemma 1.5]{cp}}] \label{l:cp} If $\iota$ is
  $H$-spherical then the $\K$-$\dim(W^{\la{h}})=1$; and
  $\Lambda_\iota^\sigma=-\Lambda_\iota$.  \qed
\end{Lem}

Using the form $B$ defined in (\ref{eq:BB}), we introduce a scalar
product on $\la{c}^*$, $\la{b}^*$, $\la{a}^*$. An element $\lambda\in
\la{c}^*$ is called {\it integral\/} if
$\frac{2\inpr{\lambda}{\alpha}}{\inpr{\alpha}{\alpha}}\in \Z$ for all
$\alpha\in \Delta_\C$, and it is called {\it dominant\/} if
$\inpr{\lambda}{\alpha}\ge 0$ for all $\alpha\in \Delta_\C$.  For
$\beta\in \Delta_\C$, we define the fundamental weights $\omega_\beta$
by
\begin{align*}
  \frac{2\inpr{\omega_\beta}{\alpha}}{\inpr{\alpha}{\alpha}}=\delta_{\alpha\beta},
  \quad \forall \alpha\in\Delta_\C,
\end{align*}
where $\delta_{\alpha\beta}$ denotes the Kronecker symbol.  Similarly,
we define these notions for $\la{b}^*$ and $\la{a}^*$.  It is
well-known that the the highest weight $\Lambda_\iota$ is integral and
dominant, and conversely, every integral dominant weight is the
highest weight of an irreducible representation of $\la{g}_\C$.  We
prove an analogous result for real spherical representations:

\begin{Prop} [cf. {\cite[Proposition 4.15]{gjt}}]\label{p:exist}
 The highest weight
  $\lambda_\iota$ is integral and dominant.  There exists $\ell\in\N$
  such that for every integral dominant $\lambda\in \la{a}^*$,
  $\ell\lambda$ is a highest weight of a real absolutely irreducible
  $H$-spherical representation of $\la{g}$.
\end{Prop}

\begin{proof}
  The fact that $\lambda_\iota$ is integral and dominant follows from
  the representation theory of $\la{sl}(2,\R)$ (see \cite[Lemma
  4.12]{gjt}).

  For $\alpha\in \Delta_\C$ (or $\alpha\in \Delta_\sigma$), we take
  $h_\alpha\in \la{c}$ and $h_\alpha^*\in \la{c}^*$ (or $h_\alpha\in
  \la{a}$ and $h_\alpha^*\in \la{a}^*$) such that
  \begin{align*}
    \inpr{h_\beta}{h_\alpha}=\beta(h_\alpha)=\inpr{\beta}{\alpha}
    \quad\text{and}\quad
    \inpr{h_\alpha^*}{\beta}=h^*_\alpha(h_\beta)=\delta_{\alpha\beta}
  \end{align*}
  for all $\alpha,\beta\in \Delta_\C$ (or $\alpha,\beta\in
  \Delta_\sigma$).  For $\lambda\in \la{c}^*$, we denote by
  $\bar\lambda\in \la{a}^*$ its restriction to $\la{a}$, and for $x\in
  \la{c}$, we denote by $\bar x\in\la{a} $ its orthogonal projection
  to $\la{a}$. It follows from (\ref{eq:rest}) and (\ref{eq:inv}) that
  for $\alpha\in \Delta_\C$ such that $\bar\alpha\ne 0$ and $\beta\in
  \Delta_\sigma$,
  \begin{align*}
    \bar h_\alpha=h_{\bar\alpha} \quad\text{and}\quad
    h_\beta^*=\sum_{\alpha\in\Delta_\C: \bar\alpha=\beta} \bar
    h^*_\alpha.
  \end{align*}
  Suppose that $\lambda=\sum_{\beta\in\Delta_\sigma}
  n_\beta\omega_\beta$ for $n_\beta\in\Z_{\ge 0}$. Since
  $\omega_\beta=\frac{1}{2}\inpr{\beta}{\beta}h^*_\beta$, we have
  \begin{equation}\label{eq:lll}
    \lambda = \sum_{\alpha\in\Delta_\C: \bar\alpha\ne 0} n_{\bar\alpha}
    \frac{\inpr{\bar\alpha}{\bar\alpha}}{\inpr{\alpha}{\alpha}}\cdot 
    \bar \omega_\alpha.
  \end{equation}
  It is well known that $\inpr{\alpha_1}{\alpha_2}\in\Q$ for
  $\alpha_1,\alpha_2\in \Sigma_\C$. Hence, using that
  \begin{align*}
    \bar\alpha=\frac{1}{4}(\alpha-\alpha^\theta-\alpha^\sigma+\alpha^{\sigma\theta}),
  \end{align*}
  we deduce that the coefficients in (\ref{eq:lll}) are rational
  numbers. We take $\ell\in\N$ such that
  \begin{align*}
    \ell\lambda=\sum_{\alpha\in\Delta_\C: \bar\alpha\ne 0} m_{\alpha}
    \bar \omega_\alpha
  \end{align*}
  for $m_\alpha\in 2\Z_{\ge 0}$ and consider an irreducible complex
  representation $\iota:\la{g}_\C\to\GL({\Wc})$ with the highest
  weight
  \begin{equation}\label{eq:Lambda}
    \Lambda=\sum_{\alpha\in\Delta_\C: \bar\alpha\ne 0}
    m_{\alpha} \omega_\alpha.
  \end{equation}

  Let $\Delta_\C^\theta=\{\alpha\in\Delta_\C:\,
  \alpha|_\mathfrak{b}=0\}$.  It was shown in \cite{sat} that there
  exists an involution $\alpha\mapsto \tilde\theta(\alpha)$ of the set
  $\Delta_\C\setminus\Delta_\C^\theta$ such that
  \begin{align*}
    -\alpha^\theta=\tilde\theta(\alpha)+\sum_{\beta\in\Delta_\C^\theta}
    u_{\alpha,\beta}\beta,\quad \alpha\in
    \Delta_\C\setminus\Delta_\C^\theta,
  \end{align*}
  for some $u_{\alpha,\beta}\in\Z_{\ge 0}$. Moreover, according to
  \cite[\S9]{on}, the involution $\tilde\theta$ is induced by an
  automorphism of the Dynkin diagram of $\Delta_\C$. In particular,
  $\big\langle\tilde\theta(\alpha),\tilde\theta(\alpha)\big\rangle =
  \inpr{\alpha}{\alpha}$ for all
  $\alpha\in\Delta_\C\setminus\Delta_\C^\theta$. Also, it is clear
  that $\alpha|_\la{a}=\tilde\theta(\alpha)|_\la{a}$. This shows that
  $m_\alpha=m_{\tilde\theta(\alpha)}$ and by \cite[\S8]{on}, the
  restriction of $\iota$ to $\la{g}$ leaves the a real form $W$ of
  ${\Wc}$ invariant.

  It remain to check that the representation $\iota$ is spherical.
  Recall that $\la{d}:=\la{c}\cap\la{q}=\{x\in\la{c}:\sigma(x)=-x\}$
  is a maximal abelian subalgebra of $\la{q}$.  Let
  $\Delta_\C^\sigma=\{\alpha\in\Delta_\C:\, \alpha|_\mathfrak{d}=0\}$.
  It was shown in \cite{cp} that there exists an involution
  $\alpha\mapsto \tilde\sigma(\alpha)$ of the set
  $\Delta_\C\setminus\Delta_\C^\sigma$ such that
  \begin{equation}\label{eq:tilde_sigma}
    -\alpha^\sigma=\tilde\sigma(\alpha)+\sum_{\beta\in\Delta_\C^\sigma}
    v_{\alpha,\beta}\beta,\quad \alpha\in \Delta_\C\setminus\Delta_\C^\sigma,
  \end{equation}
  for some $v_{\alpha,\beta}\in\Z_{\ge 0}$. Since $m_\alpha$'s are
  even, according to \cite{cp}, the representation $\iota$ is
  spherical provided that $m_\alpha=m_{\tilde\sigma(\alpha)}$ for
  $\alpha\in \Delta_\C\setminus\Delta_\C^\sigma$.  Hence, it suffices
  to show that the involution $\tilde\sigma$ is induced by an
  automorphism of the Dynkin diagram of $\Delta_\C$.  Without loss of
  generality, we may assume that $\Delta_\C^\sigma\ne \emptyset$.
  Then one can check that $\la{c}\cap \la{h}$ is a Cartan subalgebra
  of $[\la{z}_\la{h}(\la{d}), \la{z}_\la{h}(\la{d})]$ with the system
  of simple roots $\Delta_\C^\sigma$. The corresponding Weyl group
  $\mathcal{W}_\C^\sigma$ is generated by reflections
  \begin{align*}
    w_\beta(\alpha) = \alpha -
    2\frac{\inpr{\alpha}{\beta}}{\inpr{\beta}{\beta}}\beta, \quad
    \beta\in\Delta_\C^\sigma.
  \end{align*}
  This implies that for every $w\in \mathcal{W}_\C^\sigma$,
  \begin{align}
    (\Sigma^+_\C\setminus\langle\Delta_\C^\sigma\rangle)^w &
    \subset\Sigma^+_\C\setminus \langle\Delta_\C^\sigma\rangle,\label{eq:w1}\\
    \alpha^w&\in \alpha+\langle\Delta_\C^\sigma\rangle,\quad
    \alpha\in\Delta_\C. \label{eq:w2}
  \end{align}
  Take $w_0\in\mathcal{W}_\C^\sigma$ such that
  $(\Delta_\C^\sigma)^{w_0}=-\Delta_\C^\sigma$.  It follows from
  (\ref{eq:tilde_sigma}) and (\ref{eq:w1}) that the map $\alpha\mapsto
  -\alpha^{\sigma w_0}$ preserves $\Sigma_\C^+$, and hence, it induces
  an automorphism of the Dynkin diagram of $\Delta_\C$. On the other
  hand, it follows from (\ref{eq:tilde_sigma}) and (\ref{eq:w2}) that
  for $\alpha\in\Delta_\C\setminus\Delta_\C^\sigma$,
  \begin{align*}
    -\alpha^{\sigma
      w_0}\in\tilde\sigma(\alpha)+\langle\Delta_\C^\sigma\rangle.
  \end{align*}
  This implies that $\tilde\sigma(\alpha)=-\alpha^{\sigma w_0}$,
  $\alpha\in\Delta_\C\setminus\Delta_\C^\sigma$, and finishes the
  proof.
\end{proof}

\begin{rem}\label{r:q} \rm Suppose that $G={\mathbf G}(\R)^o$ for
  a semisimple algebraic $\Q$-group $\mathbf G$. Choosing the Cartan
  subalgebra $\la{c}$ to be defined over $\Q$, we have the
  $\star$-action of the Galois group $\operatorname{Gal}(\bar{\Q}/\Q)$
  on the set of simple roots $\Delta_\C$. By \cite[Theorem~3.2]{tits},
  the representation constructed in Proposition~\ref{p:exist} is
  defined over $\Q$ provided that the highest weight $\Lambda$ is in
  the root lattice, and the coefficients in (\ref{eq:Lambda}) are
  invariant under the $\star$-action. In particular, if $\mathbf G$ is
  an inner form, then the $\star$-action is trivial, and $\ell
  \lambda$ is realized as a highest weight of a representation defined
  over $\Q$ for some $\ell$.
\end{rem}

\section{Structure of Satake compactification} \label{sec:satake}

Let $G$ be a connected noncompact semisimple Lie group with a finite
center and $\iota:G\to \text{GL}_\R(W)$ an irreducible almost faithful
representation of $G$ on a finite dimensional real vector space $W$.
Let $\sigma$ be an involution of $G$ and $H$ the symmetric subgroup of
$G$ with respect to $\sigma$.  We assume that $H$ fixes a nonzero
$v_0\in W$.

We start with some basic observations: let
$$\la{n}=\sum_{\alpha\in\Sigma_\sigma^+} \la{g}_{\alpha} .$$
\begin{lem} \label{l:irr} We have $W^{\lambda_{\iota}}=\{v\in W:\la{n}
  v=0\}$.
\end{lem}

\begin{proof}
  Let $\alpha\in\Sigma_\sigma^+$. Then $\la{g}_\alpha
  W^{\lambda_\iota}\subset W^{\lambda_\iota+\alpha}$.  Since
  $\lambda_\iota$ is the highest weight in $\Phi_\iota$, we conclude
  that $W^{\lambda_\iota+\alpha}=0$. This shows that
  $\la{n}W^{\lambda_\iota}=0$.

  Now let $v\in W$ such that $\la{n}v=0$. Suppose $v\notin
  W^{\lambda_\iota}$. Then there exists $y\in
  W':=\sum_{\lambda<\lambda_\iota}W^{\lambda}$ such that $\la{n}y=0$.
  Let $\la{n}^-=\sum_{\alpha\in \Sigma_\sigma^+} \la{g}_{-\alpha}$.
  Then $\la{g}=\la{n}^-\oplus\la{g}_0\oplus\la{n}$.  Note that
  \begin{align*}
    U_0(\la{n})y=0;\quad U_0(\la{g}_0)W'\subset W'; \quad
    U_0(\la{n}^-)W'\subset W',
  \end{align*}
  where $U_0(\la{n})$ denotes the linear span of (non-constant)
  monomials formed from a basis of $\la{n}$, and the others are
  defined similarly. By Poincare-Birkhoff-Witt's theorem, it follows
  that
  \begin{align*}
    U(\la{g})y\subset W',
  \end{align*}
  where $U(\la{g})$ is the universal enveloping algebra of $\la g$.
  This is a contradiction, because $U(\la{g})y=W$ by the
  irreducibility of the $\la{g}$-action on $W$.
\end{proof}

\begin{lem}\label{l:weights}
  For every $\lambda_\iota$-connected $I\subset\Delta_\sigma$ and
  every $w\in\mathcal{W}$, there exists $\lambda\in \Phi_\iota$ such
  that $\supp \lambda=I$ and $(wv_0)^\lambda\ne 0$.
\end{lem}

\begin{proof}
  First, we consider the case of $I=\emptyset$; that is, we show that
  $(wv_0)^{\lambda_\iota}\ne 0$. We denote by
  $\sigma_w=\Ad(w)\circ\sigma\circ{\Ad}(w^{-1})$ the involution of
  $\la{g}$ corresponding to the symmetric subgroup $wHw^{-1}$. Take a
  maximal $\lambda\in \Phi_{\iota}$ such that $(wv_0)^{\lambda}\ne 0$
  and suppose that $\lambda\ne \lambda_\iota$.  Then by
  Lemma~\ref{l:irr} there exist $\alpha\in\Sigma^+_\sigma$ and
  $X\in\la{g}_\alpha$ such that $X(wv_0)^\lambda\ne 0$. Since
  $X+\sigma_w(X)$ belongs to the Lie algebra of $wHw\inv$,
  \begin{align*}
    (X+\sigma_w(X))(wv_0)=0.
  \end{align*}
  Therefore there exists $\mu\in\Phi_{\iota}$ such that $(wv_0)^\mu\ne
  0$ and $\lambda+\alpha=\mu+\alpha^{\sigma_w}$. Since
  $\sigma_w(a)=-a$ for all $a\in\la{a}$, we have
  $\alpha^{\sigma_w}=-\alpha$. Therefore
  $\mu=\lambda+2\alpha>\lambda$, which contradicts the choice of
  $\lambda$.

  Now we prove the general case.  Given a $\lambda_\iota$-connected
  $I\subset\Delta_\sigma$, there exists $w_0\in\mathcal{W}_\sigma$
  such that the weight $\lambda_\iota\circ{\Ad}(w_0)$ has support
  equal to $I$ (\cite[Lemma B.8]{gjt}).  Then by the above case
  \begin{align*}
    w_0(wv_0)^{\lambda_\iota\circ\Ad w_0}=(w_0wv_0)^{\lambda_\iota}\ne
    0.
  \end{align*}
  This proves the lemma.
\end{proof}

\subsection{Symmetric subgroup as a stabilizer}

\begin{Prop} \label{p:proper} The map $G/H\to V$ given by $gH\mapsto
  gv_0$, for all $g\in G$, is proper. In particular, the orbit $Gv_0$
  is closed.
\end{Prop}

\begin{proof}
  Take any $w\in\cW$. By Lemma~\ref{l:weights}
  $(wv_0)^{\lambda_\iota}\neq 0$. Since the representation is almost
  faithful,
  \begin{align*}
    \lambda_\iota = \sum_{\alpha\in\Delta_\sigma} m_\alpha \alpha,
  \end{align*}
  where $m_\alpha>0$ for all $\alpha$. Therefore the map
  \begin{align*}
    a\mapsto \exp(a)(wv_0)=\sum_{\lambda\in\Phi_\iota}
    e^{\lambda(a)}(wv_0)^\lambda
  \end{align*}
  from $\la{a}^+$ to $V$ is proper. Now, since $G=K\exp(\la{a}^+)\cW
  H$, the map $g\mapsto gv_0$ is proper.
\end{proof}

\begin{Prop}
  \label{p:max_stab}
  $H$ is a subgroup of finite index in $\Stab_G(v_0)$.
\end{Prop}

\begin{proof}
  Let $L=\Stab_G(v_0)$. Then $H\subset L$, and by
  Proposition~\ref{p:proper} $L/H$ is compact. Therefore, since $H$ is
  reductive and $L$ is a real almost algebraic subgroup of $G$, we
  conclude that $L$ unipotent radical of $L$ is trivial.  Hence $L$ is
  reductive.  Let $\la{l}$ denote the Lie subalgebra of $\la{g}$
  associated to $L$, and
  $\la{L}^\perp=\{X\in\la{g}:\inpr{X}{\la{l}}=0\}$. Since $\la{l}$ is
  reductive, the Killing form of $\la{g}$ restricted to $\la{l}$ is
  nondegenerate. Therefore we get
  \begin{equation}
    \label{eq:reductive-perp}
    [\la{l},\la{l}^\perp]\subset\la{l}^\perp
    \quad \text{and} \quad
    \la{G}=\la{L}\oplus\la{L}^\perp.
  \end{equation}
  Since $H$ is a symmetric subgroup and $H\subset L$, we note that
  \begin{equation}
    \label{eq:comm}
    \la{h}\cap \la{l}^\perp=\{0\}, \ 
    [\la{h}^\perp,\la{h}^\perp]\subset\la{h}, \ 
    [\la{h},\la{h}]\subset  \la{h}, \  \text{and} \ 
    [\la{h},\la{h}^\perp]\subset\la{h}^\perp.
  \end{equation}
  Put, $\la{m}=\la{h}^\perp\cap \la{l}$. Then
  \begin{align*}
    % \label{eq:com4}
    [\la{l}^\perp, \la{m}] & \subset [\la{l}^\perp,\la{h}^\perp] \cap
    [\la{l}^\perp,\la{l}] \subset [\la{h}^\perp,\la{h}^\perp]\cap
    [\la{l}^\perp,\la{l}]
    \subset \la{h}\cap\la{l}^\perp=\{0\}, \\
    [\la{h},\la{m}] &\subset [\la{h},\la{h}^\perp]\cap [\la{h},\la{l}]
    \subset \la{h}^\perp\cap\la{l}=\la{m},\\
    [\la{m},\la{m}] &\subset [\la{h}^\perp,\la{h}^\perp]\subset\la{h}.
  \end{align*}
  We put $\la{m}_1=[\la{m},\la{m}]\oplus\la{m}$. Then $\la{m}_1$ is
  stable under $\ad(\la{h})$, $\ad(\la{l}^\perp)$, and $\ad(\la{m})$.
  Since $\la{h}^\perp=\la{m}+\la{l}^\perp$, we conclude that
  $\la{m}_1$ is an ideal in $G$.  Since $\la{m}_1\subset\la{l}$, we
  have that $\la{m}_1\cdot v_0=0$. Hence
  \begin{align*}
    \la{m}_1\cdot Xv_0\subset X\cdot\la{m}_1\cdot
    v_0+[\la{m}_1,X]\cdot v_0=\{0\}, \quad \forall X\in\la{g}.
  \end{align*}
  Therefore, since $\la{G}$ acts irreducibly on $W$, $\la{m}_1$ acts
  trivially on $W$. Since $G$ acts almost faithfully on $W$, we get
  $\la{m}_1=\{0\}$, and hence $\la{l}=\la{h}$. Now the conclusion of
  the proposition follows because $L/H$ is compact.
\end{proof}

Using Proposition~\ref{p:proper} and the proof of
Proposition~\ref{p:max_stab} it is straightforward to deduce the
following.

\begin{cor} \label{cor:gen:finite} Suppose that $G$ acts linearly and
  almost faithfully on a finite dimensional vector space $E$ over
  $\K=\R$ or $\C$. Suppose that there exists $0\neq w_0\in E$ such
  that $Hw_0=w_0$ and $\K$-$\spn(Gw_0)=E$. Then the map $gH\mapsto gv$
  from $G/H$ to $E$ is proper. Moreover $H$ is a subgroup of finite
  index in $\Stab_G(w_0)$.
\end{cor}

Let $S(W)$ denote the unit sphere in $W$, and $\pi:W\setminus\{0\}\to
S(W)$ denote the radial projection.

\begin{Prop}\label{p:homeo}
  The map $\pi: V\to \pi(V)$ is a homeomorphism.
\end{Prop}

\begin{proof}
  To verify that the map $\pi$ is bijective, we suppose that
  $g_1v_0=\lambda g_2v_0$ for some $g_1,g_2\in G$ and $\lambda\ne \pm
  1$. Then it follows that for some $g\in G$ and $\lambda\in (-1,1)$,
  $gv_0=\lambda v_0$. Therefore $g^nv_0\to 0$ as $n\to\infty$, which
  contradicts the conclusion of Proposition~\ref{p:proper} that $Gv_0$
  is closed.

  It is clear that the map is continuous and $G$-equivariant.  Since
  the orbits of $G$ in the projective space of $W$ are locally closed,
  it follows that $\pi(G\cdot v_0)$ is locally compact. Hence, the map
  $\pi$ is a homeomorphism.
\end{proof}

\subsection{Satake Boundary}

We define the {\it Satake boundary}\/ $V^\infty$ of $V$ to be the set
of the limit points of the sequences $\pi(v_n)$, $v_n\in V$,
$v_n\to\infty$. Note that identifying $G/H$ with $\pi(V)$, the space
$\pi(V)\cup V^\infty$ gives a compactification of $G/H$ similar to the
Satake compactification of the Riemannian symmetric space of $G$
constructed in \cite{sat}.

We use notations from Section~\ref{sec:basic}.  For
$J\subset\Delta_\sigma$, let $\mathfrak{a}_J=\ker(J)$,
$\mathfrak{a}^J$ its orthogonal complement, and
\begin{align} \label{eq:aJplus} \mathfrak{a}^{J,+}&=\{a\in
  \mathfrak{a}^J:\,\alpha(a)\ge 0 \quad \text{for all $\alpha\in
    J$}\}.
\end{align}
The set $J$ is the system of simple roots on $\mathfrak{a}^J$, and its
Weyl group $\mathcal{W}_J$ can be identified with the subgroup of
$\mathcal{W}_\sigma$ that acts trivially on $\mathfrak{a}_J$.  We
choose a set $\mathcal{W}^J$ of representatives of the double cosets
$\mathcal{W}_J\backslash
\mathcal{W}_\sigma/\mathcal{W}_{\sigma,\theta}$.  In particular,
$\mathcal{W}=\mathcal{W}^\emptyset$.

For $J\subset\Delta_\sigma$ and $w\in\mathcal{W}$, we set
\begin{align*}
  V^\infty_{J,w} = \left\{ \lim \pi(k\exp(a)wv_0) :
    \begin{tabular}{l}
      $k\in K$, $a\in \mathfrak{a}^{+}$, \\
      $\alpha(a)\to\infty$ for $\alpha\in\Delta_\sigma \setminus J$,\\
      $\alpha(a)$ is bounded for $\alpha\in J$.
    \end{tabular}
  \right\}.
\end{align*}

The main result of this section is the following theorem, which gives
an explicit combinatorial description of the decomposition of
$V^\infty$ into $G$-orbits.

Define
\begin{align*}
  \mathcal{O}_{J,w}=\bigcup_{w_1\in \mathcal{W}_J} V^\infty_{J,w_1w}.
\end{align*}

\begin{thm}\label{thm:satake}
  The decomposition of $V^\infty$ into $G$-orbits is given by
  \begin{align*}
    V^\infty=\bigcup_{I,w} \mathcal{O}_{I,w}
  \end{align*}
  where the union is taken over all $\lambda_\iota$-connected subsets
  $I\subsetneq\Delta_\sigma$ and $w\in \mathcal{W}^I$.

  Moreover,
  \begin{equation}\label{eq:oiw}
    \mathcal{O}_{I,w}=\pi(G(wv_0)^I)
  \end{equation}
  and
  \begin{equation}\label{eq:oiw2}
    \mathcal{O}_{I_1,w_1}\cap
    \mathcal{O}_{I_2,w_2}=\emptyset \quad \text{for $I_1\ne I_2$.}
  \end{equation}
\end{thm}

We shall prove this theorem through a series of auxiliary results.

\begin{Prop}\label{p:v_inf}
  Let $J\subset \Delta_\sigma$, $w\in\mathcal{W}$, and $I\subset J$ be
  the largest $\lambda_\iota$-connected subset. Then
  \begin{align*}
    V^\infty_{J,w} = \pi(K \exp(\mathfrak{a}^{I,+})(wv_0)^I) = \pi(K
    \exp(\mathfrak{a}^{J,+})(wv_0)^J)=V^\infty_{I,w}.
  \end{align*}
\end{Prop}

\begin{proof}
  Recall that $V^\infty_{J,w}$ is the set of limit points of the
  sequences $\pi(k\exp(a_n)wv_0)$ where $k\in K$ and $\{a_n\}\subset
  \mathfrak{a}^{+}$ such that $\alpha(a_n)\to\infty$ for $\alpha\in
  \Delta_\sigma \setminus J$ and $\alpha(a_n)$ is bounded for
  $\alpha\in J$.  Passing to a subsequence, we may assume that there
  exists $a\in\mathfrak{a}^{I,+}$ such that $\alpha(a_n)\to\alpha(a)$
  for every $\alpha\in I$. Then for
  $\lambda=\lambda_\iota-\sum_{\alpha\in\Delta_\sigma}
  n_\alpha(\lambda)\alpha\in\Phi_\iota$,
  \begin{align*}
    \sum_{\alpha\in\Delta_\sigma}n_\alpha(\lambda)\alpha(a_n) \to
    \left\{
      \begin{tabular}{ll}
        $\sum_{\alpha\in\Delta_\sigma}n_\alpha(\lambda)\alpha(a)$ 
        & if $\supp(\lambda)\subseteq I$,\\
        $+\infty$ 
        & if $\supp(\lambda)\nsubseteq J$.
      \end{tabular}
    \right.
  \end{align*}
  Note that by Proposition~\ref{p:satake}, $\supp(\lambda)\subseteq J$
  iff $\supp(\lambda)\subseteq I$.  By Lemma~\ref{l:weights},
  $(wv_0)^I\ne 0$. Therefore
  \begin{align*}
    \pi(\exp(a_n)wv_0) & = \pi \left(\sum_{\lambda\in\Phi_\iota}
      \exp(\lambda(a_n))(w\cdot v_0)^\lambda
    \right)\\
    & = \pi \left(\sum_{\lambda\in\Phi_\iota} \exp \left(
        -\sum_{\alpha\in\Delta_\sigma} n_\alpha(\lambda)\alpha(a_n)      \right) (wv_0)^\lambda
    \right) \\
    & \stackrel{n\to\infty}{\longrightarrow} \pi \left(\sum_{\lambda:
        \, \supp(\lambda)\subset I} \exp \left(
        -\sum_{\alpha\in\Delta_\sigma}n_\alpha(\lambda)\alpha(a)
      \right) (w v_0)^\lambda
    \right)\\
    & = \pi(\exp(a)(wv_0)^I).
  \end{align*}
  This shows that $V^\infty_{J,w}\subset
  \pi(K\exp(\mathfrak{a}^{I,+})(wv_0)^I)$. On the other hand, given
  $a\in\mathfrak{a}^{I,+}$, one can find a sequence $\{a_n\}\subset
  \mathfrak{a}^{+}$ such that $\alpha(a_n)=\alpha(a)$ for $\alpha\in
  I$, $\alpha(a_n)$ is bounded for $\alpha\in J\setminus I$, and
  $\alpha(a_n)\to +\infty$ for $\alpha\in\Delta_\sigma\setminus J$.
  This completes the proof of the first equality.

  By Proposition~\ref{p:satake}, $(wv_0)^J=(wv_0)^I$, and using that
  $\la{a}^{J,+}\subset \la{a}_I+\la{a}^{I,+}$, we deduce that
  \begin{align*}
    \exp(\la{a}^{J,+})(wv_0)^I\subset \R^+\cdot
    \exp(\la{a}^{I,+})(wv_0)^I.
  \end{align*}
  This implies the second equality.

  The third equality is a consequence of the first two equalities.
\end{proof}

Using the same argument as in the proof Proposition~\ref{p:v_inf}, we
also deduce

\begin{Prop}\label{p:closure}
  For every $w\in\mathcal{W}$ and $J\subset \Delta_\sigma$,
  \begin{align*}
    \overline{V^\infty_{J,w}}=\bigcup_{\text{$\lambda_\iota$-connected
        $I\subset J$}} V^\infty_{I,w}.
  \end{align*}
\end{Prop}
Note that Proposition~\ref{p:closure} implies that the orbit
$\mathcal{O}_{I,w}$ is open iff $I\subsetneq \Delta_\sigma$ is a
maximal $\lambda_\iota$-connected set. In this case,
$|I|=|\Delta_\sigma|-1$.

\subsection{Notation and basic facts}

For $g\in G$ let $c_g$ denote the inner conjugation by $g$. For
$w\in\cW$, we define the involutive automorphism $\sigma_w:=c_w\circ
\sigma \circ c_w\inv$. Since $\theta(w)=w$, and $\sigma$ and $\theta$
commute, we have that $\theta$ and $\sigma_w$ also commute. Let
\begin{equation} \label{eq:b_0} \la{b}_0:=\{X\in
  \la{b}:\sigma(X)=X\}=\la{b}\cap \la{h}=\la{b}\cap\la{a}^\perp.
\end{equation}
Since $\Ad w(\la{a})=\la{a}$ and $\Ad w(\la{b})=\la{b}$, and $\Ad w$
preserves the Killing form, we have $\Ad w(\la{b}_0)=\la{b}_0$.
Therefore
\begin{equation}
  \label{eq:wb} 
  \sigma_w(\la{b})=\la{b}, \
  \la{b}_0=\{X\in\la{b}:\sigma_w(X)=X\},\  
  \la{a}=\{X\in\la{b}:\sigma_w(X)=-X\}.
\end{equation}

\subsubsection*{Parabolic subalgebra $\la{p}_J$ and and a
  decomposition of its Levi subalgebra}

Let $J\subset \Delta_\sigma$. Since $\sigma_w(\alpha)=-\alpha$ for all
$\alpha\in\Sigma_\sigma$, we have that $\sigma_w(\la{a}_J)=\la{a}_J$.
Since $\sigma_w$ preserves the Killing form on $\la{g}$, we have that
$\sigma_w(\la{a}^J)=\la{a}^J$. Let $\la{z}_{\la{g}}(\la{a}_J)$ denote
the centralizer of $\la{a}_J$ in $\la{g}$. Let
\begin{equation} \label{eq:Sigma_J}
  \Sigma_J=\{\beta\in\Sigma_\sigma:\beta=\sum_{\alpha\in J} n_\alpha
  \alpha,\, n_\alpha\in\z\} \quad \text{and} \quad
  \Sigma_J^+=\Sigma_J\cap\Sigma_\sigma^+.
\end{equation}
Define
\begin{align}
  \la{n}_J&:=\sum_{\beta\in\Sigma^+_\sigma\setminus
    \Sigma_J}\la{g}_\beta, \quad \text{and} \nonumber\\
  \label{eq:p_J}
  \la{p}_J&:=\la{z}_\la{g}(\la{a}_J)\oplus\la{n}_J,
\end{align}
which is a parabolic subalgebra of $\la{g}$. We define
\begin{align}
  \label{eq:m0}
  \la{m}_0 & =[\la{g}_0,\la{g}_0], \quad \text{where
    $\la{g}_0=\la{z}_{\la{g}}(\la{a})$ as before,} \\
  \la{m}_J & =\sum_{\beta\in\Sigma_J^+} \la{g}_{-\beta} +
  [\la{g}_{-\beta},\la{g}_{\beta}] + \la{g}_{\beta}. \label{eq:m_J}
\end{align}
Then
\begin{align*}
  [\la{z}_{\la{g}}(\la{a}_J),\la{z}_{\la{g}}(\la{a}_J)]=\la{m}_0+\la{m}_J.
\end{align*}
Note that $[\la{m}_0,\la{m}_J]=\la{m}_J$ and
$[\la{m}_J,\la{m}_J]\subset\la{m}_J$. Since
$[\la{z}_{\la{g}}(\la{a}_J),\la{z}_{\la{g}}(\la{a}_J)]$ is semisimple,
its ideal $\la{m}_J$ is semisimple. 

Since $\la{m}_0$ is semisimple for the ideal
\begin{equation} \label{eq:m-lambdai}
\la{m}_{\lambda_\iota}:=\{x\in\la{m}_0:xW^{\lambda_\iota}=0\},
\end{equation}
there exists an ideal $\la{m}_c$ such that
\begin{equation} \label{eq:m_c}
\la{m}_0=\la{m}_c\oplus\la{m}_{\lambda_\iota}.
\end{equation} 
Also, there exist ideals $\la{m}_c^J$ and  $\la{m}_{\lambda_\iota}^J$ of $\la{m}_c$ and $\la{m}_{\lambda_\iota}$ respectively, such that
$$\la{m}_c+\la{m}_J=\la{m}_c^J\oplus\la{m}_J\quad \text{ and } \quad
\la{m}_{\lambda_\iota}+\la{m}_J=\la{m}_{\lambda_\iota}^J\oplus
\la{m}_J.$$

By Remark~\ref{r:b_gamma},
$\la{b}_0=\spn\{b_{\tilde\delta}:\tilde\delta\in\Delta_0\}\subset
{\la{m}}_0$. Since ${\la{m}}_0\subset\la{g}_0$ is semisimple,
we have $\la{a}\cap\la{m}_0=\{0\}$. Therefore
\begin{equation}
  \label{eq:b_0-in-m_0}
  \la{b}_0={\la{m}}_0\cap\la{b}.
\end{equation}

Since $\sigma_w(\la{g}_\beta)=\la{g}_{-\beta}$, we conclude that
$\la{m}_J$ is $\sigma_w$-invariant. Similarly, $\la{m}_J$ is
$\theta$-invariant. Therefore
\begin{align*}
  \la{m}_J=(\la{k}\cap\la{m}_J)\oplus (\la{p}\cap (\Ad w)\la{q}\cap
  \la{m}_J)\oplus ((\Ad w)\la{h}\cap \la{m}_J).
\end{align*}
Note that $\la{a}$ is a maximal abelian subalgebra of $\la{p}\cap(\Ad
w)\la{q}$.  Next we show that 
\begin{align}
  \label{eq:a_J-m_J}
  \la{a}\cap\la{m}_J=\la{a}^J.
\end{align}

For each $\beta\in J$, let $a_\beta\in\la{a}$ be such that
$\inpr{a_\beta}{a}=\beta(a)$ for all $a\in\la{a}$. Then
$\{a_\beta:\beta\in J\}$ is a basis of $\la{a}^J$. Hence by
Remark~\ref{r:a_alpha} and \eqref{eq:m_J}, we have $\la{a}^J\subset
\la{m}_J$.  For any $\beta\in\Sigma_J^+$, if $Y_\pm\in
\la{g}_{\pm\beta}$, and $X\in\la{a}_J$, then
\begin{align*}
  \inpr{X}{[Y_-,Y_+]} = \inpr{[X,Y_-]}{Y_+} =-\beta(X)\inpr{Y_-}{Y_+}
  =0.
\end{align*}
Therefore, by \eqref{eq:m_J}, $(\la{a}\cap\la{m}_J)\perp \la{a}_J$;
that is, $\la{a}\cap\la{m}_J\subset \la{a}^J$. This justifies
\eqref{eq:a_J-m_J}. 

Note that center of $\la{z}_\la{g}(\la{a}_J)$ is contained in the
center of $\la{z}_\la{g}(\la{b})$, which in turn is contained in
$\la{c}=(\la{k}\cap\la{c})+\la{b}$. As $\la{b}=\la{a}+\la{b}_0$, for
$\la{c}_J:=\text{Center}(\la{z}_\la{g}(\la{a}_J))\cap\la{k}$,
\begin{equation} \label{eq:z_J1} 
\la{z}_\la{g}(\la{a}_J) =\la{c}_J+\la{b}+\la{m}_0+\la{m}_J =
\la{c}_J\oplus \la{a}_J \oplus \la{m}_c^J\oplus
\la{m}_{\lambda_\iota}^J\oplus \la{m}_J.
\end{equation}

Let $P_J$ denote the parabolic subgroup of $G$ associated to
$\la{p}_J$. Let $M_J$ and $N_J$ denote the analytic subgroups of $P_J$
associated to the subalgebras $\la{m}_J$ and $\la{n}_J$, respectively.
In the course of the above discussion, we have also proved the
following:

\begin{Prop}\label{p:MI-sym}
  Let $J\subset \Delta_\sigma$ and $w\in\mathcal{W}$. Then the
  semisimple group $M_J$ is invariant under $\theta$ and $\sigma_w$,
  and $\la{a}^{J}$ is a Cartan subalgebra of $M_J$ for the pair
  $(\theta, \sigma_w)$. 
\end{Prop}

Note that $\la{a}^{J}$ has the system of simple roots $J$ with
  the Weyl group $\mathcal{W}_J$, and we have the decomposition
\begin{equation}\label{eq:ccar}
  M_J=(M_J\cap K)\exp(\mathfrak{a}^{J,+})\mathcal{W}_J (wHw^{-1}\cap M_J).
\end{equation}

\begin{lem}\label{l:triv}
  Let $J\subset \Delta_\sigma$. Then the following assertions hold:
  \begin{enumerate}
  \item[(i)] $W^J=\{v\in W: av=\lambda_\iota(a)v,\,\forall a\in
    \la{a}_J\}$.
  \item[(ii)] $Z_G(A_J)\cdot W^J\subset W^J$.
  \item[(iii)] $(Z_G(A_J)\cap wHw^{-1})(wv_0)^J=(w v_0)^J$ for any
    $w\in \mathcal W$.
  \item[(iv)] $N_J$ acts trivially on $W^J$.
 \end{enumerate}
\end{lem}

\begin{proof}
  If $\lambda\in\Phi_\iota$ such that $\supp\lambda\subset J$ and
  $a\in \la{a}_J$, then by definition of $\supp\lambda$, we have that
  $\lambda(a)=\lambda_\iota(a)$. Therefore $av=\lambda_\iota(a)v$ for
  all $v\in W^\lambda$.

  Since an open subset of $\la{a}_J$ is contained in the boundary of
  $\la{a}^+$, there exits $X\in\la{a}_J$ such that $\alpha(X)>0$ for
  all $\alpha\in\Delta_\sigma\setminus J$. Therefore if
  $\lambda\in\Phi_\iota$ such that $\supp \lambda\not\subset J$, then
  $\lambda(X)<\lambda_\iota(X)$. Therefore (i) follows from the above
  two observations and the definition of $W^J$.

  Since the centralizer preserves the isotypical components, we obtain
  (ii).

  Let $X\in\la{a}_J$ be as above. Then
  \begin{equation} \label{eq:limJ} \lim_{t\to\infty}
    e^{-t\lambda_\iota(X)}\exp(tX)(wv_0)=(wv_0)^J.
  \end{equation}
  
  Since $Z_G(A_J)\cap wHw^{-1}$ fixes $wv_0$, it acts trivially on the
  $\R$-span of $A_J(wv_0)$, which contains $(wv_0)^J$ by
  \eqref{eq:limJ}. Therefore (iii) holds.

  Let $\lambda\in \Phi_\iota$ be such that $\supp\lambda\subset J$;
  that is, $\lambda=\lambda_\iota-\sum_{\alpha\in J} n_\alpha\alpha$,
  where all $n_\alpha\geq 0$. Suppose that $\gamma=\sum_{\alpha\in
    \Delta_\sigma} m_\alpha \alpha$, where all $m_\alpha\geq 0$, is
  such that $\lambda+\gamma\in \Phi_\iota$. Since $\lambda_\iota$ is
  the highest weight in $\Phi_\iota$, we have
\begin{equation*}
  \lambda_\iota - (\lambda+\gamma)\in \Sigma_\sigma^+.
\end{equation*}
Then $m_\alpha=0$ for all $\alpha\in \Delta_\sigma\setminus J$; that
is, $\gamma\in \Sigma_J^+$.  This shows that if $\gamma\in
\Sigma_\sigma^+\setminus\Sigma_J$, then
$$\la{g}_\gamma
W^{\lambda}\subset W^{\lambda+\gamma}=0 .$$ Therefore $\la{n}_JW^J=0$.
Thus (iv) holds.

\end{proof}

\begin{Prop}\label{p:irr}
  $\la{z}_J:=\la{z}_{\la{g}}(\la{a}_J)$ acts irreducibly on $W^J$.
\end{Prop}

\begin{proof}
With notation as in the proof of Lemma \ref{l:irr}, we have
$U_0(\la{n})W^{\lambda_\iota}=0$ and  $U_0(\la{n}^-)W^{\lambda_\iota}\subset W'$.
This easily implies that $\la{g}_0$ acts irreducibly on $W^{\lambda_\iota}$. 

Since $\la{n}=(\la{m}_J\cap\la{n})+\la{n}_J$, and $\la{n}_J\cdot
  W^J=0$, it follows from Lemma~\ref{l:irr}
  \begin{align}
    \label{eq:m_I}
    W^{\lambda_\iota}=\{v\in W^J:(\la{m}_J\cap\la{n})v=0\}.
  \end{align}
  Since $\la{z}_J$ is reductive, $W^J$ is a direct sum of irreducible
  $\la{z}_J$-modules, and each of them admits a nonzero subspace
  which is annihilated by $\la{m}_J\cap\la{n}$ (Engel's
  theorem). Hence by \eqref{eq:m_I} each of the $\la{z}_J$-submodules
  contains a nonzero subspace of $W^{\lambda_\iota}$ which in turn is
  invariant under $\la{g}_0\subset \la{z}_J$. Since
  $\la{g}_0$ acts irreducibly on $W^{\lambda_\iota}$, we
  conclude that $\la{z}_J$ acts irreducibly on $W^J$.
\end{proof}

\begin{lem} \label{eq:m_c-k}
$\la{m}_c\subset\la{k}\cap\la{z}_{\la{g}}(\la{b})$. 
\end{lem}

\begin{proof}
Let $\tilde{\la{m}}_0$ be maximal noncompact ideal of $\la{m}_0$.
%Define 
%\begin{equation} \label{eq:m0_tilde} 
%\tilde{\la{m}}_0 = \sum_{\tilde\gamma\in \langle\Delta_0\rangle\cap\Sigma^+}
%  \la{g}_{-\tilde\gamma}(\la{b}) +
%  [\la{g}_{-\tilde\gamma}(\la{b}),\la{g}_{\tilde\gamma}(\la{b})] +
%  \la{g}_{\tilde\gamma}(\la{b}).
%\end{equation}
It follows from \cite[Lemma~7.1.4]{sch} that $\tilde{\la{m}}_0\subset
(\Ad w)(\la{h})$. Hence, by Lemma~\ref{l:triv}(iii), 
$\tilde{\la{m}}_0\subset \la{m}_{\lambda_\iota}$. 
Since by Remark~\ref{r:b_gamma},
$\la{b}_0=\spn\{b_{\tilde\delta}:\tilde\delta\in\Delta_0\}\subset
\tilde{\la{m}}_0$, we deduce that
$\la{m}_c\subset\la{z}_{\la{g}}(\la{b})$. 
This implies the lemma.
\end{proof}

\begin{Def} \label{d:J(I)} \rm For a $\lambda_\iota$-connected $I$ of
  $\Delta_\sigma$, let
  $I'=(I\cup\{\lambda_\iota\})^\perp\cap\Delta_\sigma$; that is, the
  set of roots in $\Delta_\sigma$ which are orthogonal to
  $\lambda_\iota$ and all roots in $I$. We define
$$J(I):=I\cup I' .$$  We
note the following:
\begin{enumerate}
\item $J(I)$ uniquely determines $I$, as $I$ is the maximal
  $\lambda_\iota$-connected subset of $J(I)$.

\item If $\beta,\gamma\in \Sigma_\sigma^+$ and $\beta\perp\gamma$,
  then $\beta+\gamma\not\in \Sigma^+_\sigma$. Therefore
  $[\la{m}_I,\la{m}_{I'}]=0$, and hence $M_{J(I)}=M_IM_{I'}$ is an
  almost direct product. Also
  $\Sigma_{J(I)}=\Sigma_{I}\cup\Sigma_{I'}$.

\item $W^{J(I)}=W^I$.
\end{enumerate}
\end{Def}

\begin{Prop} \label{p:irr2} Let $I$ and $J=J(I)$ be as above.  Then
$$
(\la{m}_{\lambda_\iota}^J+\la{m}_{I'})W^{I}=0
\quad\hbox{and}\quad cw=\Lambda_\iota(c)w,\quad\forall c\in\la{c}_J,\;
w\in W^I.
$$
%$\la{c}_J$ acts trivially on $W^J$ if $\K=\R$, and it
%  acts on $W^J$ via the linear functional
%  $\Lambda_\iota:\la{c}_J\mapsto i\R$ if $\K=\C$; recall
%  Remark~\ref{r:K}.
The Lie algebra $\la{m}_c^I\oplus \la{m}_I$ acts irreducibly and faithfully on $W^I$
  over $\K$.
\end{Prop}

\begin{proof}
 By the definition of
  $I'$ and Proposition~\ref{p:satake}, for any $\gamma\in \Sigma_{I'}$
  and $\lambda\in\Phi_\iota$ with $\supp\lambda\subset I$, we have
  $\lambda+\gamma\not\in \Phi_\iota$.  Therefore,
  $\la{m}_{I'}W^J=0$.
Since $\la{m}_{\lambda_\iota}^J$ is a semisimple ideal in
$\la{z}_{\la{g}}(\la{a}_J)$, by Proposition~\ref{p:irr} and
\eqref{eq:m-lambdai}, we conclude that
$\la{m}_{\lambda_\iota}^JW^I=0$. 
  By Proposition~\ref{p:irr} and \eqref{eq:z_J1}, $\la{c}_J\oplus \la{m}_c^J\oplus
  \la{m}_I$ acts irreducibly on $W^J$ over $\K$. If $\K=\R$ then
  $W_{\C}^{\Lambda_\iota}\cap W$ is a one-dimensional
  $\la{c}_J$-invariant subspace. Since $\exp(\la{c}_J)\subset K$, we
  conclude that $\la{c}_J(W_{\C}^{\Lambda_\iota}\cap W)=0$. Hence, by
  irreducibility, $\la{c}_JW^J=0$.
  Suppose if $\K=\C$, then $\la{c}_J$ being central in
  $\la{c}_J\oplus \la{m}_c^J\oplus \la{m}_I$, by the irreducibility we
  conclude that $\la{c}_J$ acts via $\K$-scalars on $W^J$. This
  proves the first claim.

  It follows from above that $\la{m}_c^J\oplus \la{m}_I$ act irreducibly on $W^J$ over $\K$. Since
  $\la{m}_0\cap\la{m}_{I'}\subset\la{m}_{\lambda_\iota}$, we have that
  $\la{m}_c^J=\la{m}_c^I$. This proves irreducibility.
 
  By \eqref{eq:m-lambdai} and \eqref{eq:m_c}, $\la{m}_c^I$ acts
  faithfully on $W^{\lambda_\iota}$, and hence on $W^I$.
  We observe that any nonzero ideal of $\la{m}_I$ contains $a^{I_1}$
  for some $\emptyset\neq I_1\subset I$ such that $I_1\perp(I\setminus
  I_1)$.  Since $I$ is $\lambda_\iota$-connected, $I_1$ is
  $\lambda_\iota$-connected, and hence $\la{a}^{I_1}\not\subset\ker
  (\lambda_\iota)$. But then
  $\la{a}^{I_1}v_0^{\lambda_\iota}=\lambda_\iota(\la{a}^{I_1})v_0^{\lambda_\iota}\neq
  0$. Thus, $\la{m}_I$ acts faithfully on $W^I$.
\end{proof}

%\subsubsection*{Description of $G$-orbits on the boundary}

\begin{Prop}\label{p:orbit}
  For a $\lambda_\iota$-connected subset $I$ of $\Delta_\sigma$ and
  $w\in \mathcal{W}$, we have
  \begin{align*}
    \pi(G(wv_0)^I)=\pi(K\exp(\mathfrak{a}^{I,+})\mathcal{W}_I
    (wv_0)^I) =\pi(K\exp(\mathfrak{a}^{I,+})(\mathcal{W}_Iw v_0)^I).
  \end{align*}
\end{Prop}

\begin{proof}
  By the Iwasawa decomposition we have
  $G=KP_I=K\exp(\la{z}_{\la{g}}(\la{a}_I))(\la{a}_I)N_I$. Now $N_I$
  acts trivially on $W^I$. Therefore, in view of
  Proposition~\ref{p:irr2},
      $$\pi(G(wv_0)^I)=\pi(KM_I(wv_0)^I) .$$ Now the first equality follows
      from Lemma~\ref{l:triv}(iii), and \eqref{eq:ccar}.

  Since the weight spaces $W^\lambda$, $\lambda\in\Phi_\iota$, are
  orthogonal with respect to a $K$-invariant scalar product, and by
  Lemma~\ref{l:triv}, $\mathcal{W}_I\subset M_I\cap K$ preserves
  $W^I$, it follows that $\mathcal{W}_I$ preserves the orthogonal
  complement of $W^I$, and hence
  \begin{align*}
    w\cdot v^I=(w v)^I\quad \text{for all $w\in\mathcal{W}_I$ and
      $v\in W$}.
  \end{align*}
  This justifies the second equality in the proposition.
\end{proof}

\subsection{Disjointness of the $G$-orbits in the boundary and stabilizers of $(w_0v_0)^I$}

\begin{lem}\label{l:stab}
  Let $I$ be a $\lambda_\iota$-connected subset of $\Delta_\sigma$ and
  let $J=J(I)$. Then
  \begin{align*}
    \Stab_G(W^I):=\{g\in G:gW^I=W^I\}=P_J.
  \end{align*}
\end{lem}

\begin{proof}
  Let $Q=\Stab_G(W^I)$. It follows from Lemma~\ref{l:triv} that
  $Q\supset P_J$. Hence, $Q=P_S$ with $J\subset S\subset
  \Delta_\sigma$. Since $\la{z}_S:=\la{z}_{\la{g}}(\la{a}_S)\subset
  P_S$, we have $\la{z}_SW^I=W^I$. By Proposition~\ref{p:irr}
  $\la{z}_S$ acts irreducibly on $W^S$. Therefore $W^I=W^S$. Since $I$
  is $\lambda_\iota$-connected, by Proposition~\ref{p:satake} $I$ is
  the maximal $\lambda_\iota$-connected component of $S$. Hence by
  Definition~\ref{d:J(I)} $S\subset J(I)$. Hence $J=S$.
\end{proof}

Let $I$ be a $\lambda_\iota$-connected subset of $\Delta_\sigma$,
 $J=J(I)$, and $w_0\in \cW$. We consider the group $L=\{g\in
G:g(w_0v_0)^I=(w_0v_0)^I\}$ with the Lie algebra
$\la{l}=\{X\in\la{g}: X(w_0v_0)^I=0\}$.
Note that $\la{a}_J$ normalizes $\la{l}$, because $\la{a}_J\subset
\la{a}_I$, and by Lemma~\ref{l:triv} $(w_0v_0)^I$ is an eigenvector
for each element of $\la{a}_I$. Moreover by Lemma~\ref{l:triv}, we
have $\la{n}_J\subset \la{l}$. Therefore
\begin{equation} \label{eq:l}
  \la{l}=\la{n}_J+\la{l}\cap\la{z}_\la{g}(\la{a}_J)+
  \la{l}\cap \sum_{\beta\in\Sigma^+_\sigma\setminus\Sigma_J}
  \la{g}_{-\beta}.
\end{equation}
By \eqref{eq:z_J1} and Proposition~\ref{p:irr2}
\begin{equation} \label{eq:l_zJ} \la{l}\cap\la{z}_\la{g}(\la{a}_J) =
  (\la{c}_J\cap\ker\Lambda_\iota)\oplus(\la{a}_J\cap\ker\lambda_\iota)\oplus
  \la{m}_{\lambda_\iota}^J\oplus\la{m}_{I'}\oplus((\la{m}_c^I\oplus\la{m}_I)\cap\la{l}).
\end{equation}

\begin{Prop} \label{p:stab_M_I} We have
  \begin{align*}
    (\la{m}_c^I\oplus\la{m}_I)\cap\la{l} =(\la{m}_c^I\oplus\la{m}_I)\cap (\Ad w_0)\la{h}= 
(\la{m}_c^I\cap(\Ad w_0)\la{h}) \oplus (\la{m}_I\cap(\Ad w_0)\la{h}).
  \end{align*}
  In particular, the orthogonal projection of $\la{l}$ on $\la{a}^I$
  is trivial.
\end{Prop}

\begin{proof} 
  By Proposition~\ref{p:MI-sym}, applied to $I$ in place of $J$ and
  $w_0$ in place of $w$, $\la{h}_{I,w_0}:=(\la{m}_c^I+\la{m}_I)\cap\Ad
  w_0(\la{h})$ is a symmetric subalgebra of $\la{m}_c^I+\la{m}_I$. By
  Lemma~\ref{l:triv}(iii), $\la{h}_{I,w_0}(w_0v_0)^I=0$. Let $M_c^J$
  denote the analytic subgroup of $G$ associated to $\la{m}_c^I$. Due
  to Proposition~\ref{p:irr2}, we can apply
  Corollary~\ref{cor:gen:finite} to $M_c^IM_I$ in place of $G$, and
  $W^I$ in place of $E$, to obtain the first equality.
The second equality holds because $\la{m}_c^I$ and $\la{m}_I$ are
  invariant under $\sigma$. The last conclusion follows from
  Proposition~\ref{p:MI-sym} because $\la{a}^I$ is orthogonal to $\Ad
  w_0(\la{h})$.
\end{proof}

\begin{Prop}
  \label{p:cofinite}
  $\la{l}\subset\la{p}_J.$
\end{Prop}

\begin{proof}
  Suppose that $\la{l}\not\subset\la{p}_J$. Then by \eqref{eq:l} there
  exists $\beta\in\Sigma_\sigma^+\setminus \Sigma_J$, $\tilde
  \beta\in \Sigma^+$ with $\tilde\beta|_{\la{a}}=\beta$, and
 $0\neq X\in \la{l}$ such that
$$
X=X_{-\tilde\beta}+Y\quad\hbox{where $0\neq X_{-\tilde\beta}\in
\la{g}_{-\tilde\beta}(\la{b})$ and $Y\in
\sum_{\tilde\gamma\in\Sigma\backslash\{0\}:\tilde\gamma|_{\la{a}_J}=0}
\la{g}_{-\tilde\beta+\tilde\gamma}$.}
$$
Replacing $X$ by a scalar
multiple from the beginning, without loss of generality we may assume
that $B(X_{-\tilde\beta},X_{-\tilde\beta})=1$ (see~\eqref{eq:BB}).
Since $-\theta(X_{\tilde\beta}) \in \la{g}_{\tilde\beta}(\la{b})
\subset \la{g}_\beta\subset\la{l}$, we have
$[-\theta(X_{-\tilde\beta}),X]=b_{\tilde\beta}+Z\in\la{l}$, where
\begin{align*}
  b_{\tilde\beta}=[-\theta(X_{-\tilde\beta}),X_{-\tilde\beta}]\in\la{b}\quad\hbox{and}\quad
  Z=[-\theta(X_{-\tilde\beta}),Y]\in
  \sum_{\tilde\gamma\in\Sigma\backslash \{0\}}
  \la{g}_{\tilde \gamma}(\la{b}) \cap \la{z}_{\la{g}}(\la{a}_J).
\end{align*}
Therefore $b_{\tilde\beta}+Z\in \la{l}\cap\la{p}_J$ and its projection
on $\la{b}$ equals $b_{\tilde\beta}$. By \eqref{eq:l}, \eqref{eq:l_zJ}
and Proposition~\ref{p:stab_M_I},
\begin{align*}
  b_{\tilde\beta}\in \la{b}_0+(\la{a}_I\cap\ker\lambda_\iota).
\end{align*}
If $a_{\tilde\beta}$ is the projection of $b_{\tilde\beta}$ on
$\la{a}$, then
\begin{align*}
  a_{\tilde\beta}\in \la{a}_I \cap\ker\lambda_\iota.
\end{align*}
Since $\la{b}_0\perp\la{a}$, by Remark~\ref{r:b_gamma}, for all
$a\in\la{a}$,
\begin{align*}
  \inpr{a}{a_{\tilde \beta}} = \inpr{a}{b_{\tilde \beta}} =
  \tilde{\beta}(a)=\beta(a).
\end{align*}
Therefore, $\beta\perp (I\cup\{\lambda_\iota\})$. Using that scalar
products of simple roots are nonpositive, we deduce that $\beta\in
\langle I'\rangle \subset J$, which is a contradiction.
\end{proof}

From \eqref{eq:l}, \eqref{eq:l_zJ}, Proposition~\ref{p:stab_M_I} and
Proposition~\ref{p:cofinite} we deduce the following:

\begin{Cor} \label{c:l} We have
  \begin{align*}
    \la{l} = (\la{c}_J\cap\ker\Lambda_\iota) \oplus
    (\la{a}_J\cap\ker\lambda_\iota) \oplus \la{m}_{\lambda_\iota}^J
    \oplus \la{m}_{I'} \oplus (\la{m}_c^J\cap\Ad w_0(\la{h}))\oplus
    (\la{m}_I\cap\Ad w_0(\la{h})) \oplus \la{n}_J.
  \end{align*}
  In particular,
  \begin{align*}
    \operatorname{Unipotent\ Radical}(L)=N_J \quad\text{and}\quad
    L\subset P_J.
  \end{align*}
  \qed
\end{Cor}

\begin{proof}[Proof of Theorem~\ref{thm:satake}]
  Consider a divergent sequence
  \begin{align*}
    v_n=k_n\exp(a_n)w_nv_0\in V.
  \end{align*}
  where $k_n\in K$, $a_n\in \mathfrak{a}^{+}$, $w_n\in \mathcal{W}$.
  After passing to a subsequence if necessary, we can assume that
  $k_n\to k\in K$, $w_n=w\in \mathcal{W}$, and that there exists
  $J\subsetneq \Delta_\sigma$ such that $\alpha(a_n)$ is bounded for
  $\alpha\in J$ and $\alpha(a_n)\to \infty$ for $\alpha\in
  \Delta_\sigma \setminus J$. Then the limit points of the sequence
  $v_n$ are in $V^\infty_{J,w}$. This proves that
  \begin{align*}
    V^\infty=\bigcup_{J\subsetneq \Delta_\sigma, w\in\mathcal{W}}
    V^\infty_{J,w}.
  \end{align*}
  Moreover, by Proposition~\ref{p:v_inf}, it suffices to take the
  union over $\lambda_\iota$-connected subsets only. Then
  \eqref{eq:oiw} follows from Proposition~\ref{p:orbit}.
  
  By Corollary~\ref{c:l} the unipotent radical of the stabilizer of
  $g_i(w_iv_0)^{I_i}$ is $g_iN_{J_i}g_i\inv$.  Therefore if the
  $G$-orbits of $(w_iv_0)^{I_i}$ are same, then
  $g_1N_{J_1}g_1\inv=g_2N_{J_2}g_2\inv$. Since $P_{J_i}=N_G(N_{J_i})$,
  we have $g_1P_{J_1}g_1\inv=g_2P_{J_2}g_2\inv$. Hence $J_1=J_2$ and
  $I_1=I_2$. Thus (\ref{eq:oiw2}) follows.
\end{proof}

Note that Theorem~\ref{th:sat} follows from Theorem~\ref{thm:satake},
Proposition~\ref{p:v_inf} and Proposition~\ref{p:closure}.

Recall that for $\Omega\subset S(W)$, we have defined
\begin{align*}
  \Theta_{\Omega}=\{I\subset \Delta_\sigma:\text{$I$ is
    $\lambda_\iota$-connected and } \bar \Omega\cap V_I^\infty\ne
  \emptyset\}.
\end{align*}
We denote by $\Theta_{\Omega}^{\min}$ the set of minimal elements in
$\Theta_\Omega$ with respect to inclusion.

\begin{cor}[Tube lemma]\label{c:satake}\label{p:nbhd}
  For any compact set $\Omega\subset S(W)$, there exists a collection
  $\{U_J\subset \la{a}^{J,+}:\, J\in\Theta_\Omega^{\min}\}$ of compact
  sets such that
  \begin{align*}
    \Omega\cap V_I^\infty &\subset \bigcup_{J\in\Theta_\Omega^{\min}}
    \pi (K\exp(U_J+\la{a}_J^+)(\mathcal{W}v_0)^I) \quad \text{for
      every $I\subseteq \Delta_\sigma$.}
  \end{align*}
\end{cor}

\begin{proof}
  Suppose that the corollary fails. Then for any choice of compact
  sets $U_J$, there exists $v=\pi(k\exp(a)(wv_0)^I)\in \Omega$ with
  $k\in K$, $w\in W$, and $a\in\la{a}^+$, $a\notin U_J+\la{a}_J^+$ for
  every $J\in\Theta_\Omega^{\min}$.  Therefore, there exists
  $v_n=\pi(k_n\exp(a_n)(w_nv_0)^I)\in \Omega$ with $k_n\in K$,
  $a_n\in\la{a}^+$, $w_n\in W$ such that for every
  $J\in\Theta_\Omega^{\min}$, $\alpha(a_n)\to\infty$ for at least one
  $\alpha$ in $J$. Passing to a subsequence, we may assume that for
  some $I\subset \Delta_\sigma$, $\alpha(a_n)\to\infty$ if
  $\alpha\in\Delta_\sigma\setminus I$ and $\alpha(a_n)$ is bounded if
  $\alpha\in I$. Then by Proposition~\ref{p:v_inf}, the limit points
  of the sequence $\{v_n\}$ are in $V^\infty_{I_0,w_0}$,
  $w_0\in\mathcal{W}$, where $I_0$ is the largest
  $\lambda_\iota$-connected subset of $I$. Then $\Omega\cap
  V^\infty_{I_0}\ne \emptyset$ and $I_0\in \Theta_\Omega$.  On the
  other hand, $I_0\nsupseteq J$ for every $J\in\Theta_\Omega^{\min}$.
  This gives a contradiction and proves the corollary.
\end{proof}

\section{Invariant measures at infinity}
\label{sec:invariant-measures}

In the previous section, we have shown that
\begin{equation*}
  V^\infty=\bigsqcup_{\text{$\lambda_\iota$-connected $I\subsetneq
      \Delta_\sigma$}} V^\infty_I
\end{equation*}
where
\begin{equation*}
  \R^+\cdot V^\infty_I=\bigcup_{w\in\mathcal{W}}
  G(wv_0)^I.
\end{equation*}
In this section we describe an algebraic condition on $I$ so that
$V_I^\infty$ admits a $G$-invariant measure, and give a formula for
the measure. We also provide a natural class of $I$ for which the
condition holds. The results of this section are obtained mainly for
the the sake of more complete description of the boundary. They are
not essential for the proofs of the main results stated in the
introduction.

\begin{thm}
  \label{thm:measure}

  Let $I$ be a $\lambda_\iota$-connected subset of $\Delta_\sigma$.
  Then for any $w_0\in \cW$, there exists a $G$-invariant measure on
  $G(w_0v_0)^I$ if and only if
  \begin{align}
    \label{eq:a_J-lambda-rho}
    \la{a}_J\cap\ker\rho=\la{a}_J\cap\ker\lambda_\iota,
  \end{align}
  where $J=J(I)$ is as in Definition~\ref{d:J(I)}.
  
  If (\ref{eq:a_J-lambda-rho}) hold, then the $G$-invariant measure on
  $G(w_0v_0)^I$, say $\nu_{I,w_0}$, is given by (up to a constant
  multiple)
  \begin{align}
    \label{eq:measure}
    &\int\limits_{W} f\, d\nu_{I,w_0} \\
    \nonumber & = \int\limits_{K} dk \int\limits_{a\in\la{a}^{I,+}} da
    \int\limits_{\bar b\in \la{a}_J/\la{a}_J\cap\ker\lambda_\iota}
    \sum_{w\in\mathcal{W}_I} f \big(k\exp(a+b)w(w_0v_0)^{I}
    \big)\delta_I(a) e^{2\rho(b)}\,{d\bar b}
  \end{align}
  for all $f\in C_c(G(w_0v_0)^I)$, where $dk$, $da$ and $d\bar b$
  denote the Haar integrals on $K$, $\la{a}^{I}$, and
  $\la{a}_J/\la{a}_J\cap\ker\rho$, respectively, and
  \begin{align}
    \label{eq:delta_I}
    \delta_I(a):=\prod_{\alpha\in \Sigma_I^+}
    (\sinh\alpha(a))^{l_\alpha^+}(\cosh\alpha(a))^{l_\alpha^-},
    \forall a\in\la{a}^I.
  \end{align}
\end{thm}

\begin{proof}
  Let $w_0\in\cW$. Since $G$ admits no nontrivial positive real
  characters, there exists a $G$-invariant measure on $G\cdot
  (w_0v_0)^I$ {\em if and only if\/} $L:=\Stab_G((w_0v_0)^I)$ is
  unimodular. By Corollary~\ref{c:l}, $N_J$ is the unipotent radical
  of $L$. Therefore $L$ is unimodular {\em if and only if\/}
  $\abs{\det(\Ad g|_{\la{n}_J})}=1$ for all $g\in L$, {\em if and only
    if\/} $\tr(\ad(x)|_{\la{n}_J})=0$ for all $x\in\la{l}$.

  Note that $\la{c}_J\subset\la{k}$, and $\la{m}_c^J\oplus\la{m}_{\lambda_\iota}^J\oplus \la{m}_{J}$
  is semisimple. Also each of them normalizes $\la{n}_J$. Therefore
  $$\tr((\ad x)|_{\la{n}_J})=0\quad\text{ for all }x\in
  \la{c}_J+\la{m}_c^J+\la{m}_{\lambda_\iota}^J+\la{m}_J .$$ Therefore
  by Corollary~\ref{c:l}, $L$ is unimodular {\em if and only if\/}
  \begin{equation} \label{eq:rho}
    2\rho(b)=\sum_{\alpha\in\Sigma_\sigma^+\setminus\langle J\rangle}
    \tr(\ad(b)|_{\la{g}_\alpha})=0,\quad\forall
    b\in\la{a}_J\cap\ker\lambda_\iota.
  \end{equation}
  This equation is equivalent to $\ker\lambda_\iota\cap
  \la{a}_J\subset \ker\rho\cap\la{a}_J$. If $\la{a}_J\neq\{0\}$, then
  $\la{a}_J\not\subset\ker\rho$, and hence $\ker\rho\cap\la{a}_J$ is
  of codimension $1$ in $\la{a}_J$. Also
  $\ker\lambda_\iota\cap\la{a}_J$ is of codimension at most one in
  $\la{a}_J$. Therefore \eqref{eq:rho} is equivalent to
  \eqref{eq:a_J-lambda-rho}. This proves the first part of the
  theorem.

  Now to obtain the formula for the Haar integral on $G$, we suppose
  that \eqref{eq:a_J-lambda-rho} holds.  In view of \eqref{eq:p_J} and
  \eqref{eq:z_J1}, let $\tilde M_J$ be the closed subgroup of $P_J$
  associated to the Lie subalgebra
  $\la{c}_J+\la{m}_c^J+\la{m}_{\lambda_\iota}^J+\la{m}_{J}$ such that
  $P_J=\tilde M_JA_JN_J$, where $Z_G(A_J)=\tilde M_JA_J$ is a direct
  product. A right Haar integral on $P_J$ can be given by
  \begin{align}
    \label{eq:P_J:measure}
    f\mapsto \int_{\tilde M_J}dm \int_{\la{a}_J}db \int_{N_J}
    f(m\exp(b)n)e^{2\rho(b)}\, dn, \quad \forall\,f\in C_c(P_J),
  \end{align}
  where $dm$ and $dn$ denote Haar integrals on $\tilde M_J$ and $N_J$,
  respectively, and $db$ denotes the Lebesgue integral on $\la{a}_J$.
 
  Note that
  \begin{equation} \label{eq:tilde-MJ} 
    \tilde M_J=(K\cap P_J)\tilde M_J^0=(K\cap P_J)M_{\lambda_\iota}^JM_J,
  \end{equation}
  where $M_{\lambda_\iota}^J$ is the analytic subgroup of $G$ associated
  to the subalgebra $\la{m}_{\lambda_\iota}^J$. By Corollary~\ref{c:l},
  \begin{align}
    \label{eq:M_J-M_I}
    (M_{\lambda_\iota}^JM_J)/(M_{\lambda_\iota}^JM_J)\cap L\cong
    M_I/M_I\cap w_0Hw_0\inv.
  \end{align}
  By \eqref{eq:ccar} and \eqref{eq:volume} a left invariant integral
  on $M_I/M_I\cap w_0Hw_0\inv$ is given by
  \begin{align}
    \label{eq:M_I:measure}
    f\mapsto \int\limits_{K\cap M_I}dk \int\limits_{a\in\la{a}^{I,+}}
    \sum_{w\in\mathcal{W}_I} f(kaw(M_I\cap w_0Hw_0\inv)) \delta_I(a)\,
    da,
  \end{align}
  where $f\in C_c(M_I/M_I\cap w_0Hw_0\inv)$, and $dk$ denotes a Haar
  integral on $K\cap M_I$.
  
  Combining \eqref{eq:P_J:measure}, \eqref{eq:tilde-MJ},
  \eqref{eq:M_J-M_I}, and \eqref{eq:M_I:measure} we obtain that for
  all $f\in C_c(P_J)$,
  \begin{align}
    \label{eq:P_J-L}
    f\mapsto \int\limits_{K\cap P} dk
    \int\limits_{a\in\mathfrak{a}^{I,+}} da \int\limits_{\bar b\in
      \mathfrak{a}_J/\mathfrak{a}_J\cap\ker\lambda_\iota} {d\bar b}
    \sum_{w\in\mathcal{W}_I}\int\limits_L f(k\exp(b)\exp(a)wl)
    \delta_I(a)\exp(2\rho(b))\, dl,
  \end{align}
  defines a right Haar integral, say $dp$, on $P_J$, were $dk$ and
  $dl$ denote Haar integrals on $K\cap P_J$ and $L$, respectively.

  Note that a Haar integral on $G$ is given by
  \begin{align}
    \label{eq:KP_J}
    f\mapsto \int_{K} dk \int_{P_J} f(k p)\, {dp},\quad f\in C_c(G),
  \end{align}
  where $dk$ denotes a Haar integral on $K$. Combining
  \eqref{eq:P_J-L} and \eqref{eq:KP_J}, and the fact that $L$ is the
  stabilizer of $(v_0w_0)^I$ in $G$, we obtain that the formula
  \eqref{eq:measure} indeed gives a $G$-invariant measure on
  $G(w_0v_0)^I$.
\end{proof}

It turns out that the sets $I_\iota(I)$ satisfy the condition of the
above theorem; see \eqref{eq:Iiota} for the definition. To show this
we need the following:

\begin{Prop}
  \label{p:I_iota}
  If $I\subset\Delta_\sigma$ is $\lambda_\iota$-connected, then
  $I_\iota(I)$ is $\lambda_\iota$-connected. In fact, any $J\subset
  \Delta_\sigma$ containing $I_\iota(I)$ is $\lambda_\iota$-connected.
\end{Prop}

\begin{proof}
  Let $I_0\supset I$ be the largest $\lambda_\iota$-connected subset
  of $J$, and let ${S}=J\setminus {I_0}$. Then ${S}\subset
  ({I_0}\cup\{\lambda_\iota\})^\perp$. Therefore $\la{a}^{S}\subset
  \la{a}_{I_0}\cap\ker\lambda_\iota$.  Note that
  $\la{a}_{I_0}=\la{a}_J\oplus \la{a}_{\Delta_\sigma\setminus {S}}$.
  Now given $a\in \la{a}^{S,+}$, we write $a=x+y$ with
  $x\in\la{a}_{\Delta_\sigma\setminus S}$ and $y\in\la{a}_J$. Then for
  any $\alpha\in S\subset J$ we have $\alpha(a)\geq 0$ by
  \eqref{eq:aJplus}, and $\alpha(y)=0$, and hence $\alpha(x)\geq 0$.
  And for any $\beta\in \Delta_\sigma\setminus S$, we have
  $\beta(x)=0$.  Therefore $x\in\la{a}^+$; in other words,
  \begin{equation*}
    \label{eq:S-plus}
    \la{a}^{{S},+}\subset 
    (\la{a}_{\Delta_\sigma\setminus S}\cap \la{a}^+) + \la{a}_J.
  \end{equation*}
  Therefore, since $I\subset I_0\subset \Delta_\sigma\setminus {S}$
  and $I_\iota(I)\subset J$, we get
  \begin{equation}
    \label{eq:S-plus-in}
    \la{a}^{{S},+}\subset (\la{a}_I\cap \la{a}^+) +
    \la{a}_{I_\iota(I)}.
  \end{equation}
  
  By the definition of $I_\iota(I)$ as in \eqref{eq:Iiota}, exists
  $C=a_\iota(I)>0$ such that
  \begin{align}
    \label{eq:C1}
    \rho(x)&\leq C\lambda_\iota(x), \quad \forall x\in \la{a}_I\cap
    \la{a}^+,
    \qquad \text{and}\\
    \label{eq:C2}
    \rho(y)&= C\lambda_\iota(y), \quad \forall y\in
    \la{a}_{I_\iota(I)}.
  \end{align}
  
  Combining \eqref{eq:S-plus-in}, \eqref{eq:C1}, \eqref{eq:C2}, and
  since $\la{a}^S\subset\ker\lambda_\iota$, we conclude that
  \begin{equation}
    \label{eq:rho-0} 
    \rho(a) \leq C\lambda_\iota(a)=0\,\quad\forall\,a\in\la{a}^{{S},+}.
  \end{equation}

  Since $S\perp I_0$, we have $\Sigma_J^+=\Sigma_{I_0}^+\cup
  \Sigma_{S}^+$ (cf.\ Definition~\ref{d:J(I)}). Therefore given
  $a\in\la{a}^{S,+}\subset\la{a}_{I_0}$, we have
  \begin{align}
    \label{eq:S-rho}
    2\rho(a) =\tr(\ad(a)|_{\la{n}_J}) + \sum_{\alpha\in\Sigma^+_S}
    \tr(\ad(a)|_{\la{g}_\alpha}).
  \end{align}
  Recall that $\la{a}^{S}\subset \la{m}_S\subset \la{m}_J$, $\la{m}_J$
  is semisimple, and $[\la{m}_J,\la{n}_J]\subset\la{n}_J$. Therefore
  \begin{align*}
    \tr(\ad(a)|_{\la{n}_J})=0.
  \end{align*}
  Note that
  \begin{align}
    \label{eq:Splus}
    \tr(\ad(a)|_{\la{g}^\alpha}) = (\dim\la{g}_\alpha) \alpha(a)\geq
    0, \quad \forall \alpha\in\Sigma_S^+.
  \end{align}
  Now by \eqref{eq:S-rho}, we get $\rho(a)\geq 0$. Therefore
  \eqref{eq:rho-0}, we get $\rho(a)=0$. Hence by \eqref{eq:Splus},
  \begin{equation} \label{eq:S-alpha}
    (\dim\la{g}_\alpha)\alpha(a)=0,\quad \forall \alpha\in S,\ \forall
    a\in\la{a}^{S,+}.
  \end{equation}
  Now if $S\neq\emptyset$, then for any $\alpha\in S$: we have
  $\dim\la{g}_\alpha\geq 1$; and since $\la{a}^S\perp\la{a}_S$, we
  have $\alpha(a)\neq 0$ for any $0\neq a\in\la{a}^S$. This
  contradicts \eqref{eq:S-alpha}. Hence $S=\emptyset$; that is, $J$ is
  $\lambda_\iota$-connected.
\end{proof}

\begin{cor}
  \label{c:I_i-measure}
  Let $I\subset \Delta_\sigma$ be $\lambda_\iota$-connected. Then for
  any $w_0\in \cW$, the orbit $G(w_0v_0)^{I_\iota(I)}$ admits a
  $G$-invariant measure, say $\nu_{I_\iota(I),w_0}$, such that for any
  $f\in C_c(G(w_0v_0)^{I_\iota(I)})$,
  \begin{align}
    \label{eq:I_iota-measure}
    \int_W f\,d\nu_{I_\iota(I),w_0} = \int\limits_{K} dk
    \int\limits_{\bar a\in\la{c}^+}
    \sum_{w\in\mathcal{W}_{I_\iota(I)}} f
    \big(k\exp(a)w(w_0v_0)^{I_\iota(I)} \big) \xi_{I_\iota(I)}(a)\,
    {d\bar a},
  \end{align}
  where
  \begin{align}
    \label{eq:xiI}
    \xi_{I_\iota(I)}(a)&:=\delta_{I_\iota(I)}(a)
    \exp(\tr(\ad a|_{\la{n}_{I_\iota(I)}})), \\
    \la{e}^+ & = \{\bar
    a\in\la{a}/\la{a}_{I_\iota(I)}\cap\ker\lambda_\iota :
    \alpha(a)\geq 0,\ \forall \alpha\in I_\iota(I)\},
  \end{align}
  and $d\bar a$ denotes the Lebesgue integral on $\la{e}^+$.
\end{cor}

\begin{proof}
  By Proposition~\ref{p:I_iota} $J(I_\iota(I))=I_\iota(I)$. By
  \eqref{eq:C2} for any $y\in\la{a}_{I_\iota(I)}$, $\rho(y)=0$ if and
  only if $\lambda_\iota(y)=0$.  Therefore by
  Theorem~\ref{thm:measure}, $G(w_0v_0)^I$ admits a $G$-invariant
  measure.
  
  Put $\mathfrak{E}=\la{a}_{I_\iota(I)}\cap\ker\lambda_\iota$. Then
  the map $\la{a}^{I_\iota(I)}\oplus \la{a}_{I_\iota(I)}/\mathfrak{E}
  \to \la{a}/\mathfrak{E}$, given by
  \begin{align*}
    (a,b+\mathfrak{E})\mapsto (a+b)+\mathfrak{E},\quad \forall a\in
    \la{a}^{I_\iota(I)}, \,\forall b\in\la{a}_{I_\iota(I)},
  \end{align*}
  is an isomorphism. Note that $\delta_{I_\iota(I)}(b)=1$ and $\tr(\ad
  a|_{\la{n}_{I_\iota(I)}})=0$.  Therefore $\xi_{I_\iota(I)}$ is well
  defined on $\la{a}/\mathfrak{E}$.  Moreover $(a+b)+\mathfrak{E}\in
  \la{e}^+$ if and only if $a\in\la{a}^{I_\iota(I),+}$.  Therefore
  \eqref{eq:I_iota-measure} follows from \eqref{eq:measure}.
\end{proof}

\section{Volume asymptotics}
\label{sec:vol}

In this section we derive some formulas for the volume asymptotics
(see also \cite{GW} and \cite{mau} for a similar computation).

\subsection{Basic asymptotic formula}
\label{sec:b_as}

Consider a space $\mathfrak{a}\simeq \R^r$ and a map
\begin{align*}
  \phi:\mathfrak{a}\to W:\, a\mapsto \sum_{i=1}^k e^{\lambda_i(a)}w_i.
\end{align*}
where $W$ is a finite-dimensional vector space, $w_1,\ldots, w_k\in W$
are linearly independent vectors, and $\lambda_1,\ldots,\lambda_k$ are
(additive) characters.

Fix a basis $\Delta$ of the dual space $\mathfrak{a}^*$ and set
\begin{align*}
  \mathfrak{a}^+=\{a\in\mathfrak{a}:\, \alpha(a)\ge 0 \quad \text{for
    $\alpha\in\Delta$}\}.
\end{align*}
We assume that
\begin{enumerate}
\item $\lambda_1=\sum_{\alpha\in\Delta} m_\alpha\alpha$ with
  $m_\alpha>0$.

\item $\lambda_i\le \lambda_1$ for all $i$, that is,
  $\lambda_1-\lambda_i\in\sum_{\alpha\in\Delta} m_{i,\alpha} \alpha.
  $ with $m_{i,\alpha}\ge 0$.
\end{enumerate}
Let
\begin{align*}
  \supp(\lambda_i)=\{\alpha\in\Delta:\, m_{i,\alpha}>0\}.
\end{align*}
For
\begin{align*}
  \chi=\sum_{\alpha\in\Delta} v_\alpha(\chi)\alpha\in\la{a}^*,
\end{align*}
we set
\begin{align*}
  a_\chi:=\max\left\{\frac{v_\alpha(\chi)}{m_\alpha}:\,
    \alpha\in\Delta\right\},\; I_\chi :=\left\{\alpha\in \Delta:\,
    \frac{v_\alpha(\chi)}{m_\alpha}<a_\chi\right\},\; b_\chi
  :=\#(\Delta-I_\chi).
\end{align*}

Define $\ker I_\chi:=\cap_{\alpha\in I_\chi} \ker\alpha$. Then
\begin{gather*}
  \label{eq:chi-lambda1}
  \chi(a)=a_\chi\cdot\lambda_1(a),\quad  \forall a\in \ker I_\chi, \\
  \label{eq:lambda_1-lambda_i}
  \lambda_i(a)=\lambda_1(a),\quad \forall\,i:\supp \lambda_i\subset
  I_\chi,\ \forall a\in\ker I_\chi, \\
  \label{eq:ker-chi-lambda1}
  \la{d}_0:=\ker\chi\cap \ker I_\chi =\ker \lambda_i\cap \ker I_\chi,
  \quad\forall\,i:\supp \lambda_i\subset I_\chi.
\end{gather*}
Therefore we can define
\begin{align*}
  \la{d}^{+}&:=\{\bar a\in \la{a}/\la{d}_0:
  \alpha(a)\ge 0,\ \forall \alpha \in I_\chi\},\\
  \psi(\bar a)&:=\sum_{i:\supp\lambda_i \subseteq I_\chi}
  e^{\lambda_i(a)}w_i, \quad \forall\, \bar a\in\la{a}/\la{d}_0,
  \quad \text{and}\\
  L_\chi(f)&:=\int_{\mathfrak{d}^{+}} f(\psi(\bar
  a))e^{\chi(a)}\,d\bar a,\quad \forall f\in C_c(W),
\end{align*}
where $d\bar a$ denotes the Lebesgue measure on $\la{a}/\la{d}_0$.

The main result of this subsection is the following theorem
\begin{thm}\label{th:basic_as}
  For $\chi\in \la{a}^*$ and $f\in C_c(W)$,
  \begin{align*}
    \lim_{T\to\infty}\frac{1}{T^{a_\chi}(\log
      T)^{b_\chi-1}}\int_{\mathfrak{a}^{+}} f(\phi(a)/T)e^{\chi(a)}\,
    da= \kappa_\chi \cdot L_\chi(f)<\infty,
  \end{align*}
  where
  \begin{align*}
    \kappa_\chi=\vol(\la{a}^+\cap \ker(I_\chi)\cap \{\lambda_1=1\}).
  \end{align*}
\end{thm}

We start the proof with a lemma:
\begin{lem}\label{l:upper}
  \begin{enumerate}
  \item[(a)] For $T>0$, let $\mathfrak{a}_T^+:=\{a\in
    \mathfrak{a}^+:\, e^{\lambda_1(a)}\le T\}$.  Then
    \begin{align*}
      \int_{\mathfrak{a}^+_T} e^{\chi(a)}\, da\ll T^{a_\chi}(\log
      T)^{b_\chi -1}.
    \end{align*}
  \item[(b)] For $\lambda_i$ such that $\supp \lambda_i \nsubseteq
    I_\chi$ and $T,\delta>0$, set
    \begin{align*}
      \mathfrak{a}_T^+(i,\delta)=\{a\in \mathfrak{a}_T^+:\,
      e^{\lambda_i(a)}\ge \delta T\}.
    \end{align*}
    Then for some constant $C_\delta>1$ depending on $\delta$,
    \begin{align*}
      \int_{\mathfrak{a}_T^+(i,\delta)} e^{\chi(a)}\, da \leq C_\delta
      T^{a_\chi}(\log T)^{b_\chi -2}.
    \end{align*}
  \end{enumerate}
\end{lem}

\begin{proof}
  To prove (a), we use induction on $|I_\chi|$.  If
  $I_\chi=\emptyset$, then $\chi=a_\chi\cdot \lambda_1$ and
  \begin{align*}
    \int_{a\in \mathfrak{a}^+:\lambda_1(a)\le \tau} e^{\chi(a)}\, da&=
    \int_0^\tau \vol(\mathfrak{a}^+\cap\{\lambda_1=s\})e^{a_\chi
      s}\, ds\\
    &=\int_0^\tau (cs^{r-1})e^{a_\chi s}\, ds=O(\tau^{r-1}e^{a_\chi
      \tau}).
  \end{align*}
  Let $\alpha\in I_\chi$ and
  $\mathfrak{b}^+=\mathfrak{a}^+\cap\ker(\alpha)$. Then by the
  inductive assumption,
  \begin{align*}
    &\int_{a\in \mathfrak{a}^+:\lambda_1(a)\le \tau} e^{\chi(a)}\, da=
    \int_0^{\tau/m_\alpha} e^{v_\alpha(\chi)s} \left(\int_{b\in
        \mathfrak{b}^+:\lambda_1(b)\le \tau-m_\alpha s} e^{\chi(b)}\,
      db\right)\,
    ds\\
    \ll &\int_0^{\tau/m_\alpha} e^{v_\alpha(\chi)s} (\tau-m_\alpha
    s)^{b_\chi-1}e^{a_\chi (\tau-m_\alpha s)}\, ds\ll
    e^{a_\chi\tau}\tau^{b_\chi -1},
  \end{align*}
  where $C>1$ is a constant. This proves (a).

  To prove (b), we write $\lambda_i=\lambda_1-\sum_{\alpha\in
    \supp(\lambda_i)} m_{i,\alpha} \alpha$ with $m_{i,\alpha}>0$.  For
  $a\in\mathfrak{a}_T^+(i,\delta)$, we have
  \begin{align*}
    \sum_{\alpha\in \supp \lambda_i} m_{i,\alpha} \alpha(a)\le -\log
    \delta.
  \end{align*}
  Setting $\mathfrak{c}=\ker(\supp \lambda_i)$ and
  $\mathfrak{c}_T^+=\mathfrak{a}_T^+\cap\mathfrak{c}$, we get
  \begin{align*}
    \int_{\mathfrak{a}^+_T(i,\delta)} e^{\chi(a)}\, da \leq C_\delta
    \int_{\mathfrak{c}^+_T} e^{\chi(c)}\, dc,
  \end{align*}
  for some constant $C_\delta>1$ depending on $\delta$. Since $\supp
  \lambda_i\nsubseteq I_\chi$, $b_{\chi|_{\mathfrak{c}}}\le
  b_{\chi}-1$, and (b) follows from (a).
\end{proof}

\begin{proof}[Proof of Theorem~\ref{th:basic_as}]
  There exists $c=c(f)>0$ such that if $f(\phi(a)/T)\ne 0$ or
  $f(\psi(a)/T)\ne 0$, then $a\in\mathfrak{a}^+_{cT}$.

  For $\delta>0$, set
  \begin{align*}
    \mathfrak{a}_T^+(\delta)=\bigcup_{i: \supp(\lambda_i)\nsubseteq
      I_\chi} \mathfrak{a}_T^+(i,\delta).
  \end{align*}
  where $\la{a}_T^+(i,\delta)$ is defined as in Lemma~\ref{l:upper}.
  By uniform continuity, for every $\e>0$, there exists $\delta>0$
  such that for $a\in
  \mathfrak{a}_{cT}^+-\mathfrak{a}_{cT}^+(\delta)$,
  \begin{align*}
    |f(\phi(a)/T) -f(\psi(a)/T)|<\e.
  \end{align*}
  Hence, by Lemma~\ref{l:upper},
  \begin{align*}
    &\left|\int_{\mathfrak{a}^{+}} f(\phi(a)/T)e^{\chi(a)}\,da
      - \int_{\mathfrak{a}^{+}}f(\psi(a)/T)e^{\chi(a)}\, da\right|\\
    \le& \int_{\mathfrak{a}_{cT}^+-\mathfrak{a}_{cT}^+(\delta)}
    \e\cdot e^{\chi(a)}\, da
    + \int_{\mathfrak{a}_{cT}^+(\delta)}(2\norm{f}_\infty)\cdot e^{\chi(a)}\, da\\
    =&\, O_f(\e\cdot T^{a_\chi}(\log T)^{b_\chi-1})+
    O_{f,\delta}(T^{a_\chi}(\log T)^{b_\chi-2}).
  \end{align*}
  This shows that
  \begin{align*}
    \int_{\mathfrak{a}^{+}} f(\phi(a)/T)e^{\chi(a)}\, da =
    \int_{\mathfrak{a}^{+}} f(\psi(a)/T)e^{\chi(a)}\,
    da+o(T^{a_\chi}(\log T)^{b_\chi-1}).
  \end{align*}
  Let
  \begin{align*}
    \la{s}^{+}=\la{a}^{+}\cap \ker(I_\chi)\quad\text{and}\quad
    \la{t}^+=\la{a}^{+}\cap \ker(\Delta-I_\chi).
  \end{align*}
  There exists $c=c(f)>0$ such that if $f(\psi(t)e^u)\ne 0$ for some
  $t\in\la{t}^+$ and $u\in\R$, then $e^{\lambda_1(t)}\le c e^{-u}$ and
  $e^u\le c$.  Using that for some $\e>0$,
  \begin{align*}
    \chi|_{\la{t}^+}\le (a_\chi-\e)\cdot \lambda_1|_{\la{t}^+},
  \end{align*}
  we deduce that
  \begin{align}
    \nonumber \int_{\mathfrak{t}^{+}}\int_{-\infty}^\infty
    & f(\psi(t)e^{u})e^{\chi(t)+a_\chi u} |u|^l\, dudt \\
    & \ll \int\limits_{u=-\infty}^{\log c} \, du
    \int\limits_{\{t\in\mathfrak{t}^{+}:\lambda_1(t)\le \log c-u\}}
    e^{(a_\chi-\e)\lambda_1(t)+a_\chi u} |u|^l\, dt
    \label{eq:integral-finite} \\
    &\ll \int_{-\infty}^{\log c} e^{\e u} (\log c-u)^{r-b_\chi}|u|^l\,
    du<\infty \nonumber
  \end{align}
  for every $l\in \N\cup\{0\}$. In particular, putting $l=0$ we get
  \begin{align*}
    L_\chi(f) =\int_{\la{d}^+} f(\psi(\bar a))e^{\chi(a)}d\bar a =\int
    _{t\in\la{t}^+}\int_{u=-\infty}^\infty
    f(\psi(t)e^u)e^{\chi(t)+a_\chi u} dtdu < \infty.
  \end{align*}

  Therefore applying \eqref{eq:integral-finite}, we conclude that as
  $T\to\infty$,
  \begin{align*}
    \int_{\mathfrak{a}^{+}} f(\psi(a)/T)e^{\chi(a)}\,
    da&=\int_{\mathfrak{t}^{+}}\int_{\mathfrak{s}^{+}}
    f(\psi(t)e^{\lambda_1(s)}/T)e^{\chi(t)+a_\chi\lambda_1(s)}\, dsdt\\
    &=\int_{\mathfrak{t}^{+}}\int_0^\infty
    f(\psi(t)e^{u}/T)e^{\chi(t)+a_\chi u}\cdot\vol(\la{s}^+\cap \{\lambda_1=u\})\, dudt\\
    &=\int_{\mathfrak{t}^{+}}\int_{-\log T}^\infty
    f(\psi(t)e^{u})e^{\chi(t)+a_\chi u} T^{a_\chi}\cdot \kappa_\chi
    (u+\log
    T)^{b_\chi-1}\, dudt\\
    &=\kappa_\chi\cdot L_\chi(f)\cdot T^{a_\chi}(\log T)^{b_\chi-1}+
    o(T^{a_\chi}\cdot (\log T)^{b_\chi-1}).
  \end{align*}
  This completes the proof.
\end{proof}

\subsection{Volume of symmetric space}
Let $G$ be a connected noncompact semisimple Lie group with finite
center, $H$ its symmetric subgroup, and $\iota:G\to \GL(W)$ be an
almost faithful irreducible over $\R$ representation.  We assume that
$H=\Stab_G(v_0)$ for some $v_0\in W$.  We use notation from
Section~\ref{sec:basic}. In particular, $\mu$ denotes an invariant
measure on $G/H$ and
\begin{align*}
  \lambda_\iota=\sum_{\alpha\in\Delta_\sigma} m_\alpha\alpha\in
  \mathfrak{a}^*,\quad m_\alpha\in\Q^+,
\end{align*}
is the highest weight of $\iota$.  Let $a_\iota$, $b_\iota$, $I_\iota$
be defined as in (\ref{eq:defi}).

\begin{thm}\label{l:asympt0}
  For every $f\in C_c(W)$,
  \begin{equation}\label{eq:lim}
    \lim_{T\to\infty}\frac{1}{T^{a_\iota}(\log T)^{b_\iota-1}}\int_{G/H} f(gv_0/T)\,d\mu(g)=
    \int_{W} f\, d\nu_\iota,
  \end{equation}
  where $\nu_\iota$ is a locally finite $G$-invariant measure on $W$
  concentrated on $\R^+\cdot V_{I_\iota}^\infty$.

Moreover when considered
as a measure on $\br^+\cdot V_{I_\iota}^\infty$,
$\nu_\iota$ is a 
linear combination of measures $\nu_{I_\iota,
    w}$, $w\in \mathcal W^{I_\iota}$,  given in
  (\ref{eq:I_iota-measure}).
\end{thm}

\begin{proof}
  For $v\in W$, set
  \begin{align*}
    \Phi_\iota(v)=\{\lambda\in\Phi_\iota:\, v^\lambda\ne 0\}.
  \end{align*}
  By (\ref{eq:volume}),
  \begin{align}\label{eq:gt}
    \int_{G/H} f(g v_0/T)\, d\mu(g) &= \int_K\sum_{w\in\mathcal{W}}
    \int_{\mathfrak{a}^{+}}
    f(k\exp(a)w v_0/T)\xi(a)\,dadk\\
    &=\int_K\sum_{w\in\mathcal{W}} \int_{\mathfrak{a}^{+}}
    f\left(k\cdot\sum_{\lambda\in \Phi_\iota(wv_0)} e^{\lambda(a)}(w
      v_0)^\lambda/T \right)\xi(a)\,dadk,\nonumber
  \end{align}
  where
  \begin{equation}\label{eq:xi_0}
    \xi(a)=\prod_{\alpha\in \Sigma^+}
    (\sinh\alpha(a))^{l_\alpha^+}(\cosh\alpha(a))^{l_\alpha^-} = 
    \sum_{\chi\in \Xi} t_\chi e^{\chi (a)}
  \end{equation}
  for some $t_\chi\ne 0$ and $\Xi\subset \mathfrak{a}^*$. Let
  $a_\iota$, $b_\iota$, and $I_\iota$ be as defined in
  \eqref{eq:defi}. Let
  \begin{align*}
    \Xi'=\{\chi\in \Xi:\, a_\chi= a_\iota,
    b_\chi=b_\iota\}\quad\text{and}\quad\Xi''=\Xi-\Xi'.
  \end{align*}
  Note that for $\chi\in \Xi''$, we have $a_\chi\le a_\iota=a_{2\rho}$
  and if $a_\chi= a_\iota$, then $b_\chi<b_\iota=b_{2\rho}$ and
  $I_\chi\supset I_\iota=I_{2\rho}$.

  Since by Lemma~\ref{l:weights}, $\lambda_\iota\in \Phi_\iota(wv_0)$,
  the assumptions of Section~\ref{sec:b_as} are satisfied and applying
  Theorem~\ref{th:basic_as} together with the dominated convergence
  theorem, we deduce that (\ref{eq:lim}) holds with the measure
  $\nu_\iota$ given by the formula
  \begin{equation}\label{eq:nu_mu}
    \int_{W} f\, d\nu_\iota=
    \kappa_{2\rho}\int_K\sum_{w\in\mathcal{W}} \int_{\mathfrak{d}^{+}}
    f\left(k\exp(a)(w
      v_0)^{I_\iota} \right)\xi_{I_\iota} (a)\,dadk
  \end{equation}
  where
  \begin{align} \label{eq:xi-sum} \la{d}^{+}&=\{a\in
    \la{a}/(\ker(I_\iota)\cap \ker(\rho)):\,
    \alpha(a)\ge 0,\, \alpha\in I_\iota\},\nonumber\\
    \xi_{I_\iota}(a)&=\sum_{\chi\in \Xi'} t_\chi e^{\chi (a)}.
  \end{align}
  Note that by Theorem~\ref{th:basic_as} the limit in (\ref{eq:lim})
  is finite (i.e., $\nu_\iota$ is locally finite).  Also, it is clear
  from (\ref{eq:lim}) that $\nu_\iota$ is $G$-invariant and
  homogeneous of degree $a_\iota$.  It follows from
  Proposition~\ref{p:orbit} that
  \begin{align*}
    G(wv_0)^{I_\iota}=
    K\exp(\la{d}^{+})\mathcal{W}_{I_\iota}(wv_0)^{I_\iota} =\R^+\cdot
    K \exp(\la{a}^{I_\iota,+})\mathcal{W}_{I_\iota}(wv_0)^{I_\iota}.
  \end{align*}
Note that $$V_{I_\iota}^\infty=\cup_{w\in\cW}\;
\pi(  G(wv_0)^{I_\iota})=\cup_{w\in\cW^{I_\iota}}\;
 \pi( G(wv_0)^{I_\iota}) .$$

  Since for $\chi\in \Xi$, we have
 $$\chi\in\Xi'\quad\text{ if and only if }\quad \chi\in
  2\rho+\langle I_\iota \rangle,$$ it follows that the formula for
  \eqref{eq:xi-sum} for $\xi_{I_\iota}$ is same as the formula
  \eqref{eq:xiI} of Corollary~\ref{c:I_i-measure} for
  $I_\iota=I_\iota(\emptyset)$.  
Note that
each $G$-orbit $G(wv_0)^{I_\iota}$ is a closed
subset of  $\br^+\cdot V_{I_\iota}^\infty$
and hence $f\in C_c(\br^+\cdot V_{I_\iota}^\infty)$ implies
the restriction of $f$
to $G(wv_0)^{I_\iota}$ belongs to $ C_c(G(wv_0)^{I_\iota})$.  
Hence (\ref{eq:nu_mu}) is in
  agreement with (\ref{eq:I_iota-measure}) for
$f\in C_c(\br^+\cdot V_{I_\iota}^\infty)$.  

This shows that
  $\nu_\iota$, considered as
a measure on $\br^+\cdot V_{I_\iota}^\infty$,
 is a linear combination of the measures $\nu_{I_\iota,
    w}$, $w\in \mathcal W^{I_\iota}$, given in
  (\ref{eq:I_iota-measure}).
It follows that $\nu_\iota$ is concentrated on
$\br^+\cdot V_{I_\iota}^\infty$.

\end{proof}

We can easily deduce the following volume asymptotic of balls from the
above theorem:
\begin{cor} For any norm $\|\cdot \|$ on $W$,
  \begin{align*}
    \operatorname{Vol}(\{v\in V:\, \|v\|<T\}) \sim c\cdot
    T^{a_\iota}(\log T)^{b_\iota} \quad\text{as $T\to \infty$},
  \end{align*}
  where $\operatorname{Vol}$ denotes a $G$-invariant measure on $V$
  and $c>0$.
\end{cor}

\begin{rem} \label{rem:rep} \rm Theorem~\ref{l:asympt0} holds for a
  representation $\iota$ which is not irreducible. Let
  \begin{align*}
    \mathcal{P} &=\{a\in\la{a}^+:\, \lambda(a)\le 1,\lambda\in\Phi_\iota\},\\
    a_\iota &=\max\{2\rho(a):\, a\in\mathcal{P}\},\\
    b_\iota &=\dim \mathcal{P}\cap \{2\rho=a_\iota\}.
  \end{align*}
  Then for every $f\in C_c(W)$,
  \begin{align*}
    \lim_{T\to\infty}\frac{1}{T^{a_\iota}(\log T)^{b_\iota-1}}
    \int_{G/H} f(gv_0/T)\,d\mu(g)= \int_{W} f\, d\nu_\iota,
  \end{align*}
  where $\nu_\iota$ is a $G$-invariant measure concentrated on a union
  of finitely many $G$-orbits.  To adapt the proof to this case, we
  decompose the polyhedron $\mathcal{P}$ into a finite union of
  symplicial polyhedra $\mathcal{P}_i$.  the asymptotics for the
  integral over $K\exp(\mathcal{P}_i)\mathcal{W}H$ can be computed
  from Theorem~\ref{th:basic_as}.  Using the argument from
  Section~\ref{sec:equ}, we also get the asymptotics for integral
  points.
\end{rem}

For $f\in C_c(W\setminus\{0\})$ with $\pi(\supp f)\cap
V^\infty\neq\emptyset$, we define $a_\iota(f)$, $b_\iota(f)$,
$\Theta_\iota(f)$ as in (\ref{eq:iota0}) with $\Omega=\pi(\supp f)$.
Similarly, we define $\Theta_f$ and $\Theta_f^{\min}$.

\begin{thm}\label{th:asympt}
  For every $f\in C_c(W\setminus\{0\})$ with $\pi(\supp f)\cap
  V^\infty\neq\emptyset$,
  \begin{align}\label{nutf}
    \lim_{T\to\infty} \frac{1}{T^{a_{\iota}(f)}(\log
      T)^{b_{\iota}(f)-1}}\int_{G/H} f(gv_0/T)d\mu(g)=\int_W f\,
    d\nu_{\Theta_\iota(f)},
  \end{align}
  where $\nu_{\Theta_\iota(f)}$ is a $G$-invariant measure on $W$
  which is concentrated on and locally finite on
  \begin{align*}
    \cup_{I\in\Theta_\iota(f)} \R^+\cdot V_I^\infty.
  \end{align*}
  In particular, $\nu_{\Theta_\iota(f)}(\supp f)<\infty$.
\end{thm}

\begin{proof}
  Define $\tilde W=\R^+\cdot \overline{\pi(V)}$ and $S(\tilde
  W):=S(W)\cap \tilde{W}=\pi(\tilde W)$, where $S(W)$ and $\pi:W\to
  S(W)$ are defined as in \eqref{eq:S(W)}. 
Since
$f\in C_c(W\setminus\{0\})$,
 $f|_{\tilde W}\in
  C_c(\tilde W)$.
Hence it suffices to prove the theorem for $f\in C_c(\tilde W)$.
  For $I\subset \Delta_\sigma$, we set
  \begin{align*}
    \mathcal{O}_I=\bigcup_{\lambda_\iota\text{-connected }J\supset I}
    V_J^\infty.
  \end{align*}
  It follows from Proposition~\ref{p:closure} that $\mathcal{O}_I$ is
  open in $S(\tilde W)$. We take a partition of unity $\phi_I\in
  C(S(\tilde W))$, $I\in \Theta_f^{\min}$, associated to the cover
  \begin{align*}
    \pi(\supp f)\subset \bigcup_{I\in \Theta_f^{\min}} \mathcal{O}_I.
  \end{align*}
  It suffices to prove the theorem for the functions
  $f_I(v)=f(v)\phi_I(\pi(v))$, $I\in \Theta_f^{\min}$. Hence, we may
  assume that $\Theta_f^{\min}=\{I\}$ for some
  $\lambda_\iota$-connected $I\subset \Delta_\sigma$.  Then
  \begin{align*}
    (a_\iota(f),b_\iota(f))= (a_\iota(I),b_\iota(I))
    \quad\text{and}\quad \Theta_\iota(f)=\{I_\iota(I)\}.
  \end{align*}

  By Corollary~\ref{c:satake}, there exists a compact set $U\subset
  \la{a}^{I,+}$ such that
  \begin{align*}
    \pi(\supp f)\cap V^\infty & \subset \bigcup_{I\in\Theta_f}
    \pi (K\exp(U+\la{a}_{I}^+)(\mathcal{W}v_0)^I),\\
    \pi(\supp f)\cap \pi(V) &\subset
    \pi(K\exp(U+\mathfrak{a}_{I}^{+})\mathcal{W}v_0).
  \end{align*}
  Hence, as in (\ref{eq:gt}), we have
  \begin{align*}
    &\int_{G/H} f(g v_0/T)\, d\mu(g) \\
    &= \int_K\sum_{w\in\mathcal{W}} \int_U \int_{\mathfrak{a}_{I}^{+}}
    f\left(k\cdot\sum_{\lambda\in \Phi_\iota(wv_0)}
      e^{\lambda(a)}e^{\lambda(u)}(w v_0)^\lambda/T
    \right)\xi(a+u)\,dadudk.
  \end{align*}
  We apply Theorem~\ref{th:basic_as} to the integral over
  $\la{a}^+_{I}$ in place of $\la{a}$.  For
  $\chi=\sum_{\alpha\in\Delta_\sigma} v_\alpha\alpha$, we set
  \begin{align*}
    a_\chi &:= \max\left\{\frac{v_\alpha}{m_\alpha}: \alpha\in
      \Delta_\sigma
      \setminus I \right\},\\
    I_\chi &:=\{\alpha\in \Delta_\sigma\setminus I:
    \frac{v_\alpha}{m_\alpha}
    <a_\chi\},\\
    b_\chi &:=\#((\Delta_\sigma\setminus I)\setminus I_\chi).
  \end{align*}
  We write
  \begin{align*}
    \xi(a)= \sum_{\chi\in \Xi} t_\chi e^{\chi (a)}
  \end{align*}
  for some $t_\chi\ne 0$ and $\Xi\subset \mathfrak{a}^*$. Let
  \begin{align*}
    \Xi'=\{\chi\in \Xi:\, a_\chi= a_{2\rho}, b_\chi=b_{2\rho}\}
    \quad\text{and}\quad \Xi''=\Xi-\Xi'.
  \end{align*}
  Note that for $\chi\in \Xi''$, we have $a_\chi\le
  a_{2\rho}=a_\iota(I)$ and if $a_\chi= a_{2\rho}$, then
  $b_\chi<b_{2\rho}=b_\iota(I)$ and $I_\chi\supset I_{2\rho}$ where
  $I_{2\rho}\cup I=I_\iota(I)$.  By Theorem~\ref{th:basic_as} and the
  dominated convergence theorem,
  \begin{align}
    &\lim_{T\to\infty} \frac{1}{T^{a_{2\rho}}(\log
      T)^{b_{2\rho}-1}}\int_{G/H} f(g v_0/T)\, d\mu(g)\nonumber\\
    =& \kappa\int_K\sum_{w\in\mathcal{W}}
    \int_U\int_{\mathfrak{d}^{+}} f\left(k\exp(u+a)(w
      v_0)^{I_\iota(I)} \right)\xi_{I_\iota(I)}(u+a)\,da dudk\nonumber\\
    =& \kappa\int_K\sum_{w\in\mathcal{W}} \int_{\mathfrak{e}^{+}}
    f\left(k\exp(a)(w v_0)^{I_\iota(I)} \right)\xi_{I_\iota(I)}(a)\,da
    dk
    \label{eq:lim_me}
  \end{align}
  where
  \begin{align}
    \nonumber
    \kappa&=\vol(\la{a}_{I}^+\cap \ker(I_{2\rho})\cap \{\lambda_\iota=1\}),\\
    \nonumber \la{d}^{+}&=\{a\in \la{a}_{I}/(\ker(I_{2\rho})\cap
    \ker(\rho)):\, \nonumber
    \alpha(a)\ge 0,\, \alpha\in I_{2\rho}\},\\
    \nonumber \la{e}^{+}&=\{a\in \la{a}/(\ker(I_\iota(I))\cap
    \ker(\rho)):\,
    \alpha(a)\ge 0,\, \alpha\in I_{\iota}(I)\},\\
    \label{eq:XiI_iota}
    \xi_{I_\iota(I)}(a)&=\sum_{\chi\in \Xi'} t_\chi e^{\chi (a)}.
  \end{align}

Therefore (\ref{nutf}) holds with the measure
  $\nu_{\Theta_\iota(f)}$ given by the formula
  \begin{equation*}
    \int_{W} f\, d\nu_{\Theta_{\iota(f)}}=
     \kappa\int_K\sum_{w\in\mathcal{W}} \int_{\mathfrak{e}^{+}}
    f\left(k\exp(a)(w v_0)^{I_\iota(I)} \right)\xi_{I_\iota(I)}(a)\,da
    dk
  \end{equation*}

  It is clear that   $\nu_{\Theta_\iota(f)}$
 is $G$-invariant and homogeneous
  of degree $a_\iota(f)$.  Also, it follows from
  Proposition~\ref{p:orbit} that
  \begin{align*}
    K\exp(\la{e}^{+})\mathcal{W}_{I_\iota(I)}(wv_0)^{I_\iota(I)} =
    \R^+ \cdot K
    \exp(\la{a}^{I_\iota(I),+})\mathcal{W}_{I_\iota(I)}(wv_0)^{I_\iota(I)}
  \end{align*}
  is a single $G$-orbit.
  It follows from Theorem~\ref{th:basic_as} that
  $\nu_{\Theta_\iota(f)}$ is locally finite
on $\br^+\cdot V_{I_\iota(I)}^\infty$.
 
  Since $\chi\in\Xi'$ if and only if $\chi\in
  2\rho+\Sigma_{I_\iota(I)}$, the formula \eqref{eq:XiI_iota} for
  $\xi_{I_\iota(J)}$ is same as \eqref{eq:xiI} in
  Corollary~\ref{c:I_i-measure}. Hence (\ref{eq:I_iota-measure})
  agrees with (\ref{eq:lim_me}).  Since $\xi_{I_\iota(I)}\ne 0$ on a
  set of full Lebesgue measure on $\la{d}^+$, the limit measure is
  strictly positive on nonempty open subsets of
  $G(\mathcal{W}v_0)^{I_\iota(I)}$. 
This shows that  $\nu_{\Theta_\iota(f)}$ is concentrated on
$\br^+\cdot V_{I_\iota(I)}^\infty$, proving the theorem.
\end{proof}

\begin{rem} \label{r:ab} \rm For any $f\in C_c(W\setminus\{0\})$ and a
  $\lambda_\iota$-connected $I\subset\Delta_\sigma$ satisfying the
  conditions of Theorem~\ref{thm:measure}, if
  \begin{align*}
    (a_\iota(f),b_\iota(f))<(a_\iota(I),b_\iota(I))\leq (a,b),
  \end{align*}
  with respect to the lexicographic order on the pairs, then by
  Theorem~\ref{th:asympt},
  \begin{align*}
    \lim_{T\to\infty} \frac{1}{T^{a}(\log T)^{b-1}} \int_{G/H}
    f(gv_0/T)\,dg=0=\int_W f\,d\nu_I,
  \end{align*}
  where $\nu_I$ is a $G$-invariant measure on $W$ concentrated on
  $\R^+\cdot V_I^\infty$.
\end{rem}

\section{Distribution of integral points}\label{sec:equ}
Let $G$ be a connected noncompact semisimple Lie group with finite
center, $H$ a symmetric subgroup $G$, and $\iota: G\to \GL(W)$ an
almost faithful irreducible over $\R$ representation of $G$ such that
for some $v_0\in W$, $\Stab_G(v_0)=H$.

Let $\Gamma$ be an irreducible lattice in $G$ such that $H\cap \Gamma$
is a lattice in $H$. We choose Haar measures $dg$ $dh$, $d\mu$ on $G$,
$H$, $G/H$ respectively such that
\begin{align*}
  \int_G f dg=\int_{G/H} \int_{H} f(gh)\,dh\; d\mu(g),\quad f\in
  C_c(G).
\end{align*}
It is convenient to normalize the measures so that
\begin{align*}\vol(G/\Gamma)=\vol(H/(H\cap\Gamma))=1.\end{align*}
The following result was proved in \cite{EM} (see also \cite{DRS}):

\begin{Thm}\label{mix}
  For every $\phi\in C_c(G/\G)$,
  \begin{align*} \int_{H/(H\cap\Gamma)} \phi(v h)\, dh \to \int_{G/\G
    }\phi \, dg \quad \text{as $v\to\infty$ in $G/H$}.
  \end{align*}
\end{Thm}

\begin{rem} \label{rem:mod} \rm The condition that the lattice
  $\Gamma$ is irreducible in $G$ can be relaxed. In fact, it suffices
  to assume that $G=G_1\cdots G_r$ for noncompact normal subgroups
  $G_i$'s such that $\Gamma\cap G_i$ is an irreducible lattice in
  $G_i$ and $G=G_iH$ for all $i$. If this is the case, then for
  $v=g_1\cdots g_rH\to \infty$, we have $g_i\to\infty$ for all $i$ and
  Theorem~\ref{mix} holds (see \cite[Corollary~1.2]{shah}).
\end{rem}

For $T> 0$ and $f\in C_c(W)$, define
\begin{equation}\label{eq:F_T}
  F_T(g)=\sum_{\gamma\in \G/(\G\cap H)} f(g\gamma v_0/T),\quad g\in G/\Gamma.
\end{equation}

\begin{Prop}\label{p:equi}
  Let $\phi\in C_c(G/\G)$ such that $\int_{G/\Gamma} \phi\, dg=1$ and
  $f\in C_c(W\setminus\{0\})$ with $\pi(\supp f)\cap V^\infty\neq
  \emptyset$. Then
  
  \begin{align*}
    \lim_{T\to\infty} \ \frac{1}{T^{a_\iota(f)}(\log
      T)^{b_\iota(f)}}\cdot \inpr{F_T}{\phi} =\int_W
    f\,d\nu_{\Theta_\iota(f)},
  \end{align*}
  where $\nu_{\Theta_\iota(f)}$ is as given by
  Theorem~\ref{th:asympt}.  Furthermore,
  
  \begin{align*}
    \lim_{T\to\infty} \int_{G/H} f(gv_0/T)\,d\mu(g)=\infty \implies
    \inpr{F_T}{\phi} \sim \int_{G/H} f(gv_0/T)\,d\mu(g) \quad \text{as
      $T\to\infty$}.
  \end{align*}

\end{Prop}

\begin{proof}
  We have
  \begin{align*}
    \inpr{F_T}{\phi} &=\int_{G/\Gamma} \sum_{\gamma\in \G/H\cap \G}
    f(g\gamma v_0/T) \phi(g)
    \, dg=\int_{G/H\cap \Gamma}  f(g v_0/T) \phi(g) \, dg\\
    &=\int_{G/H} f(gv_0/T) \left(\int_{H/H\cap \G} \phi(gh)\,
      dh\right)\, d\mu(g).
  \end{align*}
  
  By Theorem~\ref{mix}, for every $\e>0$, there exists a compact set
  $D\subset G/H$ such that
  \begin{align*}
    \left|\int_{H/(H\cap\Gamma)} \phi(g h)\, dh - 1 \right|<\e
  \end{align*}
  for $g\in G/H\setminus D$. Then
  \begin{align}
    \label{eq:equidis}
    &\abs{\inpr{F_T}{\phi}-\int_{G/H}
      f(gv_0/T)\,d\mu(g)} \\
    \nonumber &\quad \quad\le \e \abs{\int_{G/H\setminus D}
      f(gv_0/T)\, d\mu(g)}
    +\mu(D)\norm{f}_\infty(\norm{\phi}_\infty+1).
  \end{align}
  The second part of the proposition now follows immediately. And the
  first part of the proposition follows from Theorem~\ref{th:asympt}.
\end{proof}

\begin{thm}\label{th:sector2}
  For every $f\in C_c(W\setminus\{0\})$ with $\pi(\supp f)\cap
  V^\infty\neq\emptyset$,
  \begin{align*}
    \lim_{T\to\infty} \frac{1}{T^{a_\iota(f)}(\log
      T)^{b_\iota(f)-1}}\sum_{\gamma\in \Gamma/(\Gamma\cap H)}
    f(\gamma v_0/T) =\int_W f\, d\nu_{\Theta_\iota(f)},
  \end{align*}
  where $\nu_{\Theta_\iota(f)}$ is as in Theorem~\ref{th:asympt}.
\end{thm}

\begin{proof}
  Without loss of generality, we may assume that $f\ge 0$.  Given
  $\e>0$ there exists a compact symmetric neighborhood $\cO_\e$ of $e$
  in $G$ such that
  \begin{align*}
    |f(gv)-f(v)|<\e\quad \forall g\in\mathcal{O}_\e,\,\forall v\in W,
  \end{align*}
  We can assume that $\cO_\e\subset\cO_1$. Define $f^\pm_\e\in C_c(W)$
  by
  \begin{align}
    \label{eq:f_e+}
    f_\e^+(v):=\max_{g\in \mathcal{O}_\e} f(g\cdot v) \quad \text{and}
    \quad f_\e^{-}(v)=\min_{g\in \mathcal{O}_\e} f(g\cdot v),
    \quad\forall v\in W.
  \end{align}

  For any $T>0$, let $F_T$ and $F_T^\pm$ be defined as in
  (\ref{eq:F_T}) corresponding to $f$ and $f^\pm_\e$, respectively.
  Then
  \begin{align*}
    % \label{eq:F_T_pm}
    F_T^-(g)\le F_T (e) \le F_T^+(g),\quad \forall g\in\cO_\e.
  \end{align*}
  Hence, if $\phi\in C_c(G/\Gamma)$, with $\phi\ge 0$,
  $\supp\phi\subset \Cal O_\e$ and $\int_{G/\Gamma} \phi dg=1$, then
  \begin{align}
    \label{eq:pm}
    \inpr{F_T^-}{\phi }\le \sum_{\gamma\in \Gamma/(\Gamma\cap H)}
    f(\gamma v_0/T)=F_T(e)\le \inpr{F_T^+}{\phi}.
  \end{align}

  Since each $\R^+\cdot V^\infty_I$ is $G$-invariant, we have
  \begin{gather}
    \label{eq:Theta_iota}
    \Theta_\iota(f_\e^+) = \Theta_\iota(f)\supset
    \Theta_\iota(f_\e^-) \\
    \label{eq:a_iota}
    (a_\iota(f_\e^+),b_\iota(f_\e^+))=(a_\iota(f),b_\iota(f))\ge
    (a_\iota(f_\e^-),b_\iota(f_\e^-)).
  \end{gather}

  In view of Remark~\ref{r:ab}, by Proposition~\ref{p:equi}, we get
  \begin{align}
    \label{eq:pm-limit}
    \lim_{T\to\infty} \frac{\inpr{F_T^\pm}{\phi}}{T^{a_\iota(f)}(\log
      T)^{b_\iota(f)-1}} =\int_W f^\pm_\e\, d\nu_{\Theta_\iota(f)}.
  \end{align}

  Combining \eqref{eq:pm} and \eqref{eq:pm-limit} we conclude that
  \begin{align}
    \label{eq:pm-main}
    \int_W f^-_\e\,d\nu_{\Theta_\iota(f)} &\leq\liminf_{T\to\infty}
    \frac{F_T(e)}{T^{a_\iota(f)}(\log T)^{b_\iota(f)-1}} \\ &\leq
    \limsup_{T\to\infty} \frac{F_T(e)}{T^{a_\iota(f)}(\log
      T)^{b_\iota(f)-1}} \leq \int_W f^+_\e\,d\nu_{\Theta_\iota(f)}.
  \end{align}

  By \eqref{eq:f_e+},
  \begin{align}
    \label{eq:f_e}
    \int_{W} f_\e^\pm\, d\nu_{\Theta_\iota(f)} - \int_{W} f\,
    d\nu_{\Theta_\iota(f)}\le \e \cdot \nu_{\Theta_\iota(f)}(\supp
    f_\e^+)\leq \e\cdot\nu_{\Theta_\iota(f_1^+)}(\supp f_1^+).
  \end{align}
  By Theorem~\ref{th:asympt} and \eqref{eq:Theta_iota},
  $\nu_{\Theta_\iota(f_1^+)}(\supp f_1^+)<\infty$. Since $\epsilon>0$ can be
  chosen arbitrarily small, \eqref{eq:sector:iota} follows from
  \eqref{eq:pm-main} and \eqref{eq:f_e}.
\end{proof}

Note that for any $f\in C_c(W\setminus\{0\})$, we have
$(a_\iota,b_\iota)\geq (a_\iota(f), b_\iota(f))$. Therefore using
Remark~\ref{r:ab}, from Theorem~\ref{th:sector2} and
Theorem~\ref{l:asympt0} we can deduce the following.

\begin{thm}\label{th:sector}
  For every $f\in C_c(W)$,
  \begin{align}
    \label{eq:sector:iota}
    \lim_{T\to\infty} \frac{1}{T^{a_\iota}(\log
      T)^{b_\iota-1}}\sum_{\gamma\in \Gamma/(\Gamma\cap H)} f(\gamma
    v_0/T)= \int_W f\, d\nu_\iota,
  \end{align}
  where $\nu_\iota$ is as in Theorem~\ref{l:asympt0}.  \qed
\end{thm}

\begin{proof}[Proof of Theorem~\ref{th:m}]
  Let $\Gamma\subset G(\Q)$ be an arithmetic subgroup that preserves
  the integral structure on $W(\Z)$.  Since ${\mathbf G}$ and
  ${\mathbf H}$ admit no nontrivial $\Q$-characters, by \cite{BH},
  $\Gamma$ is an irreducible lattice in $G$, $\Gamma\cap H$ is a
  lattice in $H$, and $V(\z)$ is a union of finitely many orbits of
  $\G$:
  \begin{align*}
    V(\Z)=\bigcup_{i=1}^n \Gamma g_iv_0.
  \end{align*}

  For any $T>0$, consider a locally finite measure $\tau_T$ on $W$
  defined by
  \begin{align*}
    \tau_T(f)&=\frac{1}{T^{a_\iota}(\log T)^{b_\iota-1}}\sum_{v\in
      V(\Z)} f(v/T),\quad f\in C_c(W).
  \end{align*}
  Then
  \begin{align*}
    \tau_T(f)&=\frac{1}{T^{a_\iota}(\log
      T)^{b_\iota-1}}\sum_{i=1}^n\sum_{\gamma\in\Gamma/(\Gamma\cap
      giHg_i^{-1})}
    f(\gamma g_i v_0/T)\\
    &=\frac{1}{T^{a_\iota}(\log
      T)^{b_\iota-1}}\sum_{i=1}^n\sum_{\gamma\in g_i^{-1}\Gamma
      g_i/(g_i^{-1}\Gamma g_i\cap H)} f(g_i\gamma v_0/T).
  \end{align*}
  
  Note that $g_i\Gamma g_i\subset G(\Q)$ is an arithmetic subgroup of
  $G$, and $(g_i\Gamma g_i\inv)\cap H$ is a lattice in $H$. It follows
  from Theorem~\ref{th:sector} that the limit
  \begin{align*}
    \tau=\lim_{T\to\infty} \tau_T
  \end{align*}
  exists in the weak$^*$ topology, and $\tau$ is the $G$-invariant
  measure concentrated on $G(\mathcal{W}\cdot v_0)^{I_\iota}$ which is
  given by
  \begin{align*}
    \tau=\left(\sum_{i=1}^n\frac{\vol(H/(g_i^{-1}\Gamma g_i\cap
        H))}{\vol(G/\G)} \right)\nu_\iota,
  \end{align*}
  where $\nu_\iota$ is as in Theorem~\ref{l:asympt0}.

  Let $\phi\in C(S(W))$, $\phi\ge 0$, and let $\psi$ be the
  characteristic function of $[1/2,1)$. Take $c>1$, close to $1$, and
  $\psi^-,\psi^+\in C_c((1/4,2))$ such that
  \begin{align*}
    0\le \psi^\pm \le 1,\quad \psi^-\le \psi\le \psi^+, \quad
    \psi^-|_{[c/2,c^{-1}]}=1, \quad \supp(\psi^+)\subset [c^{-1}/2,c].
  \end{align*}
  Then for $f_{\psi} (v):=\phi(\pi(v))\psi(\norm{v})$, $\forall v\in
  W$, we have
  \begin{gather*}
    \tau_T(f_{\psi^-})\le \tau_T(f_\psi)\le \tau_T(f_{\psi^+}),\\
    \tau(f_{\psi^+})\le \tau (f_\psi(c^{-1}v))= c^{a_\iota} \tau(f_\psi),\\
    \tau(f_{\psi^-})\le \tau(f_\psi(c^{-1}v))= c^{-a_\iota}
    \tau(f_\psi).
  \end{gather*}
  Taking $c\to 1$, this implies that
  \begin{align*}
    \lim_{T\to\infty}\frac{1}{T^{a_\iota}(\log
      T)^{b_\iota-1}}\sum_{v\in V(\Z): T/2\le \norm{v}<T}
    \phi(\pi(v))=\lim_{T\to\infty} \tau_T(f_\psi)=\tau(f_\psi).
  \end{align*}
  Using that
  \begin{align*}
    \#(V(\Z)\cap B_T)\ll T^{a_\iota}(\log T)^{b_\iota-1},
  \end{align*}
  the proof can be completed by an easy geometric series argument.

  We also compute an explicit formula for the limit measure
  $\mu_{\iota}$.  Let $\phi\in C(S(W))$, $\chi$ be the characteristic
  function of $(0,1)$, and define
  \begin{align*}
    f_\chi(v)=\phi(\pi(v))\psi(\norm{v}),\quad\text{for all $v\in W$}.
  \end{align*}
  It follows from (\ref{eq:nu_mu}) that for some $c_1,c_2>0$,
  \begin{align}
    \label{eq:mu_iota}
    \int_{S(W)} \phi\, d\mu_\iota &= c_1 \int_K\sum_{w\in\mathcal{W}}
    \int_{\mathfrak{a}^{I_\iota,+}} \int_{\R}
    f_\chi\left(k\exp(a)e^t(w v_0)^{I_\iota} \right)
    \xi_{I_\iota}(a)e^{a_\iota t}\,dt dadk \\
    \nonumber &= c_2 \int_K\sum_{w\in\mathcal{W}}
    \int_{\mathfrak{a}^{I_\iota,+}} \phi\left(\pi(k\exp(a)(w
      v_0)^{I_\iota}) \right)\frac{\xi_{I_\iota}(a)}{\norm{k\exp(a)(w
        v_0)^{I_\iota}}^{a_\iota}}\,dadk.
  \end{align}
\end{proof}

\begin{proof}[Proof of Theorem~\ref{th:con}]
  Theorem~\ref{th:con} follows from Theorem~\ref{th:m} approximating
  the characteristic function of the cone by continuous functions.
\end{proof}

\begin{proof}[Proofs of Theorems~\ref{th:cone_o} and~\ref{th:m_o}]
  Proofs are based on Theorem~\ref{th:asympt} and
  Theorem~\ref{th:sector2} and are similar to the proofs of
  Theorem~\ref{th:con} and Theorem~\ref{th:m}.  We skip details.  It
  follows from Theorem~\ref{th:asympt} that the measure
  $\mu_{\Theta_\iota(\phi)}$ is given by the formula
  \begin{align}\label{eq:mu_alpha}
    &\int_{S(W)} \phi\, d\mu_{\Theta_\iota(\phi)} \\
    =&\,c \sum_{I\in\Theta_\iota(\phi)}\int_K\sum_{w\in\mathcal{W}}
    \int_{\mathfrak{a}^{I,+}} \phi\left(\pi(k\exp(a)(w v_0)^{I})
    \right)\frac{\xi_{I}(a)}{\norm{k\exp(a)(w
        v_0)^{I}}^{a_{\iota}(\phi)}}\,dadk\nonumber
  \end{align}
  where $\phi\in C(S(W))$ and $c>0$.
\end{proof}

\subsection{\it Proof of Corollary~\ref{cor:dioph}}
\label{subsec:proof:cor:dioph}
Take any $\e>0$ and Consider the cone
\begin{align*}
  \mathcal{C}=\{w\in W\setminus\{0\}:\, \norm{\pi(w)-v_0}<\e\}.
\end{align*}
Since $V^\infty$ has only finitely many orbits of $G$ the cone is
generic for sufficiently small $\e>0$.

Suppose $\del \cC$, the boundary of $\cC$, has strictly positive
measure with respect to the smooth measure class on a $G$-orbit, say
$\cO_1$, contained in $V^\infty$. Since $\cO_1$ and $\del \cC$ are
real analytic varieties, we conclude that $\cO_1\subset\del \cC$.
Since $V^\infty$ has only finitely many distinct $G$-orbits, and $\del
\cC$ are disjoint for distinct $\e>0$, we conclude that $\cC$ is
admissible for sufficiently small $\e>0$. Now the corollary follows
from Theorem~\ref{mt}.  \qed

\section{Comparison with Chambert-Loir--Tschinkel
  conjecture}\label{s:tsch}

Recently, Chambert-Loir and Tschinkel proposed a general conjecture
about asymptotics of the number of integral points on algebraic
varieties.  A weaker version of this conjecture appeared in \cite{ht}.
To facilitate a comparison, we state some of our results using the
language of arithmetic algebraic geometry.

\newcommand{\eff}{\operatorname{eff}}

Let ${\mathbf G}$ be a connected $\Q$-simple adjoint algebraic group,
which is isotropic over $\R$, and ${\mathbf X}$ the wonderful
compactification of ${\mathbf G}$. The wonderful compactification was
constructed over $\C$ in \cite{cp} and over arbitrary field of odd
characteristic in \cite{cs}. It is a smooth projective variety defined
over $\Q$ such that ${\mathbf G}$ is contained densely in ${\mathbf
  X}$, and $D:={\mathbf X}\setminus {\mathbf G}$ is a divisor with
normal crossings and smooth irreducible components.  Given a field
$k\supset\Q$, we set ${\mathbf X}_k={\mathbf X}\times_\Q k$.  Let
$\Pic({\mathbf X}_\C)$ be the absolute Picard group,
$\Lambda_{\eff}(X_\C)\subset \Pic({\mathbf X}_\C)\otimes \R$ the
effective cone, and $K_{\mathbf X}$ the canonical class.  We denote by
$\Delta_\C$ the system of simple roots of $\mathbf G$.  It was shown
in \cite{cp} (cf.~\cite[Sec.~6.1]{bk}) that there is an isomorphism
$\lambda\mapsto [L_\lambda]$ between the weight lattice of $\mathbf G$
and the Picard group $\Pic({\mathbf X}_\C)$ such that the irreducible
components of the boundary divisor $D$ correspond to $L_\alpha$,
$\alpha\in\Delta_\C$.  Note that these irreducible components generate
a finite index subgroup in $\Pic({\mathbf X}_\C)$ (the root lattice).
Given $v\in {\mathbf X}(\C)$ and $[L]=\sum_{\alpha\in\Delta_\C}
q_\alpha [L_\alpha]\in {\rm Pic}({\mathbf X}_\C)$, we set
\begin{align*}
  I(v)=\{\alpha\in\Delta_\C : v\in \supp L_\alpha\}
  \quad\text{and}\quad [L]_v=\sum_{\alpha\in I(v)} q_\alpha
  [L_\alpha].
\end{align*}
We define a metric on the real projective space:
\begin{align*}
  d([x],[y])=\frac{\norm{x\wedge y}}{\norm{x}\cdot\norm{y}},
\end{align*}
where $\norm{\cdot}$ is the standard Euclidean norm.

\begin{thm}\label{th:picard}
  Let $\mathcal{G}$ be a group scheme over $\hbox{\rm Spec}(\Z)$ with
  generic fiber ${\mathbf G}$.  Then there exists $k\in\N$ such that
  for every ample metrized line bundle $\mathcal{L}=(L,H_\mathcal{L})$
  on ${\mathbf X}$ defined over $\Q$, every $v\in ({\mathbf
    X}\setminus {\mathbf G})(\R)$, and every sufficiently small
  $\e=\e(v)>0$,
  \begin{align*}
    \#\{z\in\mathcal{G}(\frac{1}{k}\Z):\, d(z,v)<\e,
    H_\mathcal{L}(z)<T\}\sim_{T\to\infty} c\cdot T^{a}(\log T)^{b-1},
  \end{align*}
  where
  \begin{align*}
    c&=c(v,\mathcal{L},\e)>0,\\
    a&=a(v,L)=\inf\{r:\, r[L]_v+[K_{\mathbf
      X}+D]_v\in\Lambda_{\eff}({\mathbf X}_\C)_v\},\\
    b&=b(v,L)=\left\{\begin{tabular}{l} the co-dimension of the face of
        $\Lambda_{\eff}({\mathbf
          X}_\C)_v$\\
        containing $a[L]_v+[K_{\mathbf X}+D]_v$
      \end{tabular}\right\}.
  \end{align*}
\end{thm}

\begin{rem} \rm Theorem~\ref{th:picard} holds with $k=1$ if we take
  $v\in \overline{\mathcal{G}(\Z){\mathbf G}(\R)^\circ}$. In
  particular, we can take $k=1$ when ${\mathbf
    G}(\R)=\mathcal{G}(\Z){\mathbf G}(\R)^\circ$.  This equality holds
  assuming that $\mathbf G$ is $\Q$-split and $\mathcal{G}$ is the
  canonical $\Z$-model of $\mathbf G$ (see \cite[Remark in Sec.
  2]{wh}).
\end{rem}

\begin{rem} 
  \rm Our results also apply to non-smooth compactifications of
  $\mathbf G$ (for example, one can take the closure of the image of
  $\mathbf G$ under the irreducible representation with the highest
  weight $\sum_\alpha n_\alpha \omega_\alpha$ with some $n_\alpha=0$).
  We expect that an analogue of Theorem~\ref{th:picard} holds with
  parameters $(a,b)$ computed with respect to the minimal resolution
  of singularities of the pair $({\mathbf X},D)$.  A basic example of
  this type was worked out in \cite{ht}.
\end{rem}

\begin{proof}[Proof of Theorem~\ref{th:picard}]
  We refer to \cite[Sec.~6.1]{bk} for standard facts about the
  wonderful compactification.  Recall that the effective cone
  $\Lambda_{\eff}({\mathbf X}_\C)$ is generated by $[L_\alpha]$ for
  $\alpha\in \Delta_\C$, the ample cone is generated by
  $[L_{\omega_\alpha}]$ for $\alpha\in\Delta_\C$ ($\omega_\alpha$'s
  are the fundamental weights), and
  \begin{align*}
    K_{\mathbf X}\sim -L_{2\rho_\C}-\sum_{\alpha\in\Delta_\C}
    L_\alpha,
  \end{align*}
  where $2\rho_\C$ is the sum of positive roots of $\Delta_\C$.  The
  support of $L_\alpha$ is isomorphic to the fibration over ${\mathbf
    G}/{\mathbf P}_\alpha\times {\mathbf G}/{\mathbf P}_\alpha$, where
  ${\mathbf P}_\alpha$ is the maximal parabolic subgroup corresponding
  to $\alpha$, with fibers equal to the wonderful compactification of
  the adjoint form of the Levi subgroup of ${\mathbf P}_\alpha$.  This
  implies that the Galois action on ${\rm Pic}({\mathbf X}_\C)$
  correspond to the twisted Galois action ($\star$-action) on
  $\Delta_\C$.

  We denote by $\Delta$ the system of restricted roots (with respect
  to a Cartan involution) chosen so that $r(\Delta_\C)=\Delta\cup
  \{0\}$ where $r$ is the restriction map.  We have the decomposition
  into disjoint $\mathbf G$-orbits:
  \begin{align*} {\mathbf X}=\bigcup_{I\subset \Delta_\C} {\mathbf
      O}_I,
  \end{align*}
  where ${\mathbf O}_{\Delta_\C}={\mathbf G}$ and ${\mathbf
    O}_I\subset \supp L_\alpha$ iff $\alpha\notin I$. The structure of
  the set ${\mathbf X}(\R)$ was described in \cite[Ch. 7]{ji}. In
  particular, we have
  \begin{align*} {\mathbf X}(\R)=\bigcup_{I\subset \Delta} {\mathbf
      O}_{r^{-1}(I)}(\R).
  \end{align*}
  The set ${\mathbf X}(\R)$ is a union of of finitely many Satake
  compactifications $\bar V$ of ${\mathbf G}(\R)^o$ so that
  $V^\infty_I\subset {\mathbf O}_{r^{-1}(I)}(\R)$ for every $I\subset
  \Delta$. It follows from the weak approximation for $\mathbf G$ that
  each connected component of ${\mathbf G}(\R)$ contains a rational
  point.  We take $k\in \N$ so that each connected component of
  ${\mathbf G}(\R)$ contains a point from
  $\mathcal{G}(\frac{1}{k}\Z)$. By Borel--Harish-Chandra theorem,
  $\mathcal{G}(\frac{1}{k}\Z)$ is a union of finitely many
  $\mathcal{G}(\Z)$-orbits, and it suffices to compute the asymptotic
  for each of these orbits. For simplicity, we consider the orbit of
  the identity.

  For $\alpha\in \Delta$, we set
  \begin{align*}
    L_\alpha:=\sum_{\beta\in\Delta_\C} L_{r^{-1}(\alpha)}.
  \end{align*}
  Let $J=r(I(v))\subset \Delta$. Then $v\in V^\infty_{\Delta\setminus
    J}$. If $L\sim L_\lambda$ for a dominant weight $\lambda$, we get
  \begin{align*} [L]_v=\sum_{\alpha\in J}
    \frac{m_\alpha}{\#r^{-1}(\alpha)}\cdot [L_\alpha],\quad
    [K_{\mathbf X}+D]_v=\sum_{\alpha\in J}
    \frac{u_\alpha}{\#r^{-1}(\alpha)}\cdot [L_\alpha],
  \end{align*}
  where $m_\alpha$'s and $u_\alpha$'s are given as in
  (\ref{eq:defi0}).

  Passing to a tensor power, we may assume that $L\sim L_\lambda$ is
  very ample and linearized, and $\lambda$ is in the root lattice.
  Let $\iota$ be the $\Q$-representation of $\mathbf G$ on
  $H^0({\mathbf X}, L)$.  Then
  \begin{align*}
    H_\mathcal{L}(z)=H(\iota(z)),\quad z\in\mathcal{G}(\Z),
  \end{align*}
  where $H$ is the standard height function with respect to a lattice
  $\Lambda\subset H^0({\mathbf X}, L)$. Passing to a finite index
  subgroup, if necessary, we may assume that
  $\iota(\mathcal{G}(\Z))\subset \hbox{Stab}(\Lambda)$. Then
  \begin{align*}
    H_\mathcal{L}(z)=\norm{\iota(z)},\quad z\in\mathcal{G}(\Z),
  \end{align*}
  where $\norm{\cdot}$ is a norm on $H^0({\mathbf X}, L)\otimes \R$.
  Since the representation $\iota$ has the unique highest weight
  $\lambda$, the results of Section~\ref{sec:intro} apply (see
  Corollary~\ref{cor:dioph} and Example~\ref{sec:group}).
\end{proof}

\end{document}